\documentclass[12pt]{amsart}
\usepackage[T1]{fontenc}
\usepackage[utf8]{inputenc}
\usepackage[english]{babel}
\usepackage{epsfig}
\usepackage{amsmath}
\usepackage{amssymb}
\usepackage{amscd}
\usepackage{latexsym}
\usepackage{tabularx}
\usepackage{a4wide}
\usepackage[usenames]{color}
\usepackage{enumerate}
\usepackage{subfigure}

\usepackage{url}
\usepackage{verbatim}

\newcommand{\PC}{\mathscr{P}}

\usepackage{circuitikz}
\ctikzset{bipoles/resistor/height=0.15}
\ctikzset{bipoles/resistor/width=0.4}

\usepackage[
colorlinks, citecolor=blue,
pdfauthor={Omid Amini, Eduardo Esteves}, 
pdftitle={Voronoi tilings, toric arrangements and degenerations of line bundles II},
pdfstartview ={FitV},
]{hyperref}

\usepackage{tikz, float} 
\usetikzlibrary {positioning}

\usetikzlibrary{calc,decorations.markings}
\usetikzlibrary{shapes,snakes}

\numberwithin{equation}{section}

\tikzstyle{Cwhite}=[scale = .8,circle, fill = white, minimum size=3mm] 
\tikzstyle{Cgray}=[scale = .4,circle, fill = gray, minimum size=3mm] 
\tikzstyle{Cblack2}=[scale = .4,circle, fill = black, minimum size=5mm] 
\tikzstyle{Cblack}=[scale = .7,circle, fill = black, minimum size=3mm]
\tikzstyle{C0}=[scale = .9,circle, fill = black!0, inner sep = 0pt, minimum size=3mm]
\tikzstyle{C1}=[scale = .7,circle, fill = black!0, inner sep = 0pt, minimum size=3mm]
\tikzstyle{Cred}=[scale = .4,circle, fill = red, minimum size=3mm]

\usepackage[matrix,arrow,tips,curve]{xy}
\usepackage{pb-diagram,pb-xy}
\usepackage{verbatim}
\usepackage{mathrsfs}
\usepackage{color}
\usepackage[leqno]{amsmath}

\newcommand{\Vor}{\mathrm{Vor}}
\newcommand{\FP}{\mathcal{FP}}

\newcommand{\E}{\mathbb E}

\newcommand{\AC}{\mathcal{CAC}}
\newcommand{\supp}{\mathrm{supp}}

\newcommand{\m}{\mathfrak m}

\small\normalsize

\renewcommand{\:}{\colon}

\newcommand{\ol}{\overline}

\renewcommand{\Im}{{\mathrm{Im}}}

\newtheorem{thm}{Theorem}[section]

\newtheorem{lemma}[thm]{Lemma}

\newtheorem{prop}[thm]{Proposition}

\newtheorem{claim}[thm]{Claim}
\newtheorem{cor}[thm]{Corollary}

\newtheorem{remark}[thm]{Remark}

\newtheorem{defi}[thm]{Definition}

\newtheorem{question}[thm]{Question}

\newtheorem{example}[thm]{Example}

\newcommand{\dl}{\mathfrak d}

\newcommand{\U}{{\mathcal U}}

\def\X{\mathcal X}

\newcommand{\Z}{\mathbb{Z}}

\newcommand{\R}{\mathbb{R}}

\newcommand{\N}{\mathbb{N}}

\renewcommand{\k}{\kappa}

\newcommand{\Spec}{\operatorname{Spec}}
\newcommand{\g}{\mathfrak g}

\newcommand{\f}{\mathfrak f}
\newcommand{\Po}{\small{\textsc{P}}}

\usepackage{marginnote}
\newcommand{\hooklongrightarrow}{\lhook\joinrel\longrightarrow}

\numberwithin{equation}{section}

\newcommand{\he}{\text{h}}  
\newcommand{\te}{\text{t}}  

\newcommand{\comp}[1]{\overbar[.7]{#1}}
   
\newcommand{\mg}{\mathscr M} 
\newcommand{\mgbar}{\comp\mg}

\usepackage{scalerel}

\newcommand{\ind}[1]{{_{\scaleto{#1}{4.8pt}}}}

\newcommand{\indm}[1]{{_{\scaleto{#1}{3.2pt}}}} 
\newcommand{\indmbar}[1]{{_{\scaleto{#1}{4.6pt}}}} 
\newcommand{\indbi}[2]{_{{\scaleto{#1}{4.8pt}}_{\hspace{-.05cm}#2}}} 

\makeatletter
\newsavebox\myboxA
\newsavebox\myboxB
\newlength\mylenA

\newcommand*\overbar[2][0.75]{%
    \sbox{\myboxA}{$\m@th#2$}%
    \setbox\myboxB\null
    \ht\myboxB=\ht\myboxA%
    \dp\myboxB=\dp\myboxA%
    \wd\myboxB=#1\wd\myboxA
    \sbox\myboxB{$\m@th\overline{\copy\myboxB}$}
    \setlength\mylenA{\the\wd\myboxA}
    \addtolength\mylenA{-\the\wd\myboxB}%
    \ifdim\wd\myboxB<\wd\myboxA%
       \rlap{\hskip 1\mylenA\usebox\myboxB}{\usebox\myboxA}%
    \else
        \hskip -0.5\mylenA\rlap{\usebox\myboxA}{\hskip 0.5\mylenA\usebox\myboxB}%
    \fi}
\makeatother

\renewcommand{\o}{\mathfrak o}
\renewcommand{\P}{\mathbf P}
\newcommand{\Gm}{\mathbf{G_m}}
\newcommand{\A}{\mathbf A}

\newcommand{\x}{\textsc{x}}

\newcommand{\p}{\textsc{p}}
\newcommand{\q}{\textsc{q}}

\newcommand{\bt}{\mathrm{bt}}

\title{Voronoi tilings, toric arrangements and degenerations of line bundles II}
 \author{Omid Amini}
 \author{Eduardo Esteves}
 \date{\today}
  \address{CNRS - Centre de math\'ematiques Laurent Schwartz, \'Ecole Polytechnique,  France}
\email{omid.amini@polytechnique.edu}

\address{Instituto Nacional de Matem\'atica Pura e Aplicada, Estrada Dona Castorina 110,
22460-320 Rio de Janeiro RJ, Brazil}
\email{esteves@impa.br}
\begin{document}

\begin{abstract}
 We describe limits of line bundles on nodal curves
  in terms of toric arrangements associated to Voronoi tilings of Euclidean spaces. 
These tilings encode information on the relationship between the 
possibly infinitely many limits, and ultimately give rise to a new
definition of \emph{limit linear series}. This article and its first
and third part companion parts are the first in a series aimed to explore this new
approach. 

In the first part, we set up the combinatorial framework and 
showed how graphs weighted with integer lengths associated to the 
edges provide tilings of Euclidean spaces by polytopes associated to
the graph itself and to its subgraphs. 

In this part, we describe the \emph{arrangements of  toric varieties}
associated to these tilings. Roughly speaking, the normal fan to each polytope in the
tiling corresponds to a toric variety, and these toric varieties are
glued together in an arrangement according to how the polytopes
meet. We provide a thorough description of these toric arrangements
from different perspectives: by using normal fans,  as unions of torus
orbits, by describing the (infinitely many) polynomial equations
defining them in products of doubly infinite chains of projective
lines, and as degenerations of algebraic tori.

These results will be of use in the third part to achieve our goal of describing all \emph{stable limits} of
a family of line bundles along a degenerating family of curves. 
Our main result there asserts that the collection of all these limits is 
parametrized by a connected 0-dimensional closed substack of 
the Artin stack of all torsion-free rank-one sheaves on the limit
curve. The substack is the quotient of an arrangement of toric
varieties as described in the present article by the torus of the same
dimension acting on it. 
\end{abstract}
\maketitle

\tableofcontents

     \addtocontents{toc}{\protect\setcounter{tocdepth}{1}}

\section{Introduction}

This is a sequel to our work~\cite{AE1}. Our aim here is to give the geometric significance of the  combinatorial results we proved there. This will be done by using toric geometry and will be of use  in~\cite{AE3} to achieve our goal of describing stable limits of line bundles on nodal curves. 
In this introduction, after briefly recalling the results we proved in~\cite{AE1}, we explain the main contributions of the present work. A reader interested in algebraic geometry shall find a more detailed exposition
and a clearer picture of the link to the algebraic geometry of curves in the introduction to the third part.

\subsection{Stable limits of line bundles}\label{stablinebdles} 
Let $X$ be a stable curve
of arithmetic genus $g$ over an algebraically closed field $\kappa$. 
Recall that this means $X$ is a reduced, projective and 
connected scheme of dimension one that has at most ordinary nodes as
singularities, a nodal curve, and in addition has an ample canonical bundle. Stable
curves appear on the boundary of the Deligne--Mumford compactification
$\mgbar_g$ of the moduli space $\mg_g$ of smooth projective curves of genus $g$. 

\smallskip

The question which motivates the study undertaken in our series of
papers starting with~\cite{AE1}, the current paper and \cite{AE3} is the following:

\begin{question}\label{ques:main} Let $X$ be a stable curve of genus $g$ and denote by $x$ the corresponding point in the Deligne--Mumford compactification $\mgbar_g$. 

\smallskip

Describe all the possible limits of linear series $\g^r_d$ over any sequence $X_1, X_2, \dots$ of smooth projective curves of genus $g$ whose corresponding points $x_1, x_2, \dots$ in $\mg_g$ converge to $x$.
\end{question}

Recall that a linear series $\g^r_d$ on a projective curve $Y$ is by definition a vector subspace of dimension $r+1$ of the space of global sections $H^0(Y, L)$ of a line bundle $L$ of degree $d$ on $Y$. 

\smallskip

The above simple looking question turns out to be a difficult and
multifaceted problem for which only partial results are known. In the
case where $X$ is a curve of compact type, a satisfactory answer was
given in the pioneering work by Eisenbud and Harris~\cite{EH86}, who
derived several interesting consequences  in the study of the geometry
of the moduli space of curves~\cite{EH87-WP,EH87-KD}. Apart from the
case of compact type curves, Medeiros and the second named
author~\cite{EM} gave an answer to the problem for curves with two
components, and, recently, Osserman extended the above mentioned
results to the case of curves of pseudo-compact
type~\cite{Os-pseudocompact, Oss16}, which cover both the cases
considered previously in~\cite{EH86, EM}. As far as the stable curves
lying on the deeper strata of the moduli space are concerned, only
little is known. The most interesting results in this direction are
probably the ones obtained in the work by Bainbridge, Chen, Gendron,
Grushevsky and M\"oller, \cite{BCGGM18}, which describes 
limits of holomorphic one forms, i.e., $\g^0_{2g-2}$ obtained as
 global sections of the
canonical sheaf, and its sequel~\cite{BCGGM19}, which studies sections of powers 
of the canonical sheaf. The above question is also 
intimately linked to tropical and hybrid geometry; we refer 
to~\cite{AB15, MUW, LM18, BJ16, Cart15, Oss19, He19, FJP} for other
results which go in the direction of partially answering the question
we posed. 

\smallskip

The main combinatorial object associated to a stable curve $X$ is its dual
graph and it is now well-understood that dual graphs and their
geometry play a central role in understanding questions of the type
described above. The dual graph $G=(V,E)$ of $X$ 
consists of a vertex set $V$ in one-to-one
correspondence with the set of irreducible components of $X$, and an edge set $E$ in
one-to-one correspondence with the set of nodes. An edge lies between
two vertices if the corresponding node lies on the two corresponding 
components. 

\smallskip

In our series of work we explore a new approach to
Question~\ref{ques:main} based on the geometry of dual graphs and the
use of  Artin stacks and their ``Hilbert schemes.'' In~\cite{AE3}, we
give an answer to the following question:

 \begin{question}\label{ques:main2} Let $X$ be a stable curve of genus $g$ and denote by $x$ the corresponding point in the Deligne--Mumford compactification $\mgbar_g$. 

\smallskip

Describe all the stable limits of line bundles over any sequence 
$X_1, X_2, \dots$ of smooth projective curves of genus $g$ whose 
corresponding points $x_1, x_2, \dots$ in $\mg_g$ converge to $x$.
\end{question}

More precisely, let $\pi\colon\mathcal X\to B$ be a one-parameter
smoothing of $X$. Here, $B$ is the spectrum of $\kappa[[t]]$, and 
$\pi$ is a projective flat morphism whose generic fiber $\mathcal X_\eta$ is
smooth and whose special fiber is isomorphic to $X$. We fix such an
isomorphism.  The total space $\mathcal X$ is regular except possibly at the nodes of
$X$. For each edge $e\in E$, one can associate to the corresponding
node $N_e$ of $X$ a positive integer $\ell_e$, called the the
singularity degree (or the thickness) of $\pi$ at $N_e$. This is
obtained by looking at 
the completion of the local ring of $\mathcal X$ at $N_e$, which is isomorphic to
$\k[[u,v,t]]/(uv-t^{\ell_e})$. If $\ell_e=1$ for an edge $e\in E$, 
then $\mathcal X$ is regular at $N_e$.

Let $L_\eta$ be an invertible sheaf on the generic fiber. If $\mathcal
X$ is
regular, it extends to an invertible sheaf $\mathcal L$ on $\X$. It is
not unique, as $\mathcal L\otimes\mathcal O_{\mathcal X}(\sum
f(v)X_v)$ is another extension, for every integer valued function
$f\colon V \to 
\mathbb Z$. (Here $X_v$ is the component corresponding to $v\in V$.) For
a general family $\pi$, without the regularity assumption, the sheaf $L_\eta$ extends to a relatively
torsion-free, rank-one sheaf $\mathcal I$ on $\X/B$, that is, a $B$-flat 
coherent sheaf on $\mathcal X$ whose fibers over $B$ are torsion-free,
rank-one. Again, there is no uniqueness for such an extension. Furthermore, one could perform a finite base change
to the family $\pi$, take the pullback of $L_\eta$ to the new generic fiber and consider its
torsion-free extensions. These will be extensions on a different total
space, but the special fibers of all the families are the same, and thus the restrictions
of all these extensions to $X$ are torsion-free, rank-one sheaves that
we call the \emph{stable limits} of $L_\eta$.

\smallskip

We describe these stable limits of line bundles on stable curves in terms of Voronoi decompositions of Euclidean spaces associated to graphs and their associated arrangements of toric varieties. 

\smallskip

\subsection{Voronoi tilings} Let $X$ be a nodal curve and let $G=(V, E)$  be its dual
graph, with self loops removed. We denote by $\mathbb E$ the set of all the possible orientations of the edges in $G$. Thus, for each
edge there are two possible arrows, pointing to the two different
vertices joined by that edge. For $e\in\mathbb E$ with end-points $u$
and $v$, we write $e=uv$ if $e$ has tail equal to $u$ and head equal
to $v$, and let $\ol e$ denote the reverse orientation of the same edge. 

For a coefficient ring $A$, we consider the following two complexes:
$$
d_A\colon C^0(G,A) \to C^1(G,A)\quad\text{and}\quad
\partial_A\colon C_1(G,A) \to C_0(G,A).
$$
Here, $C^0(G,A)$ denotes the $A$-module of functions $V\to A$ and 
$C_0(G,A)$ the free $A$-module generated by the vertices in $V$. Similarly, $C^1(G,A)$
is the $A$-module of functions $f\colon\mathbb E\to A$ which verify the 
condition  $f(\ol e)=-f(e)$ for each $e\in\mathbb E$,
and $C_1(G,A)$ is the quotient of the free $A$-module generated by 
$\mathbb E$ modulo the submodule generated by $e+\ol e$ for all 
$e\in\mathbb E$. Finally, the differentials are given by
$d(f)(e)=f(v)-f(u)$ and 
$\partial(e)=v-u$ for
each oriented edge $e=uv\in\mathbb E$. 

The above $A$-modules come with natural isomorphisms $C_0(G,A)\to C^0(G,A)$, which takes $v$ to
the characteristic function $\chi_{\indm v}$, and $C_1(G,A)\to C^1(G,A)$, which takes $e$ to $\chi_{\indm e}-\chi_{\indmbar{\ol e}}$. Also, we have  bilinear forms $\langle\,,\rangle$ on $C_0(G,A)$ and $C_1(G,A)$ such that 
$\langle v,w\rangle=\delta_{v,w}$ for $v,w\in V$ and
$\langle e,f\rangle =\delta_{e,f}-\delta_{e,\overline f}$ for $e,f\in\mathbb E$. The above
isomorphisms induce bilinear forms on $C^0(G,A)$ and 
$C^1(G,A)$ as well. Denote by $d^*_A\colon C^1(G,A)\to C^0(G,A)$ 
the homomorphism corresponding to $\partial_A$ under the
isomorphisms. Then, $d^*_A$ is the adjoint to $d_A$, that is,
$\langle f,d^*_A(h)\rangle\,=\,\langle d_A(f),h\rangle $ for each $f\in C^0(G,A)$ and $h\in C^1(G,A)$.

\smallskip

Let $H_{0,A}:=\{f\in C^0(G,A)\,|\,\sum f(v)=0\}$ and $F_A:=\Im(d_A)$. 
The Laplacian of $G$ is $\Delta_A:=d^*_Ad_A$. The homomorphism 
$d^*_\R$ induces an isomorphism $F_\R\to H_{0,\R}$. Also, 
the bilinear form $\langle\,,\rangle$ induces a norm on $F_\R$
corresponding via $d^*_{\mathbb R}$ to the quadratic form $q$
on $H_{0,\mathbb R}$ satisfying $q(f)=\langle f,\Delta_{\mathbb R}(f)\rangle$ for
each $f\in C^0(G,\mathbb R)$. 

In~\cite{AE1}, we described certain tilings of
$H_{0,\mathbb R}$ by polytopes.  For instance, let $\Lambda_A:=\text{Im}(d^*_A)$. Then 
$\Lambda_{\mathbb R}=H_{0,\mathbb R}$ and $\Lambda_{\mathbb Z}$ is a
sublattice of $H_{0,\mathbb Z}$ of finite index equal to the number of
spanning trees of $G$, by the Kirchhoff Matrix-tree Theorem. The
\emph{standard Voronoi tiling} of $G$, denoted $\Vor_G$, is the Voronoi decomposition of
$H_{0,\mathbb R}$ with respect to $\Lambda_{\mathbb Z}$ and $q$: the
tiles are 
$$
\text{Vor}_q(\beta):=\{\eta\in H_{0,\mathbb R}\,|\,
q(\eta-\beta)\leq q(\eta-\alpha)\text{ for every }\alpha\in 
\Lambda_{\mathbb Z}-\{\beta\}\}
$$
for $\beta\in \Lambda_{\mathbb Z}$. This is one of the infinitely many different tilings we
consider. The other ones are variations of this one that we call \emph{twisted mixed
Voronoi tilings} and denote $\text{Vor}^{\m}_{G,\ell}$. 
Though the standard Voronoi tiling is homogeneous,
meaning all tiles are translates of the tile centered at the origin,
$\text{Vor}_G(O)$, a twisted mixed Voronoi tiling is obtained by
putting together translations of the tiles $\text{Vor}_H(O)$
associated to a collection of connected spanning subgraphs $H$ of
$G$.  More precisely, the twisted mixed Voronoi
tiling $\text{Vor}^{\m}_{G,\ell}$ depends on $\m\in C^1(G,\Z)$ (the
``twisting'') and an edge length
function $\ell\: E\to\N$; the tiles are the polytopes 
$d^*(\dl^\m_f)+\text{Vor}_{G^\m_f}(O)$ for
$f\in C^0(G,\Z)$ with $G^\m_f$ connected, 
where $\dl^\m_f\in C^1(G,\R)$ is a modification of $d_\Z(f)$ defined by
$$
\dl^\m_f(e):=\frac{1}{2}\big(\delta^{\m,\ell}_e(f)-\delta^{\m,\ell}_{\bar
  e}(f)\big),
\text{ where }
\delta^{\m,\ell}_e(f):=\Big\lfloor\frac{f(v)-f(u)+\m_e}{\ell_e}\Big\rfloor
\text{ for each $e=uv\in\E$,}
$$
and $G^\m_f$ is the spanning subgraph of $G$ obtaining by removing the edges
$e\in\E$ for which $\dl^\m_f(e)\not\in\Z$. The construction of these tilings is briefly recalled in Section~\ref{sec:recap}. We refer to \cite{AE1}
for a thorough presentation.  

\subsection{Toric arrangements}
 
 Our aim in this paper is to associate arrangements of toric varieties
 to the twisted mixed Voronoi tilings we constructed in~\cite{AE1}. 
  
To an arrangement of rational polytopes in a real vector space, 
rational with respect to a fixed full rank lattice, we describe how 
to associate an arrangement of toric varieties, a scheme whose
irreducible components are the toric varieties associated to the
polytopes, which intersect in toric subvarieties corresponding to
the faces of those polytopes. The arrangement is called a \emph{toric
  tiling} if it arises from a tiling of the whole vector space; 
see Subsection ~\ref{sec:toricarrangement}.

Let $G =(V, E)$ be a finite connected graph without loops. 
The lattice $C^1(G,\Z)$ in $C^1(G, \R)$ gives rise to an 
arrangement $\square_G$ of hypercubes hypercubes 
$\square_ \alpha$ in $C^1(G, \R)$ for $\alpha\in C^1(G, \Z)$, 
where for each $\alpha \in C^1(G, \Z)$, the points of 
$\square_\alpha$ are  of the form $\alpha + \epsilon$ with 
$|\epsilon(e)|\leq \frac 12$ for all $e\in \E$. The arrangement of 
hypercubes $\square_G$ leads to an arrangement $\mathbf R$ 
of toric varieties $\P_\alpha$ for $\alpha \in C^1(G, \Z)$, where each 
$\P_\alpha$ is isomorphic to the product $ \prod_{e\in E}
\P^1_e$; see Subsection~\ref{tortilR}. 
(Here $\P^1_e$ denotes a copy of the projective line $\P^1$.)

Consider a length function $\ell\: E \rightarrow \N$ and a 
twisting $\m\in C^1(G, \Z)$. Let $H$ be the subdivision of $G$ 
given by the length function, and consider the Voronoi tiling 
$\Vor^\m_{G,\ell}$ of $H_{0,\R}$. The tiling gives rise to an
arrangement of toric varieties that we denote by 
$Y^{\bt}_{\ell,\m}$ and call it the \emph{basic toric tiling} associated to the 
triple consisting of the graph $G$, the edge length function $\ell$ 
and the twisting $\m$; see Subsection~\ref{bastortil}.

We describe a natural embedding of $Y^{\bt}_{\ell,\m}$ in $\mathbf
R$. We do this by describing explicitly a natural embedding in $\P_0$ 
of the toric variety associated
to $\Vor_G(O)$ in Theorem~\ref{thm:toric}, the first main result of
our paper. We use this to associate to each tile of $\Vor^\m_{G,\ell}$
a toric subvariety of a certain $\P_\alpha$. We show that the union in
$\mathbf R$ of these subvarieties has the structure of
$Y^{\bt}_{\ell,\m}$, our second main result; see Theorem~\ref{vorvor} and Corollary~\ref{Ystr}.

We introduce moduli to $Y^{\bt}_{\ell,\m}$, by deforming the equations
of its irreducible components: To the extra data of characters
$$
a\: C^1(G, \Z) \to \Gm(\k)\quad\text{ and }\quad b\: H^1(G, \mathbb Z)
\to \Gm(\k),
$$
we associate  a toric tiling $Y_{\ell, \m}^{ a, b}$ in 
$\mathbf R$, that lies in 
the rational equivalent class of $Y^{\bt}_{\ell,\m}$, and coincides
with it when the characters are trivial; see
Subsection~\ref{gentortil}.  (Here $\Gm$ denotes the multiplicative
group.) 

\smallskip

The extra data is necessary to describe the stable limits on the
stable curve $X$ of an
invertible sheaf $L_\eta$ in the situation of
Subsection~\ref{stablinebdles}. We show in \cite{AE3} that the
collection of stable limits is parameterized by a connected closed
substack of dimension 0 of the Artin stack of all torsion-free,
rank-one sheaves on $X$, and describe the substack as the quotient of
$Y_{\ell, \m}^{ a, b}$ for certain $\ell$, $\m$, $a$ and $b$ by the
natural action of a certain torus. The function $\ell$ is the one
introduced in Subsection~\ref{gentortil}. The character $a$ encodes partially
the infinitesimal data of the arc drawn by $\pi$ on the moduli of
curves. The character $b$ is an
amalgamation of glueing data arising from certain limits of $L_\eta$, whereas
the twisting $\m$ keeps track of where those limits fail to be
invertible. Furthermore, we show that all $Y_{\ell, \m}^{ a, b}$ arise
in this way.

\smallskip

The rest of the current paper is then devoted to a thorough study of
the arrangements $Y_{\ell, \m}^{ a, b}$ and the way they are embedded 
in $\mathbf R$. Here is what we do:

\begin{itemize}
\item We give a description of the toric stratification of
  $Y^{a,b}_{\ell, \m}$. More precisely, we define a natural action of the
  torus $\mathbf{G}_\m^{V}/\Gm \simeq \mathbf{G}_\m^{|V|-1}$ on each
  $Y^{a,b}_{\ell, \m}$ and characterize its orbits; see
  Subsection~\ref{action} and Theorem~\ref{unionorbits}. This is the
  crucial result used in \cite{AE3}.

\item We describe a complete set of (infinitely many) equations for
  the embedding of $Y^{a,b}_{\ell, \m}$ inside $\mathbf R$; see
  Theorem~\ref{thm:main2}. Of course, this will be crucial in
  describing the moduli of the various $Y^{a,b}_{\ell, \m}$ in the
  ``Hilbert scheme'' of $\mathbf R$.

\item We show that each $Y^{a, b}_{\ell, \m}$ can be obtained as an
  equivariant flat degeneration of the algebraic torus
  $\mathbf{G}_\m^{|V|-1}$, and thus the quotient stack is the
  degeneration of a point; see Theorem~\ref{thm:degtori}. This is
  naturally a consequence of what we do in \cite{AE3}, but it is shown
  independently here.
\end{itemize}

\subsection{Organization} The layout of the paper is as follows. 
In Section~\ref{sec:recap}, we recall the Voronoi tilings introduced in~\cite{AE1}. 
In Section \ref{sec:toricarrangements}, we describe 
the toric variety associated to a given graph, and use this in 
Section \ref{mixtortil} to define the arrangements of toric varieties 
associated to the Voronoi tilings of Section~\ref{sec:recap}, 
and to describe the way the toric varieties associated to the graph
and certain of its connected subgraphs are glued together. 
In Section~\ref{mixorb}, we describe the action of the algebraic 
torus ${\mathbf G}_{\mathbf m}^{|V|-1}$ on the arrangements of toric
varieties and give a complete description of its orbits. 
In Section~\ref{sec:eqtiling}, we work out the equations giving our
arrangements for the natural embedding in $\mathbf R$. 
Finally, in Section~\ref{sec:degtori}, we show that our toric
arrangements are all obtained as equivariant degenerations of the 
algebraic torus ${\mathbf G}_{\mathbf m}^{|V|-1}$.
\medskip

  \section{Voronoi tilings associated to graphs} \label{sec:recap}

The aim of this preliminary section is to recall the construction of
the different types of Voronoi tilings associated to graphs with edge
lengths and twisting 
carried out in detail in the first part of this work~\cite{AE1}.  
In this paper we will freely use the definitions, notations and
results of~\cite{AE1}.  
For the convenience of the
reader, we reproduce a few of them here and later refer to 
appropriate parts of~\cite{AE1} for the others.

\medskip

Let $G=(V,E)$ be a finite connected graph without loops.  We denote by
$\mathbb E$ the set of all the orientations of the edges in $G$. We
use $\{u,v\}$ to denote a non-oriented edge connecting vertices $u$ and
$v$, even if not unique, and $uv$ or $vu$ for the two possible
orientations of that edge. 
For an oriented edge $e=uv\in \E$, we denote by $\ol e=vu$ the edge
$\{u,v\}$ with the reverse  orientation. We call $u$ and $v$ the tail and head of the edge $e=uv$, respectively, and sometimes denote them by $\te_e$ and $\he_e$.   In particular, we have $\te_{\ol e} = \he_{e}$ and $\he_{\ol e} = \te_{e}$.
For a subset $X \subseteq V$, we denote by $\E(X, V- X)$ the set of
oriented edges of the form $uv$ 
with $u\in X$ and $v\in  V- X$. These sets are called (oriented) cuts in the graph.  

Let  $A$ be one of the rings $\R$, $\mathbb Q$ or $\Z$. We use the
 standard notations $C^0(G, A)$ and $C^1(G, A)$ for the $A$-modules
 of 0- and 1-cochains with values in $A$, respectively: the first 
consists of all the functions $f: V \rightarrow A$, whereas the second
consists of all the functions $g: \mathbb E \rightarrow A$ which
verify 
$g(e) = -\, g(\bar
e)$ for each oriented edge $e$ in $\E$. Similarly, we use $C_0(G,A)$
and $C_1(G, A)$ for the 
$A$-modules of 0- and 1-chains on $G$: the first is the free module
generated by the vertices of $G$, the generator associated to $v\in V$
denoted by $(v)$, 
whereas the second is the quotient of the free module generated by the
 oriented edges of $G$ by the submodule generated by the relations 
$(e)=-(\ol e)$ for each $e\in\E$, where $(e)$ denotes the element
associated to $e$ in the free module and the quotient.

The cochain complex is the complex 
$C^{\bullet}(G,A): C^0(G,A)\stackrel {d}{\longrightarrow} C^{1}(G,A)$, where 
the differential $d$ is defined by
\[
d(f)(e) := d_{uv}f:=
f(v)-f(u)
\]
for each $f\in C^0(G,A)$ and $e=uv\in\E$. On the other hand,
the chain complex is the
complex
$C_{\bullet}(G,A): C_1(G,A)\stackrel{\partial}{\longrightarrow} C_0(G,A)$,
where the boundary map $\partial$ is defined by
\[
\partial (e) := (v)-(u)
\]
for each $e=uv\in\E$.

As defined above, it is clear that the spaces $C_i(G,A)$ and
$C^i(G,A)$ are canonically dual for
$i=0,1$. Moreover,  there are natural scalar products
$\langle\,,\rangle$ on
$C_0(G,A)$ and $C_1(G,A)$ that naturally identify the space $C_i(G,A)$ with
$C^i(G,A)$ for $i=0,1$. Because of this, we simply use $\alpha_e$ for the value of $\alpha\in
 C^1(G,A)$ at an oriented edge $e\in \E$.  Also, $\partial$ becomes identified with the adjoint $d^*$
of $d$, and we get $\partial d = d^* d = \Delta$, 
 the \emph{Laplacian} of the finite graph $G$, which is defined by 
\[
\Delta(f)(v) = \sum_{\substack{e\in\mathbb E\,\\ \he_e=v}}
f(\he_e)-f(\te_{e})
\]
for each $f\in C^0(G,A)$ and $v\in V$.

At some point in the text it will be more convenient to fix an
orientation of the edges. By an orientation we mean a map $\mathfrak
o: E \to \E$ which for each edge $\{u,v\}$ of $G$  takes one of the
oriented edges $uv$ or $vu$ as value. In other words, $\mathfrak o$ is
a right inverse to the natural forgetful map $\E\to E$. We denote by $E^{\mathfrak o}$ the
image of $\mathfrak o$.  Given $e\in\E$, by an abuse of the notation, we will also denote by $e$
its image in $E$. And given $e\in E$, we let $e^{\mathfrak o}$ denote $\mathfrak o(e)$.

Using an orientation $\mathfrak o$, we can identify the inclusion of
the lattice $C^1(G, \Z)$ in $C^1(G, \R)$ with the standard lattice
$\Z^E$ in $\R^E$; in this way, 
the scalar product on $C^1(G, \R)$ becomes the
standard Euclidean product on $\R^E$. We denote the corresponding norm
by $\|\,\,\|$. Furthermore, 
under this identification, the Voronoi decomposition of 
$C^1(G, \Z) \subset C^1(G, \R)$ with
respect to  $\|\,\,\|$ is the standard tiling of
$\R^E$ by hypercubes $\square_{\alpha}$, one for each element 
$\alpha \in \Z^E$,  defined by
 \[
\square_{\alpha}:=\Bigl\{\,x\in \R^E\,\Bigl|\,
 |x_e-\alpha_e| \leq \frac
 12\, \textrm{ for all }e \in E \Bigr\}.
\]
 
The submodule $F_A$ of $C^{1}(G,A)$ defined by the image of $d$ is called \emph{the module of  $A$-valued cuts}.  By definition, it is generated by elements of the form $d(\chi_{\ind{X}})$ for subsets $X \subseteq V$, where $\chi_{\ind X}$ is the characteristic function of $X$. One easily verifies that 
\begin{align*}
d(\chi_{\ind{X}})(e) :=
\left\{
\begin{array}{rl}
-1  & \quad \mbox{if } e\in \E(X,V-X),\\
1 &   \quad \mbox{if } \ol e \in \E(X , V-X),\\
0 & \quad \mbox{otherwise.}
\end{array}
\right.
\end{align*}

We call $F_\Z$ \emph{the cut lattice}  and $F_\R$ \emph{the cut space}
of $G$, respectively. Elements of the form $d(\chi_{\ind{X}})$ are
called \emph{cut elements}.

Define $\Lambda_A$ as the image by $d^*$ of $F_A$,  and let
$H_{0,A}$ be defined by 
\begin{align*}
H_{0,A} :=& \bigl\{f\in C^0(G, A)\,\Bigl|\, \sum_{v\in
  V}f(v)=0\bigr\}.
\end{align*}
It is easy to see that $\Lambda_A \subseteq H_{0,A}$, with equality if
$A=\mathbb R$.  The sublattice $\Lambda_\Z \subseteq H_{0,\Z}$
is called \emph{the Laplacian lattice} of the graph. 

We denote by $\Vor_{\|.\|}(F_\Z)$ the Voronoi decomposition of $(F_\R,
 \|.\|)$ with respect to the cut lattice $F_\Z$. The map $d^*$ induces
 an isomorphic Voronoi decomposition of 
 $\Lambda_\R$ with respect to the sublattice $\Lambda_\Z$, that we
 denote by $\Vor_G(\Lambda_\Z)$, or simply $\Vor_G$. In fact,
 the cells of $\Vor_G$ are the projections under $d^*$ of the hypercubes
 $\square_\alpha$ for $\alpha \in F_\Z$; see \cite[Thm.~3.14]{AE1}.
 
 In~\cite{AE1} we gave a detailed description of the Voronoi decomposition of 
$F_\R$ with respect to $F_\Z$. Let us briefly recall it here. 
  
Denote by $\Vor_F(\beta)$ the Voronoi cell of a 
lattice element $\beta \in F_\Z$. Let $\FP$ be the face poset of the
polytope $\Vor_F(O)$ of the origin in $F_\Z$; as a set, $\FP$ consists
of  the faces in $\Vor_F(O)$, and the inclusion between faces defines
the partial order. 
We explain now the combinatorics of $\FP$.
   
First, by a spanning subgraph $G'$ of $G$ we mean a subgraph with
$V(G')=V(G)$.   
A \emph{cut subgraph} of $G$ is by definition a spanning subgraph $G'$ of $G$ 
for which we can find a partition $V_1,\dots, V_s$ of $V$ such that
the edges of $G'$ 
are all those edges in $G$ that connect a vertex of $V_i$ to a vertex
of $V_j$ for $j\neq i$.  
Recall as well that an acyclic orientation of a graph $G$ is an
orientation which does not 
contain any oriented cycle. In~\cite{AE1} we called an orientation $D$ of a subgraph $G'$ of $G$ 
\emph{coherent acyclic} if $G'$ is a cut subgraph of $G$ given by a
partition $V_1,\dots, V_s$ of $V(G)$ that is ordered consistently
with $D$, that is, all the edges between $V_i$ and $V_j$ for $i<j$ get
orientation in $D$ from $V_i$ to $V_j$.

A coherent acyclic orientation $D$ can be viewed as the subset 
$\E(D)\subseteq\E$ of its oriented edges. The corresponding set of
non-oriented edges is denoted $E(D)$. We define $\AC$ as the set
consisting of all 
the coherent acyclic 
orientations of cut subgraphs of $G$. The set $\AC$ has a natural
poset structure: for $D_1$ and $D_2$ in $\AC$, set $D_1\preceq D_2$ 
  if $\E(D_2)\subseteq\E(D_1)$. We proved the following theorem in~\cite{AE1}.

\begin{thm}[{\cite[Thm.~3.19]{AE1}}]\label{thm:vor} The two posets $\FP$ and $\AC$ are isomorphic. 
 \end{thm}

We briefly recall how the isomorphism is defined. First, for an
element $x \in F_\R$, recall that the positive support of $x$, 
denoted $\supp^+(x)$, is defined as
 the set of all oriented edges $e\in \mathbb E$ with $x_e >0$. A \emph{bond
 element} in $F_\Z$ is any element of the form $d(\chi_{\ind X})$ for
 a nonempty $X\subsetneqq V$ such that
 both the graphs $G[X]$ and $G[V-X]$ induced on $X$ and 
$V-X$, respectively,  are connected. (The graph $G[X]$ induced by $G$ on $X$ has vertex set $X$ and edge set all those edges of $G$ with both endpoints in $X$.) 

We showed in~\cite{AE1} that bond elements form a system of generators
for $F_\Z$. 
Moreover, they define a hyperplane arrangement in $F_\R$ as follows.  Let $\beta$ be a bond element of $F_\Z$. The affine hyperplane $F_\beta$  of $F_\R$ is defined by
$$
F_\beta:=\bigl\{\,x\in F_\R\,|\, 2\langle x,\beta  \rangle =
\|\beta\|^2\,\bigr\}.
$$
Moreover, for each lattice point $\mu \in F_\Z$ and bond element
$\beta$, we consider the affine hyperplane $F_{\mu,\beta} :=
\mu+F_\beta$. The cell $\Vor_F(\mu)$ is the open cell containing $\mu$
of the hyperplane arrangement given by the $F_{\mu,\beta}$ for all
bond elements $\beta$.

Consider now the Voronoi cell  $\Vor_F(O)$ and let $\frak f$ be a face
of this polytope. 
We define  $\mathcal U_{\frak f}$ as the set of all the bond elements $\beta$ in $F_\Z$ such that $\frak f \subset F_\beta$.

With these definitions, we proved that the isomorphism 
$\phi: \FP \rightarrow \AC$ is simply given by 
$$
\phi(\frak f )
:= \bigcup_{\beta \in
  \mathcal U_{\frak f}} \supp^+(\beta).
$$

 The tiling given by  $\Vor_{\|.\|}(F_\Z)$ or $\Vor_G(\Lambda_\Z)$
is regular in the sense that all
 the Voronoi cells are translations of each other. We recall now the
 generalization of the above 
picture to the case of graphs equipped with an integer length function 
$\ell: E \to \mathbb N$ and a twisting $\mathfrak m\in C^1(G, \Z)$.
 
 Let $G=(V, E)$ be a connected loopless finite 
graph with integer edge lengths $\ell: E \to \mathbb N$, and let $H$ be the subdivision of $G$ where each edge $e$ is subdivided $\ell_e -1$ times. 
Let $\mathfrak m\in C^1(G, \Z)$. We call $\mathfrak m$ a \emph{twisting}.

There is a natural subcomplex $\square_{H}^{\mathfrak m}$ of the
standard Voronoi tiling of $C^1(G, \R)$ by hypercubes, which is
defined as follows. For an element $f\in C^0(G, \Z)$, and for an edge 
$e=uv$ in $\E$, set $\delta^\m_e(f) := \lfloor
\frac{f(v)-f(u)+\m_e}{\ell_e}\rfloor$, where $e$ is viewed in $E$ via
the natural forgetful map $\E\to E$. 
 
We reproduce the following definitions from~\cite[Section 5]{AE1}.
For each $f\in C^0(G, \Z)$, let:
\begin{itemize}
\item $\dl^\m_f\in C^1(G, \R)$ given by
$$
\dl^\m_f(e) :=   \begin{cases} 
\delta^\m_e(f) &  \textrm{ if } \ell_e \,\bigl|\, f(v) - f(u)+\m_e,\\
\delta^\m_e(f)+\frac 12&  \textrm{ otherwise}                                                
\end{cases}
$$
for each oriented edge $e=uv\in \E$;
\item $G^\m_f$ be the spanning subgraph of $G$ containing the support
  of all oriented edges $e\in\E$ for which $\dl^\m_{f}(e) \in \Z$;
\item $\square_{\dl^\m_f}$ be the hypercube of dimension $|E(G^{\m}_f)|$ defined by 
\[
\square_{\dl^\m_f}  :=  \dl^\m_f + \bigl\{ x \in C^1(G, \R) \,\bigl|
\,\,   |x_e | \leq \frac 12 \textrm{ for all $e$, and $x_e = 0$ for all }
e\in \E(G)-\E(G^{\m}_f)
\bigr\};
\]
\item $\square_{G, \ell}^\m$ be the arrangement of hypercubes $ \square_{\dl^\m_f}$ in $C^1(G,\R)$, so 
$$
\square^\m_{G, \ell} := \bigcup_{f\in C^0(G, \Z)}\,\square_{\dl^\m_f}.
$$
\end{itemize}

 We use $\square_f^\m$ instead of $\square_{\dl^\m_f}$ and
$\square_H^\m$ instead of  $\square_{G, \ell}^\m$, if there is no risk
of confusion.

 We consider
now the projection map $d^*: \square_H^\m \rightarrow
H_{0,\R}$.   For each 
$f \in C^0(G, \R)$ with connected subgraph $G^\m_f$, define
$\Vor_{H}^\m(f)$ as the image 
$d^*(\square^{\m}_f)$. We showed in~\cite[Prop.~5.7]{AE1} that 
$$
\Vor_{H}^\m(f)= d^*(\dl^\m_f) + \Vor_{G^\m_f}(O),
$$
where
$\Vor_{G^\m_f}(O)$ denotes the Voronoi cell of the origin in
$H_{0,\R}$ for the graph $G^\m_f$. Furthermore, we proved the following theorem:

 \begin{thm}[{\cite[Thm.~5.9]{AE1}}]\label{thm:projection3} The set
   of polytopes $\Vor^\m_H(f)$ for $f\in C^0(G, \Z)$ with connected
   spanning subgraph $G^\m_f$ is a tiling  of $H_{0,\R}$. Each $\Vor^\m_H(f)$
   in this tiling is congruent to the Voronoi cell of $G_f^\m$.
 \end{thm}

We denote this tiling by $\Vor^\m_{G,\ell}$ or simply $\Vor^\m_H$ if
there is no risk of confusion. If $\m=0$ we drop the upper index $\m$, and if $\ell=1$ we
drop the lower index $\ell$.

  \section{Toric tilings associated to graphs}\label{sec:toricarrangements}

The aim of this section is to introduce tilings by toric varieties associated to the tilings by polytopes of the previous section.  

\subsection{Toric tilings} \label{sec:toricarrangement}  Let $N$ be a
lattice and $M  := N ^\vee := \mathrm{Hom}(N, \Z)$ be its dual.  Let
$\langle\, ,\rangle$ denote the natural pairing of $M_\R$ and $N_\R$.
A rational polyhedral complex $\PC$ in $M_\R$ is by definition a (possibly infinite) collection of rational polyhedra in $M_\R$ verifying the following two properties:
\begin{itemize}
\item  A face  $\f$ of a polyhedron $\Po \in \PC$ belongs to $\PC$.
\item For any pair of polyhedra $\Po_1$ and $\Po_2 $ in $\PC$, their
  intersection is either empty or a common face to both $\Po_1$ and $\Po_2$. 
\end{itemize}

We consider only polyhedral complexes which are locally finite, that is,
for which each point $x$ on $M_\R$ lies in a neighborhood which
intersects only finitely many members of $\PC$. 

\begin{example}\label{extorarr}\rm Our two
  natural examples are the tiling of $C^1(G,\R)$ by hypercubes with
  lattice $C^1(G,\Z)$, and
  the polyhedral complex
  given by the arrangement of polytopes $\Vor^\m_H$ in $H_{0,\R}$ with
  lattice $M=\Lambda_\Z$, or the isomorphic arrangement in $F_\R$
  with lattice $M = F_\Z$ induced by the isomorphism $d^*\:
  F_\Z\to\Lambda_\Z$. 
\end{example}

Let $\PC$ be a rational polyhedral complex and let $\Po$ be a polyhedron in
$\PC$. Denote by $M_{\Po,\R}$ the linear subspace of $M_\R$ parallel to
$\Po$, thus $M_{\Po, \R}$ is the tangent space to any point in the
relative interior of $\Po$. Define the lattice $M_{\Po} :=  M \cap
M_{\Po,\R}$, and note that by the rationality of $\PC$, the lattice $M_{\Po}$ is of full
rank in $M_{\Po, \R}$. 

For each subset $\sigma\subseteq M_\R$ let
$$
\sigma^\vee:=\{v\in N_\R\,|\,\langle v,w\rangle \geq 0\text{ for every
}w\in\sigma\}.
$$
Analogously, define $\sigma^\vee\subseteq M_\R$ for each
$\sigma\subseteq N_\R$. Notice that $\sigma^\vee$ is a convex cone. If
$\sigma$ is a linear subspace, so is $\sigma^\vee$. Set $N_{\Po,\R}:=M_{\Po,\R}^\vee$ for each polyhedron
$\Po$ in $\PC$, and put $N_{\Po}:=N\cap N_{\Po,\R}$. Again by the
rationality of $\PC$, the lattice $N_{\Po}$ is of full rank in
$N_{\Po,\R}$. Also, if $\f$ is a face of $\Po$ then $N_{\Po,\R}\subset N_{\f,\R}$.

To the polyhedron $\Po$ is associated its normal fan $\Sigma_{\Po}$, which
is a rational fan in the 
vector space $N_{\R}/N_{\Po,\R}$ with respect to the lattice
$N/N_{\Po}$. Precisely, to each
face $\f$ of $\Po$ we associate the normal cone $\sigma_{\f,\Po}$ in $N_\R/N_{\Po,\R}$  given
by
$$
\sigma_{\f,\Po}:=\{v\in N_\R\,|\,\langle v,w_1-w_2\rangle\geq 0\text{ for every }w_1\in\f
\text{ and }w_2\in \Po\}/N_{\Po,\R}.
$$
(Notice that $\sigma_{\f,\Po}\subseteq N_{\f,\R}/N_{\Po,\R}$.) Then
$\Sigma_{\Po}$ is the collection of the $\sigma_{\f,\Po}$ as $\f$ runs
through the faces of $\Po$.

Let $\P_{\Po}$ be the  toric variety corresponding to the normal fan
$\Sigma_{\Po}$. Note that $\P_{\Po}$ is complete if and only if  $\Sigma_{\Po}$
has full support in $N_\R/N_{\Po,\R}$, if and only if $\Po$ is a polytope, that
is, a bounded polyhedron. 

Let $\Po$ be a polyhedron in $\PC$, and let
$\f$ be a face of $\Po$. Consider the normal cone $\sigma_{\f,\Po}$  in
$N_{\R}/N_{\Po,\R}$. The orbit closure associated to $\f$
in the toric variety $\P_{\Po}$ is isomorphic to the toric variety
associated to the fan 
$$
\Bigl\{\,\bigl(\sigma_{\f',\Po}+N_{\f, \R} \bigr)
/N_{\f, \R} \,\Bigl|\, \f' \textrm{ a face of }\Po \textrm{ contained in
} \f\,\Bigr\}.
$$
The natural pairing induces an isomorphism $ (N_\R/ N_{\f,\R})^\vee
\simeq M_{\f,\R}$, and the dual of the cone $\bigl(\sigma_{\f',\Po}+N_{\f, \R} \bigr)
/N_{\f, \R}$ in $M_{\f,\R}$ is the dual of $\sigma_{\f',\f}$. It
follows that the orbit closure 
associated to $\f$ in $\P_{\Po}$ is isomorphic to $\P_\f$, and we 
thus get a natural embedding $\P_\f \hookrightarrow \P_{\Po}$.  

\begin{defi}\label{deftorarr} \rm The \emph{toric arrangement} $\P_{\PC}$ associated to $\PC$
  is defined as follows: Consider the disjoint union of the toric
  varieties $\P_{\Po}$ for all $\Po$ in $\PC$, and identify $\P_{\f}$ with
  the orbit closure associated to $\f$ in $\P_{\Po}$, as described above,
  for each face $\f$ of each polyhedron $\Po$. 
\end{defi}

When $\PC$ is a tiling of $M_\R$ by polytopes, that is, when the
polyhedra are polytopes covering $M_\R$, we say that
$\P_{\PC}$ is a \emph{toric tiling}. 

The following proposition is straightforward.

\begin{prop}\label{tor-arr} Assume that $\PC$ is a locally finite rational polyhedral
  complex in $M_\R$. Then the toric arrangement $\P_{\PC}$ is naturally a
  scheme locally of finite type.  
\end{prop}

We are not careful enough in Definition~\ref{deftorarr}
to define completely the scheme structure of $\P_{\PC}$, as we
do not specify the scheme structure in neighborhoods of the
intersections of the various $\P_{\Po}$. However,
whatever structure is given under the conditions imposed by
Definition~\ref{deftorarr}, Proposition~\ref{tor-arr}
holds. Furthermore, we will
see definite scheme structures being given in the cases mentioned
in Example~\ref{extorarr}, in Subsection~\ref{tortilR} below and in
Subsection~\ref{bastortil}.

\subsection{Toric tiling $\mathbf R$ associated to the tiling of
  $C^1(G, \R)$ by hypercubes} \label{tortilR} Consider a finite
connected loopless graph $G=(V,E)$. Let $M:=C^1(G,\Z)$ and apply
the construction above to the associated
Voronoi tiling of  $M_\R$, the tiling by
hypercubes $\square_\alpha$ for $\alpha \in C^1(G, \Z)$. We get a
toric tiling, denoted $\mathbf R$, that we describe below. 

Since we have a
natural pairing $\langle.\,,.\rangle$ on
 $M_\R$, we may identify $M$ with the dual lattice $N:=M^\vee$.
We may thus view the normal fan $\Sigma_\alpha$ of each hypercube
$\square_\alpha\subseteq M_\R$ in the same space $M_\R$. 
It is the standard fan given by the
basis axes. More precisely, consider the family
$\mathcal A$ of all subsets $A\subset \E(G)$ which have the property that
for each edge $\{u,v\}\in E(G)$ at most one of the two possible
orientations $uv$ or $vu$ are in $A$, in other words, $A$ is an
orientation of a subset, denoted $E(A)$, of edges in $E$. The
\emph{positive cone} $\sigma_A$ associated to $A$ is defined as
$\sigma_A := \sum_{e\in A} \R_{\geq 0}(\chi_{\indm e}-\chi_{\indmbar{\ol e}})$. Then
$\Sigma_\alpha$ consists of all the cones $\sigma_A$. 

It follows that each maximal dimensional toric variety $\P_\alpha$ in the
tiling $\mathbf R$ is isomorphic to the product $\prod_{e\in E} \P^1$. 
And the tiling has the following description: 
Fix an orientation $\mathfrak o$ for the set of edges of the
graph. Consider an infinite number of copies
$\P^1_{e,i}$ of $\P^1$, indexed by the $e\in E^{\mathfrak o}$ and  integers $i\in \mathbb Z$. We
distinguish two points $0_{e,i}$ and
$\infty_{e,i}$ on each $\P^1_{e,i}$. In fact, for later use, we
give coordinates $(\x_{e,i}: \x_{\ol e,i})$ to each $\P^1_{e,i}$, with
$\x_{e,i}=0$ corresponding to $0_{e,i}$ and $\x_{\ol e,i}=0$ to
$\infty_{e,i}$. For each edge $e\in E^{\mathfrak o}$, we define 
${\mathbf R}_e$ as the locally finite scheme 
obtained by gluing $\P^1_{e,i}$ and $\P^1_{e,i+1}$ over the points $0_{e,i}$ and $\infty_{e,i+1}$ for each 
$i\in \mathbb Z$. The toric tiling $\mathbf R$ is then the
product of the $\mathbf R_e$ over the edges $e\in E^{\mathfrak o}$. 
In this description,
$$
\P_\alpha = \prod_{e\in E^{\mathfrak o}} \P^1_{e,
  \alpha_e}\quad\text{for each }\alpha\in C^1(G,\Z).
$$

We may likewise describe the other toric varieties in the tiling
$\mathbf R$, the ones
associated to the proper faces of the hypercubes. First,
for each edge $e\in E^{\mathfrak o}$ and half integer  
$a \in
\frac 12\Z-\Z$, we set  $\P_{e,a} := \{0_{e, \lfloor a\rfloor}\} =
\{\infty_{e, \lceil a \rceil}\} \subset \mathbf R_e$.  Then we extend the
definition of $\P_\alpha$ by setting 
$$
\P_\alpha := \prod_{e\in E^\mathfrak o} \P_{e,\alpha_e}\quad\text{for
  each  }\alpha\in C^1(G,\frac{1}{2}\Z).
$$
The whole tiling $\mathbf R$ is completely 
described by the toric varieties $\P_\alpha$ and the natural 
existing inclusions between them.

Removing from each $\P_\alpha$ all
those $\P_\beta$ contained in it, we obtain the interior of
$\P_\alpha$, an open subscheme denoted $\P_\alpha^*$. The
$\P_\alpha^*$ for $\alpha\in C^1(G,\frac{1}{2}\Z)$ form a
stratification of $\mathbf R$.
 
Finally, we may describe $\P_0$ as a gluing of affine schemes.
For each oriented edge $e\in \E$,
we denote by $\A^1_e$ the open subscheme of
 $\P^1_{e^{\mathfrak o},0}$ corresponding to $\x_{\ol e,0} \neq 0$ and
 by $\mathbf G_{\mathbf m,e}$ the open subscheme
 corresponding to $\x_{e,0}\x_{\ol e, 0} \neq 0$. Notice that
 $\mathbf G_{\mathbf m,e}=\mathbf G_{\mathbf m,\ol e}$,
 and we may thus abuse the notation by
 forgetting the orientation of $e$.
 For each $A\in \mathcal A$, let $S_A := \sigma_A^{\vee} \cap M$,
 the semigroup of integral points of the dual cone
 $\sigma_A^{\vee}$. The toric variety $\P_0$ is covered by the affine
 open subsets
 $\prod_{e\in A} \A^{1}_e \times
 \prod_{e\in E - E(A)}\mathbf G_{\mathbf m,e}$,
 each isomorphic to $\mathrm{Spec}\bigl(\k[S_A]\bigr)$.

 \subsection{Toric variety associated to the Voronoi cell of a graph}\label{sec:toricvar1}
 
Consider a finite connected loopless graph $G$ with uniform edge 
lengths equal to one with fixed orientation $\mathfrak o$, 
and denote by $\Vor_F(O)$ the Voronoi cell of the origin in $F_\R$ with respect to 
 $F_\Z$. In~\cite[Section 3]{AE1}, 
we described the face poset of $\Vor_F(O)$ in terms of the coherent acyclic orientations of
subgraphs of $G$, a description we reviewed in Section~\ref{sec:recap}. 
In this section, we use those results to describe the normal fan of
$\Vor_F(O)$, and 
to provide a complete description of 
 the corresponding toric variety, which we denote by $\P_G$. 
 
\smallskip

Recall the description of the face poset $\FP$ of the Voronoi cell
$\Vor_F(O)$: The faces are in bijection with coherent acyclic
orientations $D$ of cut subgraphs of $G$, and under this bijection, the face $\frak
 f_D$ associated to $D$ has 
support in the affine plane $F_D$ given by $F_D = \bigcap_{\beta \in \mathcal U_D} F_\beta$, where 
 $F_\beta =\bigl\{x\in F_\R,\bigl |\, 2\langle x, \beta\rangle
 =\|\beta\|^2\bigr\}$ and 
$\mathcal U_D$ is the set of all bond elements $\beta$ of $F_\Z$ with 
$\supp^+(\beta)\subseteq\E(D)$; see \cite[Subsection 3.4]{AE1}.
 
The pairing $\langle.\,,.\rangle$ on $M_\R=C^1(G, \R)$ gives the
 orthogonal decomposition $M_\R= F_\R \oplus \ker(d^*)$,
 which allows us to identify further $F_\R$ with its dual
 $F_\R^{\vee}$. Under this identification, the normal cone $\sigma_D$ 
 to the face $\mathfrak f_D$ of $\Vor_F(O)$ is the cone in $F_\R$
 generated 
by the $\beta\in \mathcal U_D$, that is, 
 \[
\sigma_D = \sum_{\beta \in \mathcal U_D}\R_{\geq 0}\beta.
\]
Therefore, we get the following
result. 

 \begin{prop} The normal fan $\Sigma_G$ of the Voronoi polytope
   $\Vor_F(O)$ 
consists of the cones $\sigma_D$ for $D\in \mathcal \AC$, 
 where $\sigma_D$  is the cone generated by all the $\beta\in \mathcal
 U_D$. In particular, 
the rays of $\Sigma_G$ are in bijection with
 the bond elements of $F_\Z$.
 \end{prop}
 
For each $D\in \AC$, consider the corresponding cone $\sigma_D \in
 \Sigma_G$. The dual cone 
$\sigma_D^{\vee} \subseteq F_\R$ consists of all the elements of
$F_\R$ which are non-negative on $\sigma_{D}$. Denote by $S_D :=
\sigma_D^{\vee}\cap F_\Z$ the semigroup of integral points in 
$\sigma_D^{\vee}$, and let 
 $U_D :  =\Spec(\k[S_D])$. 
The toric variety $\P_G$ is obtained by gluing the affine varieties $U_D$.   

Consider the inclusion $F_\R \hookrightarrow M_\R$. We have:

 \begin{prop}\label{prop:surjlattice} The lattice $F_\Z$ coincides with the lattice $M \cap F_\R$. In particular, the dual map $N \to F_\Z^{\vee}$ is surjective. 
 \end{prop}
 \begin{proof}
  The elements of $M \cap F_\R$ are of the form $d(f)$ for $f\in
  C^0(G, \R)$ such that $f(v) - f(u) \in \Z$ for each $e=uv\in
  \E(G)$. Since $G$ is connected, $f$ has to be of the form $h+c$ for
  $h\in C^0(G, \Z)$ and a constant $c\in \R$, and so $d(f) =d(h) \in
  F_\Z$. The second statement follows from the first, as the first implies
  that $M/F_\Z$ is torsion-free.
 \end{proof}

\begin{remark}
\rm The pairing $\langle.\,,.\rangle$ is integral on $M=C^1(G, \Z)$
and moreover identifies $M$ with its dual $N$. Under this
identification, the surjection $N\to F_\Z^{\vee}$ gets identified with
the restriction $M\to F_\Z$ of the orthogonal projection map.
\end{remark}

For each $D\in \AC$, recall that the set of oriented edges of $D$ 
is denoted by $\E(D)$, and gives the cone $\sigma_{\E(D)}$ in the normal
fan $\Sigma_0$ to $\square_0$. Since $\supp^+(\beta)\subseteq\E(D)$ for each
   $\beta\in\mathcal U_D$, the inclusion $F_\R \hookrightarrow M_\R$
   restricts to an inclusion of cones 
  $\sigma_D \hookrightarrow \sigma_{\E(D)}$. 

  \begin{prop}\label{sdd} We have $\sigma_D = F_\R \cap
    \sigma_{\E(D)}$. Also, $\R\sigma_D = F_\R \cap \R\sigma_{\E(D)}$.
  \end{prop}

    \begin{proof} Let $f$ be an element of $C^0(G, \mathbb R)$ such
      that $d(f)$ belongs to $\R\sigma_{\E(D)}$. Adding a constant to $f$ if necessary, we may assume that  the minimum
of $f$ over the vertices is zero. Let $0<\alpha_1< \alpha_2< \dots
<\alpha_l$ be all the values taken by $f$, and denote by $X_i$ the set
of all vertices $v$ with $f(v_i) \geq \alpha_i$ for each $i=1, \dots,
l$. Notice that $X_1 \supset X_2 \supset \dots \supset X_l$. We can write 
    \[f = \alpha_1 \chi_{\indbi{X}{1}} +  (\alpha_2-\alpha_1)\chi_{\indbi{X}{2}} + \dots + (\alpha_l - \alpha_{l-1}) \chi_{\indbi{X}{l}},\]
    so that 
    \[d(f) = \alpha_1 d(\chi_{\indbi{X}{1}}) + (\alpha_2-\alpha_1)d(\chi_{\indbi{X}{2}})+ \dots +(\alpha_l - \alpha_{l-1}) d(\chi_{\indbi{X}{l}}).\] 
Since the $X_i$ form a decreasing chain of  subsets, 
it follows that the supports $\supp^+(d(\chi_{\indbi{X}{i}}))$ are consistent
in their orientations. 
Moreover, since $\alpha_1>0$ and the $\alpha_i$ form an increasing sequence,
$\supp^+(d(\chi_{\indbi{X}{i}}))$ is a subset of $\supp^+(d(f))$. Each $d(\chi_{\indbi{X}{i}})$ is an element of $F_\Z$ and
is thus a sum of bond elements $\beta_{i,j} $ with
$\supp^+(\beta_{i,j}) \subseteq \supp^+(d(\chi_{\indbi{X}{i}}))$; see
\cite[Prop.~3.34]{AE1}.  If $d(f)\in \sigma_{\E(D)}$, that is, $\supp^+(d(f))\subset \E(D)$
then  all the $\beta_{i,j}$ belong to $\U_D$. This implies that $d(f)$ is a sum of bond elements belonging
  to $\U_D$, that is, 
$d(f) \in \sigma_D$, and the first statement follows.

As for the second statement, notice that at any rate all the
$\beta_{i,j}$ belong to $F_\R\cap\R\sigma_{\E(D)}$. Thus we may
suppose that $f=\chi_{\ind X}$ for some cut $X\subset V$. Since $D\in\AC$,  there is a partition
$V_1,\dots,V_s$ of $V$ such that $e\in\E(D)$ if and only if $e$ is an
edge connecting $V_i$ to $V_j$ with $i<j$. We may suppose $G[V_i]$ is
connected for each $i$. Let $X_i:=X\cap V_i$ for each
$i$. Any edge connecting $V_i-X_i$ to $X_i$ would be in
$\supp^+(d(f))$ but would not be supported over $E(D)$. Thus
$\E(V_i-X_i,X_i)=\emptyset$ for each $i$. Furthermore, since $G[V_i]$
is connected, either $X_i=\emptyset$ or $X_i=V_i$. It follows that $X$ is a
union of some of the $V_i$. We need only prove now that
$d(\chi_{\indbi{V}{i}})\in\R\sigma_D$ for each $i$. But this follows because 
$d(\chi_{\indbi{V}{i}}+\cdots+\chi_{\indbi{V}{s}})\in F_\R\cap\sigma_{\E(D)}$, and thus 
$d(\chi_{\indbi{V}{i}}+\cdots+\chi_{\indbi{V}{s}})\in \sigma_D$ for each $i$ by the
first statement. 
  \end{proof}

  We prove now the surjectivity of the dual map.

  \begin{prop}\label{prop:surjlatticesigma} Notations as above, the map of dual cones
    $\sigma_{\E(D)}^{\vee} \rightarrow \sigma_D^{\vee}$ is
    surjective. Furthermore, the induced map on integral points
    $\sigma_{\E(D)}^{\vee} \cap M \rightarrow \sigma_D^{\vee} \cap
    F_{\mathbb Z}$ is surjective.  
  \end{prop}

  \begin{proof} Let $\phi$ be an element of the dual cone
    $\sigma_D^{\vee}$ in $F_\R$ (resp.~in 
$F_{\mathbb Z}$).  We want to show the existence of an element $\eta
\in \sigma_{\E(D)}^{\vee}$ (resp.~$\eta \in \sigma_{\E(D)}^{\vee} \cap
M$), thus an element $\eta\in C^1(G,\R)$ (resp.~$\eta\in C^1(G,\Z)$)  non-negative on every edge of $\E(D)$, whose 
orthogonal projection to $F_\R$ is $\phi$. 
  
  Since $\phi$ is in $\sigma_D^{\vee}$, for each bond element $\beta
  \in \U_D$ we have $\langle \phi, \beta\rangle \geq 0$. Now, it
  follows from \cite[Prop.~3.34]{AE1} that any element $\alpha \in
  F_\R$ with $\supp^{+}(\alpha) \subset \E(D)$ is a linear combination
  with positive coefficients of bond elements in $\U_D$. Then
$$
\langle \phi, \alpha \rangle \geq 0 \qquad \textrm{for all }\alpha\in
F_\R \textrm{ with }  \supp^{+}(\alpha) \subset \E(D).
$$
This shows that for each $f\in C^0(G, \R)$ with $\supp^{+}(d(f)) \subset \E(D)$, we have 
\[\langle \phi, d(f) \rangle  = \langle d^*(\phi), f \rangle \geq 0.\] 
In particular, for each cut $X\subseteq V$  with
$\E(V- X, X) \subset \E(D)$, we have $d^*(\phi)(X) \geq
0$. Furthermore, if
$\phi\in F_\Z$ then $d^*(\phi)\in C^0(G,\Z)$.

It follows thus from the next proposition that there exists 
$\eta \in C^1(G, \R)$ (resp.~$\eta \in C^1(G, \mathbb Z)$) with
$\eta(e) \geq 0$ 
for all $e\in \E(D)$ such that 
\[d^*(\eta) = d^*(\phi). \]
Since $\eta-\phi\in\text{ker}(d^*)$, it follows that $\eta$ projects
to $\phi$. 
\end{proof}

Recall that $H_{0,\R}$ is the set of all elements $h \in C^0(G,\R)$
with $\sum_{v\in V} h(v)=0.$ Moreover, the map $d^*: C^1(G,\mathbb R)
\to C^0(G,\mathbb R)$ surjects onto $H_{0,\mathbb R}$. We are
interested in characterizing those elements of $H_{0,\mathbb R}$ which
are in the image of elements $\eta\in C^1(G,\mathbb R)$ that are positive on all oriented edges in $\E(D)$. We have the following theorem.

\begin{prop} Let $G=(V,E)$ be a finite loopless graph and $D$ an orientation of a
subset of $E$. Let $h\in C^0(G,
  \R)$ be an element in $H_{0,\R}$. The following two statements are equivalent:
\begin{enumerate}
\item For all subsets $X \subseteq V$ with $\E(V- X, X) \subseteq \E(D)$, we have $h(X) \geq 0$.
\item There exists $\eta \in C^1(G, \R)$ with $d^*(\eta)=h$ such that $\eta(e) \geq 0$ for all $e\in \E(D)$.
\end{enumerate}
Moreover, in the case $h \in C^0(G, \mathbb Z)$, we can modify $(2)$ by requiring $\eta \in C^1(G, \mathbb Z)$.
\end{prop}

\begin{proof} The implication  $(2) \Rightarrow (1)$ follows from observing that for an element $\eta$ as in (2), we have  for each subset $X \subseteq V$:
\[h(X) = d^*(\eta)(X) = \eta\bigl(\E(V- X, X)\bigr) \geq 0.\]

We prove the other implication, $(1) \Rightarrow (2)$. 
Denote by $V^+$ the set of all vertices $v$ with $h(v)>0$, by $V^-$
the set of all vertices $v$ with $h(v)<0$, and by $V^0$ the set of all vertices $v$ with $h(v)=0$. We proceed by induction on $|V^+| + |V^-|$.  If this quantity is zero, then $h \equiv 0$, and  $\eta=0$ does the job.

Suppose now we have $|V^+| + |V^-| =n>0$. Since $h\in H_{0,\R}$, both
$V^+$ and $V^-$ are nonempty. 
We claim there exists an oriented path $P$ compatible with $D$ which connects a vertex 
$v$ of $V^-$ to a vertex $u$ of $V^+$. Indeed, let $X$ be the set of vertices that can be connected to a vertex of
$V^-$ by an oriented path compatible with $D$, in the sense that no
oriented edge $e$ in $P$ satisfies $\ol e\in\E(D)$. Then $\E(V- X, X)
\subset \E(D)$. If $X\cap V^+=\emptyset$ then $h(X) = h(V^-) <0$, which would contradict the assumption on $h$. 

Let $P$ be such an oriented path from a vertex $v$ in $V^-$ to a
vertex $u$ in $V^+$. We view $P$ as an element of $C^1(G,\Z)$ taking
value one at each oriented edge $e \in \E(P)$. For each $\epsilon \geq 0$,
define $h_\epsilon := h -\epsilon d^*(P) = h- \epsilon \chi_{\indm u} +
\epsilon \chi_{\indm v}$.  So $h_\epsilon(u) = h(u) -\epsilon$, $h_\epsilon(v)
= h(v)+\epsilon$, and $h_\epsilon(w) = h(w)$ for any $w$ different
from $u$ and $v$. Let $\epsilon$ be the greatest nonnegative real number
such that $h_\epsilon(u)\geq 0$, $h_\epsilon(v)\leq 0$ and 
$h_\epsilon(X)\geq 0$ for every cut $X\subseteq V$ with $\E(V- X, X) \subset \E(D)$.
If $h \in C^0(G, \mathbb Z)$ then $\epsilon$ is an integer, 
and so $h_\epsilon \in C^0(G, \Z)$.

If $h_\epsilon(u)=0$ or $h_\epsilon(v)=0$, then the value of $|V^+| +
|V^-|$ for $h_\epsilon$ has 
decreased, hence by the induction hypothesis there exists
$\eta_\epsilon\in C^1(G,\R)$ with $\eta_\epsilon$ nonnegative on
$\E(D)$ such that 
$d^*(\eta_\epsilon) = h_\epsilon$. Also, $\eta_\epsilon$ is integral if
so is $h$. Let $\eta := \eta_\epsilon + \epsilon P$.  Then $\eta$ is
nonnegative on $\E(D)$ and $d^*(\eta) = h$. Furthermore, if $h$ is
integral, so is $\eta$.  This finishes the proof in this case.

In the remaining case, we get a nonempty proper subset $X\subsetneq V$
of vertices with 
$\E(V- X, X) \subset \E(D)$ such that  $h_\epsilon(X) = 0$. 
In this case, one verifies that both the restrictions of $h_\epsilon$
to $X$ and to $V- X$ 
satisfy the conditions of the proposition. A second induction on the
size of $V$ then shows 
the existence of $\eta_1 \in C^1(G[X], \R)$ and $\eta_2\in C^1(G[V-
X], \R)$ such that $\eta_1$ is nonnegative on all edges of $\E(D)$ 
with both endpoints in $X$, and $\eta_2$ is nonnegative 
on all edges of $\E(D)$ with both endpoints in $V-X$, 
and such that $d^*(\eta_1) = h_{\epsilon |X}$ in $G[X]$  and $d^*(\eta_2) =
h_{\epsilon |V- X}$ in $G[V- X]$. 
Furthermore, we may take $\eta_1$ and $\eta_2$ integral in the case $h$,
and so $h_\epsilon$ is integral.  
Define $\eta_\epsilon \in C^1(G, \R)$ to be equal to $\eta_1$ on
$G[X]$, 
equal to $\eta_2$ on $G[V- X]$, and equal to zero on all edges of 
$\E(V- X, X)$. We have $d^*(\eta_\epsilon) =h_\epsilon$, and
$\eta_\epsilon$ is obviously nonnegative 
on all edges of $\E(D)$. Moreover, it is integral if $h$ is
integral. Taking now $\eta := \eta_\epsilon + \epsilon P$, we conclude
as before.
\end{proof}

Consider now the map of cones $\sigma_{\E(D)}^{\vee} \to \sigma_D^{\vee}$.   Using Proposition~\ref{prop:surjlatticesigma}, we get a surjective map of semigroups 
  \[
S_{\E(D)} = \sigma_{\E(D)}^{\vee} \cap M \longrightarrow
S_{D}=\sigma_D^{\vee} \cap  F_\Z,
\]
which induces a closed embedding of $U_D\to U_{\E(D)}$  in the affine variety 
$$
U_{\E(D)}:=\Spec[\k[S_{\E(D)}]] \simeq \prod_{e\in \E(D)} \A^1_e \times 
  \prod_{e\in E - E(D)} \mathbf G_{\mathbf m,e},
$$
which is the affine open subvariety of $\P_0$ corresponding to the
cone $\sigma_{\E(D)}$. 

The maps $U_D\to U_{\E(D)}$ glue to a map $\alpha_0\:\P_G\to\P_0$ satisfying 
$\alpha_0^{-1}(U_{\E(D)})=U_D$ for each coherent acyclic orientation
$D$ of a cut subgraph of $G$. Thus $\alpha_0$ is a closed embedding. Its
image is described by the following theorem:

\begin{thm}\label{thm:toric} Notations as above, the map $\alpha_0$ 
induces an isomorphism from $\P_G$ to the subvariety of
$\P_0=\prod_{e\in E^{\mathfrak o}} \P_{e,0}$ given by the equations
 \[
\forall \textrm{ oriented cycle } \gamma \textrm { in $G$}, \quad
\prod_{e\in \gamma\cap E^{\mathfrak o}}
\x_{e,0}\prod_{e\in\bar\gamma\cap E^{\mathfrak o}} \x_{\ol e,0} = 
\prod_{e\in\bar\gamma\cap E^{\mathfrak o}} \x_{e,0}
\prod_{e\in \gamma\cap E^{\mathfrak o}}
\x_{\ol e,0}.
\]
  \end{thm}

Before proving the theorem we need to introduce notation and prove a proposition.

For each subset $A \subset \E(G)$ which is an orientation of a subset
of edges of $G$,  
let $\k[\x_A]$ be the polynomial ring generated by the variables 
$\x_e$ for $e\in A$, and define
$$
\A^{A}: = \Spec[\k[\x_A]]=\prod_{e\in A} \A^1_{e}.
$$

If $D$ is an acyclic orientation of the whole graph $G$, then
$U_{\E(D)}=\A^{\E(D)}$, and we get that $\alpha_0$ restricts to an
embedding
$U_D \hookrightarrow \A^{\E(D)}$. In this case, we have the
following proposition.

\begin{prop}\label{UDED}
Let $D$ be an acyclic orientation of the whole graph $G$. Then the 
embedding $U_D \hookrightarrow \A^{\E(D)}$ identifies $U_D$ with the
subvariety of the affine space $\A^{\E(D)}$ given by the equations
$$
(*)\quad \forall \,\, \textrm{ oriented cycle }\, \gamma \textrm{ in }
G, \,\, \prod_{e\in \gamma \cap \E(D)} \x_e = \prod_{e\in \bar \gamma
  \cap \E(D)} \x_e.
$$
  \end{prop}

  \begin{proof} Let $\alpha$ denote the surjection $\k[S_{\E(D)}] \to \k[S_D]$. 
We show first that  for each oriented cycle $\gamma$  in $G$, the
difference  
$\prod_{e\in \gamma \cap \E(D)} \x_e - \prod_{e\in \bar \gamma \cap \E(D)} \x_e$ is in the kernel of $\alpha$.
  
We view each oriented cycle $\gamma$ as the element
$\sum_{e\in\gamma}(\chi_{\indm e}-\chi_{\indmbar{\ol e}})\in M=C^1(G,\Z)$, which we
will also denote by $\gamma$. Also, we let
$$
\gamma_+:=\sum_{e\in\gamma\cap\E(D)}(\chi_{\indm e}-\chi_{\indmbar{\ol e}})\quad\text{ and }\quad
\gamma_- := \sum_{e\in\bar\gamma\cap\E(D)}(\chi_{\indm e}-\chi_{\indmbar{\ol e}}).
$$
Thus $\supp^+(\gamma_\pm) \subseteq \E(D)$ and $\gamma = \gamma_+ -
\gamma_-$. 

Now, if $\gamma$ is an oriented
cycle, $\langle \gamma, \beta \rangle =0$ for every $\beta\in F_\R$. 
Thus $\langle \gamma_+, \beta\rangle = \langle
\gamma_-, \beta\rangle$ for every $\beta\in F_{\R}$, that is, 
the orthogonal projections of $\gamma_+$ and $\gamma_-$ in $F_\R$
coincide.

It follows that $\alpha(\prod_{e\in \gamma \cap \E(D)}\x_e )= \alpha(\x^{\gamma_+}) = \alpha(\x^{\gamma_-}) = \alpha(\prod_{e\in \bar \gamma\cap \E(D)}\x_e),$ which shows that the kernel of $\alpha$ contains the difference $\prod_{e\in \gamma \cap \E(D)} \x_e - \prod_{e\in \bar \gamma \cap \E(D)} \x_e$.
 
\medskip
 
We now prove that the differences 
$\prod_{e\in \gamma \cap \E(D)} \x_e - \prod_{e\in \bar \gamma \cap
  \E(D)} \x_e$ generate the kernel of $\alpha$. Indeed, the kernel of $\alpha$ is generated by the elements of the form 
  $\x^{\eta_1} - \x^{\eta_2}$ for $\eta_1$ and $\eta_2$ in $M$ verifying
  \begin{enumerate}
   \item$\supp^+(\eta_i) \subseteq \E(D)$ for $i=1,2$ and
   \item $\langle \eta_1, \beta\rangle = \langle \eta_2, \beta\rangle$
     for every $\beta \in \mathcal U_D$.
  \end{enumerate}
Cancelling out the common factors, we may further assume that
$\supp^+(\eta_1)\cap \supp^+(\eta_2) =\emptyset$.  Note that, since
$D$ is an acyclic orientation of the whole graph, the
$\beta\in \mathcal U_D$ generate $F_\R$, whence (2) implies that
$\eta_1-\eta_2$ is orthogonal 
to the cut lattice. 

Consider the oriented graph $G'$ on the vertex set $V(G)$ which for each oriented edge 
$e\in \E(D)$ contains $\eta_1(e)$ oriented edges parallel to $e$ and
$\eta_2(e)$ oriented edges parallel to $\ol e$. Since $\eta_1-\eta_2$ is
orthogonal to the cut lattice, we get that 
$\langle \eta_1-\eta_2, d(\chi_{\indm v})\rangle =0$ for the cut defined by
each vertex $v$ in $G$, which implies that the graph
$G'$ is Eulerian, that is, at each vertex $v$ the in-degree of 
$v$ is equal to its out-degree. Then the (oriented) edges of the
Eulerian oriented graph $G'$ can be partitioned into an 
edge-disjoint union of oriented cycles $\gamma_1, \dots,
\gamma_N$. Notice that for each $i=1,\dots,N$ and each $e\in\E$
there is at most one oriented edge of $\gamma_i$ among those parallel to
$e$, and thus we can view $\gamma_i$ as an oriented cycle of $G$ as
well. Thus, we have 
$$
\eta_1=\sum_{i=1}^N\gamma_{i,+} \quad \textrm{and} \quad
\eta_2=\sum_{i=1}^N \gamma_{i,-}.
$$
Let $\textsc{A}_i :=\x^{\gamma_{i,+}}$ and $\textsc{B}_i := \x^{\gamma_{i,-}}$ for $i=1,
\dots, N$. Put $\textsc{B}_0:=1$ and $\textsc{A}_{N+1}:=1$. Since
\begin{align*}
\x^{\eta_1} - \x^{\eta_2} =  \textsc{A}_1\cdots \textsc{A}_N - \textsc{B}_1\cdots \textsc{B}_N = 
\sum_{i=0}^{N-1} \textsc{B}_0\cdots \textsc{B}_i (\textsc{A}_{i+1}-\textsc{B}_{i+1})\textsc{A}_{i+2}\cdots \textsc{A}_{N+1},
\end{align*}
it follows that $\x^{\eta_1} - \x^{\eta_2} \in
(\x^{\gamma_{1,+}}-\x^{\gamma_{1,-}},\dots,
\x^{\gamma_{N,+}}-\x^{\gamma_{N,-}})$, finishing the proof.
\end{proof}

We can now prove Theorem~\ref{thm:toric}.

 \begin{proof}[Proof of Theorem~\ref{thm:toric}]  By Proposition~\ref{UDED}, the equations giving
 $U_D$ in $U_{\E(D)}$ are just the dehomogenizations of the equations
 given in the theorem, if $D$ is an acyclic orientation of the whole graph. Now, 
every $D\in\AC$ can be completed to an
acyclic orientation of the whole graph $G$. Indeed, as the two posets
$\FP$ and $\AC$ are isomorphic by Theorem~\ref{thm:vor}, 
and the vertices of $\text{Vor}_F(O)$
correspond to acyclic orientations of the whole graph by \cite[Lem.~3.42]{AE1}, 
the statement is equivalent to saying that every face of
$\text{Vor}_F(O)$ contains a vertex of $\text{Vor}_F(O)$. (Or we may
use Lemma~\ref{completeorient} below.) Thus, the
$U_D$ for acyclic orientations $D$ of the whole graph 
form a covering of $\P_G$. It follows that the points on $\P_G$ satisfy
 the equations given in the theorem. 

Moreover, to conclude the proof
 it is enough to show that every point $p$ on $\P_0$ satisfying the
 equations given in the theorem lies on $U_{\E(D)}$ for a certain acyclic orientation $D$ of the whole 
 graph $G$. This is the case, as it follows from the equations that
 for each oriented cycle $\gamma$ in $G$, the product
 $\prod_{e\in\gamma}\x_{e,0}$ vanishes on $p$ if and only if so
 does $\prod_{e\in\gamma}\x_{\ol e,0}$. Thus, if we let $A\subseteq \E$ be
 the subset of all oriented edges $e$ such that $\x_{e,0}$ vanishes on
 $p$ we have that $A$ is an orientation of a subset of edges of $G$
 satisfying the following property: If $\gamma$ is an oriented cycle
 in $G$, then $\gamma\cap A=\emptyset$ if and only if $\gamma\cap\ol
 A=\emptyset$. The following lemma finishes the proof.
 \end{proof}

\begin{lemma}\label{completeorient} Let $A\subseteq\E$ be an orientation of a subset of
  edges of $G$ such that for each oriented cycle $\gamma$ of $G$, we
  have that $\gamma\cap A=\emptyset$ if and only if 
$\ol\gamma\cap A=\emptyset$. Then $A$ can be completed to an acyclic
 orientation of all the edges of $G$.
\end{lemma}

\begin{proof} Let $G'$ be the subgraph of $G$ obtained by removing all
  the edges supporting $A$. Give $G'$ any acyclic orientation, for
  instance, following an ordering of the vertices, and
  extend it to an orientation $\mathfrak o$ of $G$ by adding $A$. We
  claim $\mathfrak o$ is acyclic. Indeed, if $\gamma$ is an oriented
  cycle of $G$ following this orientation, since $G'$ has an induced
  acyclic orientation, then there must be $e\in\gamma\cap A$. But
  then, by hypothesis, there is $e'\in\gamma\cap\ol A$, hence
  $e'\in\gamma$ which does not follow the orientation $\mathfrak o$,
  a contradiction.
\end{proof}

\subsection{Toric orbits}\label{sec:toricorb} 
The toric variety $\P_G$ is naturally stratified by its toric orbits. We describe now this stratification of $\P_G$.  More precisely, let
$D$ be a coherent acyclic orientation of a cut subgraph of $G$. 
Denote by $T_{D}$ the torus orbit
associated to the normal cone 
$\sigma_D \in \Sigma_G$ to $\frak f_D$. 
Our aim is to describe the closure  $\overline T_D$ of $T_D$ in $\P_G$.

Let $Y_1,\dots, Y_d$ be the connected components of the graph obtained
by removing all edges of $E(D)$.
By ordering the $Y_i$
appropriately, we may and will assume that an oriented edge $e\in\E(G)$ is in $D$ if
and only if $e$ connects $Y_i$ to $Y_j$ for $i<j$. For each $i=1,\dots,d$, let $M_i:=C^1(Y_i,\Z)$ and $N_i:=M_i^\vee$.  Also, let $F_{i,\R}$ and $F_{i,\Z}$ be the cut space and the cut
lattice of $Y_i$. Denote by $\Sigma_i$ the normal fan of
$\Vor_{F_i}(O)\subset F_{i,\R}$. Also, let $E_i$ be the set of edges
of $Y_i$ and $\mathfrak o_i$ the orientation induced from the
orientation $\mathfrak o$ of $G$. As before, we will identify $M_{i,\R}$ with $N_{i,\R}$ for each
$i$. We may also view each $M_{i,\R}$ in $M_\R$ under extension by zero.

Consider first the subvariety $\P_{D,0}$ of $\P_0$, which consists of
all the points $x\in \P_0$ 
whose $e$-th coordinates 
$(x_{e,0}:x_{\ol e,0})$ for $e\in \E(D)$ satisfy $x_{e,0}=0$. This is the
closure of the torus associated to the cone $\sigma_{\E(D)}$ in the
normal fan of the hypercube $\square_0$. It is thus the toric variety associated to the fan $\Sigma_{D,0}$ in
$C^1(G,\R)/\R\sigma_{\E(D)}$ consisting of the cones 
$$
\frac{\sigma_{\E(D')}+\R\sigma_{\E(D)}}{\R\sigma_{\E(D)}}
$$
for $D'\in\AC$ with $D' \preceq D$, or equivalently, $\E(D') \supseteq\E(D)$.

The subvariety $\P_{D,0}$ is naturally isomorphic to the product
$\prod_{i=1}^{d}\P_{i,0}$ where each $\P_{i,0}$ itself is isomorphic
to the product $\prod_{e\in E_i^{\mathfrak o_i}} \P^1_e$.  The isomorphism is induced by the isomorphism
\begin{equation}\label{isotot}
\bigoplus_{i=1}^{d}C^1(Y_i,\R) \longrightarrow C^1(G,\R) \longrightarrow
\frac{C^1(G,\R)}{\R\sigma_{\E(D)}},
\end{equation}
as it takes the product of the standard fans $\Sigma_{i,0}$ associated
to the $Y_i$ onto $\Sigma_{D,0}$.  In particular, denoting by $\P_{Y_i}$ the toric varieties associated
to the fans $\Sigma_i$, Theorem~\ref{thm:toric} yields embeddings $\P_{Y_i} \hookrightarrow \P_{i,0}$. Taking the product 
of these embeddings, we thus get an embedding 
\[\alpha_D:\prod_{i=1}^{d}\P_{Y_i} \hooklongrightarrow \P_{D,0}.\]

\begin{prop}\label{strat0} Notations as above, the image of $\alpha_D$ is
  $\P_G\cap\P_{D,0}$.
\end{prop}

\begin{proof} The equations defining $\P_G\cap\P_{D,0}$ in $\P_0$ are
  those given by Theorem~\ref{thm:toric} and the equations $\x_{e,0}=0$
  for all $e\in\E(D)$. Now, if $\gamma$ is an oriented cycle in $Y_i$,
  for each $i$, the equation prescribed by Theorem~\ref{thm:toric}
  for $\P_G$ is the same as the equation prescribed for $\P_{Y_i}$ in
  $\P_{i,0}$. Thus we need only show that the equation given by an
  oriented cycle $\gamma$ in $G$ not supported entirely in any of the
  $Y_i$ is a consequence of the equations $\x_{e,0}=0$
  for all $e\in\E(D)$. 

Indeed, let $\gamma$ be such a cycle. Since $\gamma$ is
not entirely supported in any $Y_i$, there is $e\in\gamma$ connecting
$Y_i$ to $Y_j$ for $i\neq j$. Assume, without loss of generality as it
will be clear, that $i<j$. Since $\gamma$ is a cycle
there will be $f\in\gamma$ connecting $Y_l$ to $Y_m$ for $l>m$. Thus 
$e\in\E(D)$ and $\ol f\in\E(D)$, whence $\x_{e,0}=0$ and $\x_{\ol
  f,0}=0$ on $\P_G\cap\P_{D,0}$. But then the equation associated to
$\gamma$ given by Theorem~\ref{thm:toric} is satisfied on
$\P_G\cap\P_{D,0}$ because both sides are equal to 0. 
\end{proof}

\begin{thm}\label{thm:stratification} Notations as above, the embedding
  $\alpha_D$ induces an isomorphism from the product of toric varieties associated to the connected components of 
$G - E(D)$ to the closure of the orbit 
$T_D$ in $\P_G$.
\end{thm}

\begin{proof} Composition~\eqref{isotot} induces an isomorphism from 
$\oplus_i F_{i,\R}$ onto the image of $F_\R$ in the
quotient. This image is isomorphic to the quotient of $F_\R$ by
$F_\R\cap\R\sigma_{\E(D)}$. But the latter is $\R\sigma_D$, by Proposition~\ref{sdd}. So we
obtain a natural isomorphism
\begin{equation}\label{isoFF}
\bigoplus_{i=1}^{d}F_{i,\R} \longrightarrow \frac{F_\R}{\R\sigma_D}.
\end{equation}

The orbit closure of $T_D$ is isomorphic to the toric 
variety $\P(\Sigma_D)$ defined by the fan $\Sigma_D$ in
$F_\R/\R\sigma_D$ 
consisting of the cones 
$$
\tau_{D'}:=\Bigl(\sigma_{D'}+\R\sigma_D\Bigr) /\R\sigma_D
$$ 
for $D'\in \AC$ with 
 $D' \preceq D$, or equivalently, $\E(D') \supseteq\E(D)$.  To finish the proof, we will
show that the fan $\Sigma_D$ is the sum of the fans $\Sigma_i$ under the
isomorphism~\eqref{isoFF}.

For each $D'\in\AC$ with $\E(D') \supseteq\E(D)$, if we let $D'_i$ denote
the orientation of the edges in the graph $Y_i$ given by $D'$ for each
$i=1,\dots,d$, it follows that $D'_i$ is a coherent acyclic
orientation of $Y_i$. Let $\sigma_i$ denote the corresponding cone
in $F_{i,\R}$. We claim that, under the isomorphism~\eqref{isoFF}, the
sum of cones $\oplus_i\sigma_i$ is taken to $\tau_{D'}$. 

Indeed, each $\sigma_i$ is generated by elements of the form
$\beta|_{\E(Y_i)}$ where $\beta:=d(\chi_{\ind Z})$ for
  cuts $Z\subseteq V$ contained in $V(Y_i)$ such that 
$\supp^+(\beta|_{\E(Y_i)})\subseteq \E(D'_i)$. Let
$W$ be the union of $Z$ with the $V(Y_l)$ for $l>i$ and put 
$\alpha:=d(\chi_{\ind W})$. Clearly,
$\alpha|_{\E(Y_i)}=\beta|_{\E(Y_i)}$. 
Also, $\alpha\in\sigma_{\E(D')}$ and thus
$\alpha\in\sigma_{D'}$ by Proposition~\ref{sdd}. Furthermore, if
$\beta^*$ denotes the extension by zero of
$\beta|_{\E(Y_i)}$ to $C^1(G,\R)$, then
$\alpha-\beta^*\in\R\sigma_{\E(D)}$. Hence the image of 
$\beta|_{\E(Y_i)}$ in $F_{\R}/\R\sigma_D$ is in $\tau_{D'}$. 

Conversely, $\sigma_{D'}$ is generated by elements of the form 
$\beta:=d(\chi_{\ind Z})$ for cuts $Z\subseteq V$ such that
$\supp^+(\beta)\subseteq\E(D')$. Put $Z_i:=Z\cap V(Y_i)$ for each $i$ and
set $\beta_i:=d(\chi_{\indbi{Z}{i}})$. Then
$\supp^+(\beta_i|_{\E(Y_i)})\subseteq \E(D'_i)$ and thus 
$\beta_i|_{\E(Y_i)}\in\sigma_i$ by Proposition~\ref{sdd}. Let
$\beta^*_i\in C^1(G,\R)$ be the extension by zero of $\beta_i|_{\E(Y_i)}$ for
each $i$. Then $\beta-\sum_i\beta^*_i\in\R\sigma_{\E(D)}$.

On the other hand, if $D'_i$ is a coherent acyclic orientation of
a cut subgraph of $Y_i$ for each $i=1,\dots,d$, let $D'$ denote the orientation in the
graph $G$ given by the $D'_i$ and the $D_i$.  It is clear that $D'$ is
a coherent acyclic orientation of a cut subgraph of $G$ that induces $D'_i$ on each
graph $Y_i$ and satisfies $\E(D') \supseteq\E(D)$. It follows that 
the fan $\Sigma_D$ is the sum of the fans $\Sigma_i$ under the
isomorphism~\eqref{isoFF}. 

Therefore, the toric variety
$\P(\Sigma_D)$ associated to $\Sigma_D$ is isomorphic to the product
of the toric varieties associated to the $\Sigma_i$, namely $\prod \P_{Y_i}$.
\end{proof}

\section{Mixed toric tilings I: Scheme structure}\label{mixtortil}
 
\subsection{Basic toric tiling associated to $\ell$ and
  $\m$}\label{bastortil}
In this
section, we describe the 
toric tiling associated to the mixed tiling defined by a graph $G$, an
edge length function $\ell$ 
and a twisting $\m$.
  
Consider thus the general situation described
in~\cite[Section~5]{AE1} and reviewed in Section~\ref{sec:recap}. Let $G$ be a
loopless connected graph. Let $\ell: E
\rightarrow \N$ be an 
edge length function, and $\m\in C^1(G, \Z)$ a twisting. Let 
$\Vor_H^\m$ be the Voronoi decomposition of $H_{0,\R}$. 
The maximal dimensional cells of $\Vor_H^\m$ are of the form 
$\Vor^\m_H(f)$ for $f\in C^0(G, \Z)$ with connected $G^\m_f$. We have 
$\Vor^\m_H(f) = d^*(\dl^\m_f) + \Vor_{G^\m_f}(O)$.
By Theorem~\ref{thm:toric} applied to $G^\m_f$, the toric variety
$\P_{G^\m_f}$ embeds naturally in $\P_{\dl^\m_f}$, and is given by the
equations associated to the oriented cycles in $G^\m_f$. We denote
this subvariety of $\P_{\dl^\m_f}$ by $\P_{\ell,\m,f}$. Denote by $Y_{\ell, \m}^{\bt}$ the union of the
$\P_{\ell,\m,f}$ in $\mathbf R$. We call this union the \emph{basic
  toric tiling} or the \emph{basic
  toric arrangement} associated to $G$, $\ell$ and $\m$. The names are
justified by the following theorem and its corollary.

First, we introduce notation and a lemma. 
For each $\mathfrak m\in C^1(G,\Z)$, each $f\in C^0(G,\Z)$ and each
$n\in\mathbb N$, set
$$
\delta^{\mathfrak m, n}_e(f):=
\Big\lfloor\frac{f(\he_e)-f(\te_e)+n\mathfrak m_e}{n\ell_e}\Big\rfloor
\quad\text{for each }e\in\E,
$$
and let $\mathfrak d_f^{\mathfrak m,n}\in C^1(G,\frac{1}{2}\Z)$ be
defined by putting
$$
\mathfrak d_f^{\mathfrak m,n}(e)=\frac{1}{2}
(\delta^{\mathfrak m,n}_e(f)-\delta^{\mathfrak m,n}_{\ol e}(f))
$$
for each $e\in\mathbb E$. It is the same definition given in Section
\ref{sec:recap} for $n=1$, in which case the superscript is dropped.

\begin{lemma}\label{e1e2} Let $f,h\in C^0(G,\Z)$ and $n,p\in\N$. Let $\gamma$ be
  an oriented cycle in $G$ such that $\dl^{\m,p}_h(e)\in\Z$ for each
  $e\in\gamma$.  Then there is $e_1\in\gamma$ such that $\dl^{\m,n}_f(e_1)<\dl^{\m,p}_h(e_1)$ if and only if there is
$e_2\in\gamma$ such that $\dl^{\m,n}_f(e_2)>\dl^{\m,p}_h(e_2)$.
\end{lemma}

\begin{proof} Since $\mathfrak d^{\mathfrak
  m,p}_h(e)\in\Z$ for each $e\in\gamma$, it follows that 
\begin{equation}\label{sumgpre}
\sum_{e\in\gamma}\mathfrak d^{\mathfrak
  m,p}_h(e)=\sum_{e\in\gamma}\frac{\mathfrak m_{e}}{\ell_{e}}.
\end{equation}
Now, also
\begin{equation}\label{sumfpre}
\sum_{e\in\gamma}\frac{f(\he_e)-f(\te_e)+n\mathfrak
  m_{e}}{n\ell_e}=\sum_{e\in\gamma}\frac{\mathfrak m_{e}}{\ell_{e}}.
\end{equation}

Suppose there is $e_1\in\gamma$ such that
$\dl^{\m,n}_f(e_1)<\dl^{\m,p}_h(e_1)$. Since $\dl^{\m,p}_h(e_1)$ is an
integer, it follows that the summand on the left-hand side of
Equation~\eqref{sumfpre} corresponding to $e_1$ is smaller than 
$\dl^{\m,p}_h(e_1)$. But then it follows from 
Equations~\eqref{sumgpre}~and~\eqref{sumfpre} that there is
$e_2\in\gamma$ such that the summand on the left-hand side of
Equation~\eqref{sumfpre} corresponding to $e_2$ is bigger than
$\dl^{\m,p}_h(e_2)$. Since the latter is an integer, it follows that 
$\dl^{\m,n}_f(e_2)>\dl^{\m,p}_h(e_2)$.

If there is $e_1\in\gamma$ such that
$\dl^{\m,n}_f(e_1)>\dl^{\m,p}_h(e_1)$, apply the above argument to
$\ol\gamma$ to conclude.
\end{proof}

\begin{thm}\label{vorvor} Let $f_1,f_2\in C^0(G,\Z)$ with 
both $G^\m_{f_1}$ and $G^\m_{f_2}$ connected. Then $\P_{\dl^\m_{f_1}}\cap
  \P_{\dl^\m_{f_2}}\neq\emptyset$ if and only if the
Voronoi cells $\Vor^\m_H(f_1)$ and $\Vor^\m_H(f_2)$
intersect. Furthermore, in this case, 
\begin{equation}\label{ppp}
\P_{\ell,\m,f_1}\cap\P_{\ell,\m,f_2}=
\P_{\ell,\m,f_1}\cap\P_{\dl^\m_{f_2}}=
\P_{\ell,\m,f_2}\cap\P_{\dl^\m_{f_1}}
\end{equation}
and the intersection $\P_{\ell,\m,f_1}\cap\P_{\ell,\m,f_2}$ is 
the closure of the torus orbit of the corresponding face of
$\Vor^\m_H(f_i)$ for $i=1,2$. Also, letting
$G^\m_{f_1,f_2}$ be the graph obtained from the intersection of the
spanning subgraphs $G^\m_{f_1}$ and $G^\m_{f_2}$ of $G$ by keeping
only the edges
$e$ with $\dl^\m_{f_1}(e)=\dl^\m_{f_2}(e)$, there is a natural isomorphism
$$
\prod_{i=1}^d\P_{Z_i}\to \P_{\ell,\m,f_1}\cap\P_{\ell,\m,f_2},
$$
where $Z_1,\dots,Z_d$ are the components of $G^\m_{f_1,f_2}$.
\end{thm}

\begin{proof} Let $X_1,\dots,X_q\subseteq V$ be the level subsets of
  $f_2-f_1$, in increasing order. Let $D_1$ be the coherent acyclic orientation of
the cut subgraph of $G^\m_{f_1}$ induced by the ordered partition
$X_1,\dots,X_q$ of $V$, and $D_2$ that of the cut subgraph of
$G^\m_{f_2}$ 
induced by the same partition in the reverse order,
$X_q,\dots,X_1$. 

If $e=uv\in\E(D_1)$ then 
\begin{equation}\label{D1>}
\frac{f_2(v)-f_2(u)+\m_e}{\ell_e}>
\frac{f_1(v)-f_1(u)+\m_e}{\ell_e}.
\end{equation}
As the right-hand side above is an integer, it follows that
$\dl^\m_{f_2}(e)>\dl^\m_{f_1}(e)$. Conversely, if $e=uv\in\E(G^\m_{f_1})$
satisfies $\dl^\m_{f_2}(e)>\dl^\m_{f_1}(e)$ then \eqref{D1>} holds,
and thus $u\in X_i$ and $v\in X_j$ for $i<j$, that is, $e\in\E(D_1)$. To summarize,
$$
\E(D_1):=\{e\in\E(G^\m_{f_1})\,|\,\dl^\m_{f_2}(e)>\dl^\m_{f_1}(e)\}\quad\text{and}\quad
\E(D_2):=\{e\in\E(G^\m_{f_2})\,|\,\dl^\m_{f_1}(e)>\dl^\m_{f_2}(e)\},
$$
the second equality by analogy with the first. It follows as well 
that $G^\m_{f_1,f_2}$ is the subgraph
of $G^\m_{f_i}$ obtained by removing the edges $e$ in $E(D_i)$
for each $i=1,2$. In addition, the collection of connected
components of $G^\m_{f_i}[X_1],\dots,G^\m_{f_i}[X_q]$ is
$Z_1,\dots,Z_d$ for $i=1,2$.

If $\Vor^\m_H(f_1)$ and $\Vor^\m_H(f_2)$
  intersect then, by \cite[Prop.~5.12]{AE1}, 
$$
\dl^\m_{f_1}+\frac 12 \chi_{\indbi{D}{1}} = \dl^\m_{f_2}+\frac 12 \chi_{\indbi{D}{2}},
$$
where $\chi_{D_i}$ is the characteristic function of $D_i$, taking
value $+1$ at $e\in\E(D_i)$, value $-1$ at $e\in\E$ with $\bar
e\in\E(D_i)$, and value $0$ elsewhere, for $i=1,2$. 
In particular, $|\dl^\m_{f_2}(e)-\dl^\m_{f_1}(e)|\leq 1$ for each
$e\in\E$ with equality
only if $e\in\E(G^\m_{f_1})\cap\E(G^\m_{f_2})$, in fact, if and only
if $e\in\E(D_1)$ and $\ol e\in\E(D_2)$ or $\ol e\in\E(D_1)$ and
$e\in\E(D_2)$. At any rate, it follows that $\P_{\dl^\m_{f_1}}\cap
  \P_{\dl^\m_{f_2}}\neq\emptyset$.

Assume now that $\P_{\dl^\m_{f_1}}\cap
  \P_{\dl^\m_{f_2}}\neq\emptyset$. Then $|\dl^\m_{f_2}(e)-\dl^\m_{f_1}(e)|\leq 1$,
  with equality only if $\dl^\m_{f_1}(e)$ and $\dl^\m_{f_2}(e)$ are integers. Also
  $\dl^\m_{f_2}(e)=\dl^\m_{f_1}(e)$ if neither $\dl^\m_{f_1}(e)$ nor $\dl^\m_{f_2}(e)$
  is an integer. Thus $\P_{\dl^\m_{f_1}}\cap
  \P_{\dl^\m_{f_2}}=\P_\alpha$, where
  $\alpha_e:=1/2(\dl^\m_{f_1}(e)+\dl^\m_{f_2}(e))$ for each $e\in\E(G)$, unless 
$|\dl^\m_{f_2}(e)-\dl^\m_{f_1}(e)|=1/2$, in which case $\alpha_e$ is the half
integer between $\dl^\m_{f_1}(e)$ and $\dl^\m_{f_2}(e)$. 

Observe that 
$G^\m_{f_1,f_2}$ is the spanning subgraph of $G$
obtained by keeping all oriented edges $e$ such
that $\alpha_e\in\Z$. Also,
\begin{equation}\label{af1f2}
\alpha=\dl^\m_{f_1}+\frac{1}{2}\chi_{D_1}=\dl^\m_{f_2}+\frac{1}{2}\chi_{D_2}.
\end{equation}
In particular, by \cite[Prop.~5.12]{AE1}, the Voronoi cells 
$\Vor^\m_H(f_1)$ and $\Vor^\m_H(f_2)$ intersect.

Now,  $\P_{\ell,\m,f_1}\cap\P_{\dl^\m_{f_2}}$ is the subvariety of
$\P_{\ell,\m,f_1}$ given by the equations $\x_{e,\dl^\m_{f_1}(e^{\mathfrak o})}=0$ for
each $e\in\E(D_1)$. It follows from Proposition~\ref{strat0} and
Theorem~\ref{thm:toric} that $\P_{\ell,\m,f_1}\cap\P_{\dl^\m_{f_2}}$ is defined in 
$\P_{\alpha}$ by the following equations: For each $i=1,\dots,d$ and
each oriented cycle $\gamma$ in $Z_i$,
\begin{equation}\label{needlater}
\prod_{e\in \gamma\cap E^{\mathfrak o}}
\x_{e,\alpha_e}\prod_{e\in\bar\gamma\cap E^{\mathfrak o}} \x_{\ol e, \alpha_e} = 
\prod_{e\in\bar\gamma\cap E^{\mathfrak o}} \x_{e, \alpha_e}
\prod_{e\in \gamma\cap E^{\mathfrak o}}
\x_{\ol e, \alpha_e}.
\end{equation}
As the same description applies to
$\P_{\ell,\m,f_2}\cap\P_{\dl^\m_{f_1}}$, we obtain \eqref{ppp} and the
last statement of the theorem. 

Furthermore, by
Theorem~\ref{thm:stratification}, the intersection
$\P_{\ell,\m,f_1}\cap \P_{\ell,\m,f_2}$ is the closure of the torus
orbit of the face $\mathfrak f_i$ of $\Vor^\m_H(f_i)$ corresponding to
$D_i$ for each $i=1,2$. By \cite[Prop.~5.12]{AE1}, the two faces
$\mathfrak f_1$ and $\mathfrak f_2$ coincide with
$\Vor^\m_H(f_1)\cap\Vor^\m_H(f_2)$. The proof of the theorem is finished.
\end{proof}

From Theorem~\ref{vorvor}, we get directly the following corollary.

\begin{cor}\label{Ystr}
 The structure on the standard toric tiling $Y_{\ell, \m}^{\bt}$ is obtained as follows: 
Take the disjoint union 
 $$
\bigsqcup_{\substack{f\in C^0(G,\Z )\\ G^\m_f \mathrm{ connected}}} \,
\P_{\ell,\m,f}.
$$ 
For $f$ and $h$ in $C^{0}(G, \Z)$ with both $G^\m_f$ and $G^\m_h$ connected such that the corresponding Voronoi cells share a face 
 $\f := \Vor^\m_H(f) \cap \Vor^\m_H(h) \neq \emptyset$, identify the closure in $\P_{\ell,\m,f}$ of the torus orbit associated 
 to the face $\f$ of 
 $\Vor^\m_H(f)$  with the closure in $\P_{\ell,\m,h}$ of the torus orbit associated to the face $\f$ of  $\Vor^\m_H(h)$.
 \end{cor}

\subsection{General toric tiling $Y_{\ell, \m}^{a,b}$}\label{gentortil}
We are now finally in position to define the most general form of our toric tilings. 
Let $\ell: E \rightarrow \N$ be an edge length function, and $\m\in
C^1(G, \Z)$ a twisting. 
Let $a: C^{1}(G,\mathbb Z) \to \Gm(\k)=\k^*$ and $b: H^1(G,\Z) \to
\Gm(\k)=\k^*$ be two characters. 
To the quadruple $(\ell, \m, a,b)$ we are going to
 associate the toric tiling $Y_{\ell,\m}^{a,b}$, which is a
 modification of the standard tiling $Y_{\ell, \m}^{\bt}$.
 
To simplify the notation, for each $e\in\E$ we will let
$a_e:=a(\chi_{\indm e}-\chi_{\indmbar{\ol e}})$. Then $a_ea_{\ol e}=1$ for each
$e\in\E$. In addition, we will also denote by $b$ any character of
$C^1(G,\Z)$ extending $b$, and let $b_e:=b(\chi_{\indm e}-\chi_{\indmbar{\ol e}})$ for
each $e\in\E$ as before.

To each $f \in C^0(G, \mathbb Z)$ we
associate the subvariety $\P_{\ell,\m,f}^{a,b}$ of $\P_{\dl^{\m}_f}$ given by the equations  
\[
\forall \textrm{ oriented cycle } \gamma \textrm { in $G^\m_f$},
\]
\[
\prod_{e\in\bar\gamma\cap E^{\mathfrak o}}b_ea_e^{\dl^\m_f(e)}\prod_{e\in \gamma\cap E^{\mathfrak o}} \x_{e,\dl^\m_f(e)} \prod_{e\in\bar\gamma\cap E^{\mathfrak o}}\x_{\ol e,\dl^\m_f(e)}
=\prod_{e\in \gamma\cap E^{\mathfrak o}} b_ea_e^{\dl^\m_f(e)} \prod_{e\in \bar \gamma\cap E^{\mathfrak o}} \x_{e,\dl^\m_f(e)}\prod_{e\in \gamma\cap E^{\mathfrak o}} \x_{\ol e,\dl^\m_f(e)},
\]
where $G^\m_f$ is the spanning subgraph of $G$ consisting 
of all the edges of $G$ for which $\dl^\m_f$ is an integer. 
It is the same $\P_{\ell,\m,f}$ as before,
if $a$ and $b$ are the trivial characters. The equation corresponding
to $\gamma$ may also be written in the format:
\begin{equation}\label{gammapre}
\prod_{e\in\gamma}\x_{e,\mathfrak d^{\mathfrak m}_f(e^{\mathfrak o})}
=\prod_{e\in \gamma} b_ea_e^{\mathfrak d^{\mathfrak m}_f(e^{\mathfrak o})} 
\prod_{e\in \bar \gamma} \x_{e, \mathfrak d^{\mathfrak m}_f(e^{\mathfrak o})},
\end{equation}
where we recall that $e^{\mathfrak o}=e$ if $e\in E^{\mathfrak o}$ and
$e^{\mathfrak o}=\ol e$ otherwise.

We denote by
$Y_{\ell,\m}^{a,b}$ the union of the
$\P_{\ell,\m,f}^{a,b}$ for those $f\in C^0(G,\Z)$ for which $G^\m_f$ is
connected. 

The statements in Theorem~\ref{vorvor} applies to the
$\P_{\ell,\m,f}^{a,b}$, with the same proof given there, with the only
obvious modification in Equation~\eqref{needlater}.  As a consequence we
obtain the analogous to Corollary~\ref{Ystr}. 

\begin{thm} The structure on $Y_{\ell, \m}^{a,b}$ is given as follows: Take the disjoint union 
$$
\bigsqcup_{\substack{f\in C^0(G,\Z )\\ G^\m_f \mathrm{ connected }}}
\, \P_{\ell,\m,f}^{a,b}.
$$ 
For $f$ and $h$ in $C^{0}(G, \Z)$ with both $G^\m_f$ and $G^\m_h$ 
connected such that the corresponding Voronoi cells share a face 
 $\f := \Vor^\m_H(f) \cap \Vor^\m_H(h) \neq \emptyset$, identify the 
closure in $\P_{\ell,\m,f}^{a,b}$ of the torus orbit associated 
to the face $\f$ of $\Vor^\m_H(f)$  with the closure in
$\P_{\ell,\m,h}^{a,b}$ of the torus orbit associated to the face 
$\f$ of  $\Vor^\m_H(h)$. 
\end{thm}

The restriction that $G^\m_f$ be connected for $\P^{a,b}_{\ell,\m,f}$
to be part of $Y_{\ell, \m}^{a,b}$ is in fact not necessary,
as the next proposition shows. 

\begin{prop}\label{propnotcon} 
For each $f\in C^0(G,\Z)$ there is $h\in C^0(G,\Z)$ with
  $G^\m_h$ connected such that $\P^{a,b}_{\ell,\m,f}=\P^{a,b}_{\ell,\m,h}\cap\P_{\dl^\m_f}$.
\end{prop}

\begin{proof} First we claim that there is $h\in C^0(G,\Z)$ such that 
$G^\m_h$ is connected, contains $G^\m_f$ as a subgraph, 
and $|\dl^\m_h(e)-\dl^\m_f(e)|\leq 1/2$ for each
$e\in\E(G)$. And then we will see that this $h$ works.

Let $G_1,\dots,G_n$ be the connected components of
  $G^\m_f$. Let $V_i:=V(G_i)$ for each $i=1,\dots,n$. Up to
  reordering, assume by induction that there is $h\in C^0(G,\Z)$ such 
that 
\begin{equation}\label{dhf}
|\dl^\m_{h}(e)-\dl^\m_f(e)|\leq 1/2\quad\text{for each }e\in\E(G),
\end{equation}
such that $G^\m_f$ is a subgraph of $G^\m_h$, and that
$G_1,\dots,G_i$ are contained in the same connected component of
$G^\m_h$. (For $i=1$, simply take $h:=f$.) Set $S_i:=V_1\cup\cdots\cup
V_i$ and $T_i:=V_{i+1}\cup\cdots\cup V_n$. Let $q$ 
be the smallest nonnegative integer such that
$g:=h+q\chi_{\indbi{T}{i}}$ satisfies $\dl^\m_g(e)\in\Z$ for some
oriented edge $e\in\E(S_i,T_i)$. Let $e'$ be such an
edge. Its target belongs to $V_j$ for some $j>i$. We may assume
$j=i+1$. Clearly, $\dl^\m_{g}(e)=\dl^\m_h(e)$ for every $e\in\E(G)$
connecting two vertices either both inside $S_i$ or outside. Thus
$G^\m_f$ is a subgraph of $G^\m_g$ and $G_1,\dots,G_i$ are contained in
the same connected component of $G^\m_g$. Furthermore, since
$e'\in\E(G^\m_g)$, also $G_{i+1}$ is in the same component. 

We will now prove that 
\begin{equation}\label{dhg}
|\dl^\m_g(e)-\dl^\m_f(e)|\leq 1/2
\end{equation}
for each $e\in\E(G)$. As observed before, $\dl^\m_{g}(e)=\dl^\m_h(e)$, and thus 
\eqref{dhg} holds for all $e\in\E(G)$ by \eqref{dhf}, unless 
$e$ or $\ol e$ belongs to $\E(S_i,T_i)$. Let $e\in\E(S_i,T_i)$. Then 
$\dl^\m_f(e)$ is not an integer. If $\dl^\m_h(e)\in\Z$ then $q=0$, and
thus $\dl^\m_{g}(e)=\dl^\m_h(e)$. If $\dl^\m_h(e)\not\in\Z$ then 
$\dl^\m_h(e)=\dl^\m_f(e)$ by \eqref{dhf} and
$\dl^\m_h(e)\leq\dl^\m_g(e)\leq\dl^\m_h(e)+1/2$. In any case, 
we conclude that \eqref{dhg} now holds for every $e\in\E(G)$.

By induction, the claim is proved. It follows immediately from it that 
$\P_{\dl^\m_f}\subseteq\P_{\dl^m_h}$. Furthermore, since $G^\m_h$
contains $G^\m_f$ as a subgraph, which implies that every oriented cycle
$\gamma$ in $G^\m_f$ is an oriented cycle of $G^\m_h$ satisfying 
$\dl^\m_h(e)=\dl^\m_f(e)$ for every $e\in\gamma$, all the equations
defining $\P^{a,b}_{\ell,\m,f}$ in $\P_{\dl^\m_f}$ are equations
defining $\P^{a,b}_{\ell,\m,h}$ in $\P_{\dl^\m_h}$, whence 
$\P^{a,b}_{\ell,\m,f}\supseteq \P^{a,b}_{\ell,\m,h}\cap \P_{\dl^\m_f}$.

Let now $\gamma$ be an oriented cycle in $G^\m_h$, and consider the
corresponding equation \eqref{gammapre} among those defining $\P^{a,b}_{\ell,\m,h}$.
It remains to show that the equation is satisfied on $\P^{a,b}_{\ell,\m,f}$. 

If $\dl^{\mathfrak m}_f(e)\in\Z$ for every $e\in\gamma$ then $\gamma$
is in $G^\m_f$. Since also $\mathfrak d^{\mathfrak m}_f(e)=\mathfrak d^{\mathfrak
  m}_h(e)$ for each $e\in\gamma$ by \eqref{dhf}, it follows that the
equations corresponding to $\gamma$ for $\P^{a,b}_{\ell,\m,f}$ and
for $\P^{a,b}_{\ell,\m,h}$ are the same, and thus \eqref{gammapre} is
satisfied on $\P^{a,b}_{\ell,\m,f}$.

On the other hand, if there is $e\in\gamma$ such that 
$\mathfrak  d^{\mathfrak m}_f(e)\not\in\Z$, then Lemma~\ref{e1e2}
yields $e_1,e_2\in\gamma$ such that 
$$
\mathfrak d^{\mathfrak m}_f(e_1)=\mathfrak d^{\mathfrak
  m}_h(e_1)-1/2\quad\text{and}\quad
\mathfrak d^{\mathfrak m}_f(e_2)=\mathfrak d^{\mathfrak m}_h(e_2)+1/2.
$$
But then $\x_{\ol e_1, \mathfrak d^{\mathfrak m}_h(e_1^{\mathfrak o})}=0$ 
(resp.~$\x_{e_2, \mathfrak d^{\mathfrak m}_h(e_2^{\mathfrak o})}=0$) on
$\P^{a,b}_{\ell,\m,f}$, which implies that the right-hand side
(resp.~left-hand side) of
\eqref{gammapre} vanishes on $\P^{a,b}_{\ell,\m,f}$.  Thus
\eqref{gammapre} is satisfied on $\P^{a,b}_{\ell,\m,f}$.
\end{proof}

 \begin{prop}
 The toric tilings $Y_{\ell, \m}^{a,b}$ and $Y_{\ell,\m}^{\bt}$ have the same combinatorial support in $\mathbf R$, in the sense that they have the same set of supporting $\P_{\dl^\m_f}$. Moreover,  they are rationally equivalent. 
 \end{prop}

\begin{proof}
The claim about the set of supporting $\P_{\dl^\m_f}$ is obvious from the
definition. The second statement is clear as in each $\P_{\dl^\m_f}$
there is a rational deformation of $\P^{a,b}_{\ell,\m,f}$ to
$\P_{\ell,\m,f}$ obtained by deforming the characters $a$ and $b$ to
the trivial characters.
\end{proof}

\section{Mixed toric tilings II: Orbit structure}\label{mixorb}
In this section we describe actions of ${\mathbf G}_{\mathbf m}^{|V|-1}$ on
$\mathbf R$ and on the arrangements of toric varieties $Y_{\ell,
  \m}^{a,b}$, and provide a complete description of the orbits.

\subsection{The action of ${\mathbf G}_{\mathbf m}^{|V|-1}$}\label{action} 
There is a natural action of the group of characters of $C^1(G,\Z)$ on
$\mathbf R$. Namely, each character $a\: C^1(G,\Z)\to \Gm(\k)=\k^*$ acts on 
the point $p$ of $\P_\alpha=\prod_{e\in E^{\mathfrak
    o}}\P^1_{e,\alpha_e}$ with coordinates $(x_{e,\alpha_e}:x_{\ol
  e,\alpha_e})$, for each $\alpha\in C^1(G,\Z)$, by taking it to the
point on the same $\P_\alpha$ with coordinates $(a_e x_{e,\alpha_e}:x_{\ol
  e,\alpha_e})$. Clearly, the action takes $\P_\alpha$ to itself for
each $\alpha\in C^1(G,\frac{1}{2}\Z)$. 

The boundary map $d^*\: C^1(G,\Z)\to C^0(G,\Z)$ induces thus a natural
action of the group of characters of $C^0(G,\Z)$ on $\mathbf R$. Given
a character $c\: C^0(G,\Z)\to \Gm(\k)$, since $d^*(\chi_{\indm e}-\chi_{\indmbar{\ol
  e}})=\chi_{\indm v}-\chi_{\indm u}$ for each $e=uv\in\E$, we have that $c$ acts on the
point $p$ of $\P_\alpha=\prod_{e\in E^{\mathfrak
    o}}\P^1_{e,\alpha_e}$ with coordinates $(x_{e,\alpha_e}:x_{\ol
  e,\alpha_e})$, for each $\alpha\in C^1(G,\Z)$, by taking it to the
point on the same $\P_\alpha$ with coordinates $(c_v x_{e,\alpha_e}: c_ux_{\ol
  e,\alpha_e})$. Here, as before, $c_v:=c(\chi_{\indm v})$ for each $v\in V$. 

Clearly, if $c_v=c_u$ for all $u,v\in V$, the action of $c$ is
trivial. Thus we need only consider characters induced from those 
of the quotient $C^0(G,\Z)/\Z\chi_{\ind V}$. The latter has character group
isomorphic to ${\mathbf G}_{\mathbf m}^{|V|-1}$. We will loosely refer to this
action as the action of ${\mathbf G}_{\mathbf m}^{|V|-1}$ on $\mathbf R$. 

Finally, it follows from the equations defining the subvariety
$\P_{\ell,\m,f}^{a,b}$ of $\P_{\dl^{\m}_f}$ for each $f\in C^0(G,\Z)$ that the action of
${\mathbf G}_{\mathbf m}^{|V|-1}$ on $\mathbf R$ 
restricts to an action on $\P_{\ell,\m,f}^{a,b}$, and thus induces an
action on $Y_{\ell,
  \m}^{a,b}$ for any choices of edge length function $\ell: E
\rightarrow \N$, twisting $\m\in
C^1(G, \Z)$ and characters $a: C^{1}(G,\mathbb Z) \to \k^*$ and $b:
H^1(G,\Z) \to \k^*$.

Moreover, the action on the open locus of $\P_{\ell,\m,f}^{a,b}$ where
all the $e$-th coordinates for $e\in \E(G^{\mathfrak m}_f)\cap
E^{\mathfrak o}$ are nonzero is transitive, if $G^{\mathfrak m}_f$ is connected. Indeed, setting
$x_e:=\x_{e,\dl^\m_f(e)}/\x_{\ol e,\dl^\m_f(e)}$ for each such $e$, the
equations defining $\P_{\ell,\m,f}^{a,b}$ on the open locus become
$$
\prod_{e\in\bar\gamma\cap E^{\mathfrak o}}b_ea_e^{\dl^\m_f(e)}\prod_{e\in \gamma\cap E^{\mathfrak o}} x_e
=\prod_{e\in \gamma\cap E^{\mathfrak o}} b_ea_e^{\dl^\m_f(e)} \prod_{e\in \bar \gamma\cap E^{\mathfrak o}} x_{e}
$$
for all oriented cycles $\gamma$ in $G^\m_f$. For each point $p$ on
the open locus, we may use the $x_e$ to
define the character $y\: C^1(G^\m_f,\Z)\to  \k^*$ for which, following
notation, $y_e:=x_e(p)b_e^{-1}a_e^{-\dl_f^\m(e)}$ for each $e\in
\E(G^\m_f)\cap E^{\mathfrak o}$. The equations for the $x_e$
correspond thus to requiring that $y$ restricts to the trivial
character on $H^1(G^\m_f,\Z)$. Since $d^*\: C^1(G^\m_f,\Z)\to C^0(G^\m_f,\Z)$ has
kernel $H^1(G^\m_f,\Z)$, it follows that $y$ is induced from a
character of the image of $d^*$, which, since $G^{\mathfrak m}_f$ is
connected, 
is the subgroup of all $h\in
C^0(G^\m_f,\Z)$ with $\sum h(v)=0$. Thus the character of the image of $d^*$ extends
to a character of $C^0(G^\m_f,\Z)=C^0(G,\Z)$, which we denote
$c\:C^0(G,\Z)\to \k^*$. It follows now that $p$ is the result of the
action of $c$ on the point $q$ on the open locus of $\P_{\ell,\m,f}^{a,b}$
satisfying $x_e(q)=b_ea_e^{\dl_f^\m(e)}$ for each $e\in \E(G^{\mathfrak m}_f)\cap
E^{\mathfrak o}$.

 
\subsection{$Y_{\ell, \m}^{a,b}$ as a union of
  orbits} \label{sec:uorbits} In this section we describe $Y_{\ell,
  \m}^{a,b}$ as a union of orbits under the action of $\mathbf
G_{\mathbf m}^{|V|-1}$. We need first two preliminary results.

\begin{prop}\label{f1f2} Let $f_1,f_2\in C^0(G,\Z)$ and $n\in\N$. 
If $\mathfrak
  d^{\mathfrak m,n}_{f_2}-\mathfrak
  d^{\mathfrak m,n}_{f_1}\in H^1(G,\Z)$ then $\mathfrak
  d^{\mathfrak m,n}_{f_2}=\mathfrak
  d^{\mathfrak m,n}_{f_1}$.
\end{prop}

\begin{proof} Since $\mathfrak d^{\mathfrak m,n}_{f_2}-\mathfrak
  d^{\mathfrak m,n}_{f_1}\in C^1(G,\Z)$, it follows that $\mathfrak
  d^{\mathfrak m,n}_{f_2}(e)\in\Z$ if and only if $\mathfrak
  d^{\mathfrak m,n}_{f_1}(e)\in\Z$, for each $e\in\E$. Then
$\mathfrak
  d^{\mathfrak m,n}_{f_2}(e)-\mathfrak
  d^{\mathfrak m,n}_{f_1}(e)=\delta_e^{\mathfrak
    m,n}(f_2)-\delta_e^{\mathfrak m,n}(f_1)$.

Let $h:=f_2-f_1$ and $\gamma:=\mathfrak
  d^{\mathfrak m,n}_{f_2}-\mathfrak
  d^{\mathfrak m,n}_{f_1}$. Suppose by contradiction that $\gamma\neq
  0$. Let $G'$ be
  the subgraph of $G$ obtained by keeping only the edges in the
  support of $\gamma$. Then
  there is a connected component $G''$ of $G'$ containing an edge. Let
  $v$ be a vertex of $G''$ where $h(v)$ is maximum. 
Since $\gamma\in H^1(G,\Z)$, we have
\begin{equation}\label{sum0}
\sum_{\substack{e\in\E(G'')\\ \he_e=v}}\gamma_e=\sum_{\substack{e\in\E(G)\\
      \he_e=v}}\gamma_e=0.
\end{equation}

On the other hand, for each $e=uv\in\E(G'')$ write
$$
f_i(v)-f_i(u)+n\mathfrak m_e=\delta_e^{\mathfrak
  m,n}(f_i)n\ell_e+\rho_i(e),
$$
where $0\leq\rho_i(e)<n\ell_e$ for $i=1,2$. Then
$$
h(v)-h(u)+\rho_1(e)=(\delta_e^{\mathfrak
    m,n}(f_2)-\delta_e^{\mathfrak m,n}(f_1))n\ell_e+\rho_2(e).
$$
Hence, since $h(v)\geq h(u)$, we have that 
$$
\gamma_e=\delta_e^{\mathfrak m,n}(f_2)-\delta_e^{\mathfrak m,n}(f_1)\geq
0\quad\text{for each }e=uv\in\E(G'').
$$
But then Equation~\eqref{sum0} yields that $\gamma_e=0$ for each
$e\in\E(G'')$ with $\he_e=v$, an absurd.
\end{proof}

\begin{lemma}\label{fbg} Let $n\in\N$ and $f\in C^0(G,\Z)$. Let $g\in
  C^0(G,\Z)$ be defined by $g(v)=\lfloor f(v)/n\rfloor$ for each $v\in
  V(G)$. Let $e=uv\in\E(G)$. Then the following statements hold:
\begin{enumerate}
\item If $\mathfrak d^{\mathfrak m,n}_f(e)\in\Z$ then 
$\mathfrak d^{\mathfrak m}_g(e)=\mathfrak d^{\mathfrak m,n}_f(e)$.
\item If $\mathfrak d^{\mathfrak m,n}_f(e)\not\in\Z$ then 
$$
|\mathfrak d^{\mathfrak m}_g(e)-\mathfrak d^{\mathfrak m,n}_f(e)|\leq\frac{1}{2}.
$$
\end{enumerate}
\end{lemma}

\begin{proof} Observe first that 
\begin{equation}\label{fgf}
\frac{f(v)-f(u)+n\mathfrak
  m_e}{n\ell_e}-\frac{1}{\ell_e}<\frac{g(v)-g(u)+\mathfrak
  m_e}{\ell_e}<\frac{f(v)-f(u)+n\mathfrak
  m_e}{n\ell_e}+\frac{1}{\ell_e}
\end{equation}
for each  $e=uv\in\E$. Observe as well that the middle
number above is in $(1/\ell_e)\Z$. 

If 
$\mathfrak d^{\mathfrak m,n}_f(e)\in\Z$ then 
$$
\frac{f(v)-f(u)+n\mathfrak  m_e}{n\ell_e}\in\Z
$$
and thus it follows from Inequalities \eqref{fgf} and the
observation thereafter that
$$
\frac{g(v)-g(u)+\mathfrak
  m_e}{\ell_e}=\frac{f(v)-f(u)+n\mathfrak  m_e}{n\ell_e},
$$
and thus $\mathfrak d^{\mathfrak m}_g(e)=\mathfrak d^{\mathfrak m,n}_f(e)$.

On the other, suppose $\mathfrak d^{\mathfrak m,n}_f(e)\not\in\Z$. If
$\mathfrak d^{\mathfrak m}_g(e)\not\in\Z$ either then 
$\mathfrak d^{\mathfrak m,n}_f(e)=\mathfrak d^{\mathfrak
  m}_g(e)$. Indeed, if this is not the case there is an integer
$c_e$ satisfying either the inequalities
$$
\frac{f(v)-f(u)+n\mathfrak
  m_e}{n\ell_e}<c_e<\frac{g(v)-g(u)+\mathfrak
  m_e}{\ell_e}
$$
or the reverse inequalities. If the displayed inequalities hold then
$$
c_e<\frac{g(v)-g(u)+\mathfrak
  m_e}{\ell_e}< \frac{f(v)-f(u)+n\mathfrak
  m_e}{n\ell_e}+\frac{1}{\ell_e} < c_e + \frac{1}{\ell_e},
$$
an absurd, as the second number above is in $(1/\ell_e)\Z$. The
reverse inequalities yield a similar contradiction.

Finally, suppose $\mathfrak d^{\mathfrak m,n}_f(e)\not\in\Z$ but 
$\mathfrak d^{\mathfrak m}_g(e)\in\Z$. Then Inequalities \eqref{fgf}
are equivalent to
$$
\mathfrak d^{\mathfrak m}_g(e)-\frac{1}{\ell_e}<\frac{f(v)-f(u)+n\mathfrak
  m_e}{n\ell_e}<\mathfrak d^{\mathfrak m}_g(e)+\frac{1}{\ell_e},
$$
from which follows that $\mathfrak d^{\mathfrak m,n}_f(e)=\mathfrak
d^{\mathfrak m}_g(e)-1/2$ if the middle term above is smaller than $\mathfrak
d^{\mathfrak m}_g(e)$, or $\mathfrak d^{\mathfrak m,n}_f(e)=\mathfrak
d^{\mathfrak m}_g(e)+1/2$ otherwise. In any case, 
$$
\mathfrak d^{\mathfrak m,n}_f(e)=\mathfrak d^{\mathfrak
  m}_g(e)\pm\frac{1}{2}.
$$
\end{proof}
 
\begin{thm}\label{unionorbits} 
Let $\mathfrak o$ be an orientation for the edges of $G$. Let $\ell\: E\to\mathbb N$ be a length function,
  $\mathfrak m\in C^1(G,\Z)$, and let
$$
a\: C^1(G,\Z)\to \k^*\quad\text{and}\quad
b\: C^1(G,\Z)\to \k^*
$$
be characters. For each $n\in\mathbb N$ and $f\in C^0(G,\Z)$,
let $p^n_{f}$ be the point on $\mathbf P_{\mathfrak
  d^{\mathfrak m,n}_f}$ given by the coordinates 
$$
(b_ea_e^{\mathfrak d^{\mathfrak m,n}_f(e)}:1) \quad\text{for each }
e\in E^{\mathfrak o}\text{ with }\mathfrak d^{\mathfrak
  m,n}_f(e)\in\Z.
$$
Then $Y_{\ell,\m}^{a,b}$
is the union of the orbits of the $p^n_f$ under the action of 
$\mathbf G_{\mathbf m}^{|V|-1}$.
\end{thm} 

(We are also denoting by $b$ its restriction to $H^1(G,\Z)$.)

\begin{proof} As defined
in Section~\ref{gentortil}, 
the toric arrangement $Y_{\ell,\m}^{a,b}$ is the union of the toric subvarieties 
$\mathbf P^{a,b}_{\ell,\m,g}$ of $\mathbf P_{\mathfrak d^{\mathfrak m}_g}$
as $g$ varies in $C^0(G,\Z)$. (By Proposition~\ref{propnotcon}, we do
not need to restrict to those $g$ for which $G^\m_g$ is connected.) 
Each subvariety $\mathbf P^{a,b}_{\ell,\m,g}$ is given by the
equations:
\[
\forall \textrm{ oriented cycle } \gamma \textrm { in $G^\m_g$},
\]
\begin{equation}\label{gamma}
\prod_{e\in\gamma}\x_{e,\mathfrak d^{\mathfrak m}_g(e^{\mathfrak o})}
=\prod_{e\in \gamma} b_ea_e^{\mathfrak d^{\mathfrak m}_g(e^{\mathfrak o})} 
\prod_{e\in \bar \gamma} \x_{e, \mathfrak d^{\mathfrak m}_g(e^{\mathfrak o})}.
\end{equation}

Let $f\in C^0(G,\Z)$ and $n\in\N$. Let $g\in C^0(G,\Z)$
satisfying $g(v):=\lfloor f(v)/n\rfloor$ for each $v\in
V$. Lemma~\ref{fbg} yields 
$\mathbf P_{\mathfrak
  d^{\mathfrak m,n}_f}\subseteq \mathbf P_{\mathfrak
  d^{\mathfrak m}_g}$. In particular, $p^n_f\in \mathbf P_{\mathfrak
  d^{\mathfrak m}_g}$. We will now show that $p^n_f\in \mathbf P^{a,b}_{\ell,\m,g}$.

Let $\gamma$ be an oriented cycle in
$G^{\mathfrak m}_g$. We show that Equation~\eqref{gamma} 
is satisfied on $p^n_f$. First, if $\mathfrak  d^{\mathfrak m,n}_f(e)\in\Z$ for each $e\in\gamma$,
then $\mathfrak d^{\mathfrak m,n}_f(e)=\mathfrak d^{\mathfrak
  m}_g(e)$ for each $e\in\gamma\cup\bar\gamma$ by Lemma~\ref{fbg}. It 
follows that $(\x_{e, \mathfrak d^{\mathfrak m}_g(e)}: \x_{\ol e, \mathfrak
  d^{\mathfrak m}_g(e)})=(b_ea_e^{\mathfrak d^{\mathfrak m}_g(e)}:1)$
on $p^n_f$ for each $e\in E^{\mathfrak o}$ with
$e\in\gamma\cup\bar\gamma$, and then Equation~\eqref{gamma} is satisfied on
$p^n_f$. 

Suppose now there is $e\in\gamma$ such that 
$\mathfrak  d^{\mathfrak m,n}_f(e)\not\in\Z$. By Lemma~\ref{fbg} and
Lemma~\ref{e1e2}, there are $e_1,e_2\in\gamma$ such that 
$$
\mathfrak d^{\mathfrak m,n}_f(e_1)=\mathfrak d^{\mathfrak
  m}_g(e_1)-1/2\quad\text{and}\quad 
\mathfrak d^{\mathfrak m,n}_f(e_2)=\mathfrak d^{\mathfrak
  m}_g(e_2)+1/2. 
$$
But then $\x_{\ol e_1, \mathfrak d^{\mathfrak m}_g(e^{\mathfrak o}_1)}=0$
(resp.~$\x_{e_2, \mathfrak d^{\mathfrak
    m}_g(e^{\mathfrak o}_2)}=0$) on
$p^n_f$, which implies that the right-hand side (resp.~left-hand side)
of Equation~\eqref{gamma} vanishes on 
$p^n_f$. Since both sides vanish, the equation is satisfied.

\smallskip

Let now $p$ be a point on $\mathbf P^{a,b}_{\ell,\m,g}$
for some $g\in C^0(G,\Z)$. We will show that $p$ is on the orbit of
$p^n_f$ for certain $n\in\N$ and $f\in C^0(G,\Z)$. This will
finish the proof. 

\smallskip

First, let $\E(p)\subseteq \E$ be
  the set of $e\in\E$ such that $\mathfrak
  d^{\mathfrak m}_g(e)\in\Z$ and $\x_{e,\mathfrak
  d^{\mathfrak m}_g(e^{\mathfrak o})}(p)=0$. Let $E(p)\subseteq E(G)$
be the support of $\E(p)$. Notice that if $e\in\E(p)$,
then $\x_{\ol e,\mathfrak
  d^{\mathfrak m}_g(e^{\mathfrak o})}(p)\neq 0$, and thus $\bar
e\not\in\E(p)$. Furthermore, since $p$ satisfies Equation~\eqref{gamma}
for each oriented cycle $\gamma$ in $G^{\mathfrak m}_g$, for
each such cycle, if $e\in\gamma\cap \E(p)$, then there is
$e'\in\bar\gamma\cap \E(p)$. It follows that there is no oriented
cycle in $\E(p)$. In short, $\E(p)$ is an acyclic orientation of $E(p)$.

We claim there are $n\in\N$ and $f\in C^0(G,\Z)$ such that 
$\mathfrak d^{\mathfrak m,n}_f(e)=\mathfrak d^{\mathfrak m}_g(e)$ for
each $e\in\E$ supported away from $E(p)$ and $\mathfrak d^{\mathfrak m,n}_f(e)=\mathfrak
d^{\mathfrak m}_g(e)+1/2$ for each $e\in\E(p)$. In particular,
$G^{\mathfrak m,n}_f$, the subgraph of $G$ obtained by removing all
edges $e\in\E$ for which $\mathfrak d^{\mathfrak m,n}_f(e)\not\in\Z$, 
is also obtained from $G^{\mathfrak m}_g$ by removing
all edges of $E(p)$.

Indeed, we will let $\rho\in C^0(G,\Z)$ be positive and small enough, and let $n\in\N$ be large enough such that $f\in C^0(G,\Z)$ with
$f(v)=ng(v)+\rho(v)$ for each $v\in V$ satisfies the desired
conditions. First, if $(1/n)\rho$ is small enough we may assume
that $\mathfrak d^{\mathfrak m,n}_f(e)=\mathfrak d^{\mathfrak m}_g(e)$
whenever $\mathfrak d^{\mathfrak m}_g(e)\not\in\Z$. Second, we set
$\rho$ to be constant on each connected component of the subgraph obtained from
$G^{\mathfrak m}_g$ by removing all edges in $E(p)$. Furthermore, the
complementary subgraph, 
obtained by keeping
only the edges in $E(p)$, has $\E(p)$ as orientation, and
since there is no oriented cycle in $\E(p)$, we may choose $\rho$ such that $\rho(v)>\rho(u)$ for
each $e=uv\in\E(p)$. If $(1/n)\rho$ is small enough, it follows that 
$\mathfrak d^{\mathfrak m,n}_f(e)=\mathfrak d^{\mathfrak m}_g(e)+1/2$
for each $e\in\E(p)$. And if $e=uv\in\E$ is supported away from $E(p)$
and $\mathfrak d^{\mathfrak m}_g(e)\in\Z$, then $e$ is an edge of
$G^{\mathfrak m}_g$ supported away from $E(p)$, whence
$\rho(v)=\rho(u)$ and so $\mathfrak d^{\mathfrak m,n}_f(e)=\mathfrak
d^{\mathfrak m}_g(e)$. The proof of the claim is finished.

\smallskip

Let $n$ and $f$ be as in the claim. Since $|\mathfrak d^{\mathfrak m,n}_f(e)-\mathfrak
d^{\mathfrak m}_g(e)|\leq 1/2$, with equality only if $\mathfrak
d^{\mathfrak m}_g(e)\in\Z$, for each $e\in\E$, we have that $\mathbf P_{\mathfrak
  d^{\mathfrak m,n}_f}\subseteq \mathbf P_{\mathfrak
  d^{\mathfrak m}_g}$, and thus $p^n_f\in \mathbf P_{\mathfrak
  d^{\mathfrak m}_g}$. 

\smallskip

We claim that an $e$-th
coordinate vanishes on $p$ if and only if it vanishes on
$p^n_f$, for each $e\in\E(G^{\mathfrak m}_g)\cap E^{\mathfrak o}$. Indeed, suppose
first that $e\in\E(p)$. Then $\x_{e,\mathfrak d^{\mathfrak
    m}_g(e)}(p)=0$. Now, $\mathfrak d^{\mathfrak
  m,n}_f(e)=\mathfrak d^{\mathfrak m}_g(e)+1/2$. Thus, $\x_{e,\mathfrak d^{\mathfrak
    m}_g(e)}(p^n_f)=0$ as well. And if $\ol e\in\E(p)$, then
$\x_{\ol e,\mathfrak d^{\mathfrak
    m}_g(e)}(p)=0$. Since $\mathfrak d^{\mathfrak
  m,n}_f(e)=\mathfrak d^{\mathfrak m}_g(e)-1/2$, it follows that
$\x_{\ol e,\mathfrak d^{\mathfrak
    m}_g(e)}(p^n_f)=0$ as well. Finally, suppose $e$ is not supported in $E(p)$, that is, $e\in
\E(G^{\mathfrak m,n}_f)$. Then the $e$-th coordinates of $p^n_f$ are
nonzero. So are those of $p$: since
$e\not\in\E(p)$, we have $\x_{e,\mathfrak d^{\mathfrak
    m}_g(e)}(p)\neq 0$, and since $\ol e\not\in\E(p)$, we
have $\x_{\ol e,\mathfrak d^{\mathfrak
    m}_g(e)}(p)\neq 0$. 

\smallskip

Let $G'$ be a connected component of $G^{\mathfrak m,n}_f$. By what we
saw above, all the $e$-th coordinates of $p$ and
$p^n_f$ are nonzero for $e\in\E(G')\cap E^{\mathfrak o}$. 
In addition, for each oriented cycle $\gamma$ in
$G'$, Equation~\eqref{gamma} is satisfied on $p$ and on $p^n_f$. But
then the projections of $p$ and $p^n_f$ on 
the product
$$
\prod_{e\in\E(G')\cap E^{\mathfrak o}}\P^1_{e,\mathfrak
  d^{\mathfrak m,n}_f(e)}
$$
lie on the open torus inside the toric subvariety given by
Equations~\eqref{gamma} for oriented cycles 
$\gamma$ in $G'$. These two projections lie on the same orbit of the
product by the action of ${\mathbf G}_{\mathbf m}^{|V(G')|-1}$, as
observed in Subsection~\ref{action}.

As the above holds for each connected component of $G^{\mathfrak
  m,n}_f$, it follows that $p$ is on the orbit of $p^n_f$ by the
action of ${\mathbf G}_{\mathbf m}^{|V|-1}$. 
\end{proof}

\begin{cor} Given $\alpha\in C^1(G,\frac 12 \Z)$, we have that
  $\P_\alpha^*\cap Y^{a,b}_{\ell,\m}\neq\emptyset$ if and only if
  $\alpha=\dl^{\m,n}_f$ for certain $n\in\mathbb N$ and $f\in C^0(G,\Z)$, and in this case
  $\P_\alpha^*\cap Y^{a,b}_{\ell,\m}$ is an orbit of the action of
  $\mathbf G_{\mathbf m}^{|V|-1}$ on $Y^{a,b}_{\ell,\m}$ which is
  dense in $\P_\alpha\cap Y^{a,b}_{\ell,\m}$.
\end{cor}

\begin{proof} If $\alpha=\dl^{\m,n}_f$ for certain $n$ and $f$ then
  $p^n_f\in \P_\alpha^*\cap Y^{a,b}_{\ell,\m}$ by
  Theorem~\ref{unionorbits}. Conversely, if $p\in \P_\alpha^*\cap
  Y^{a,b}_{\ell,\m}$ then $p$ is on the orbit of $p^n_f$ for certain
  $n$ and $f$ by Theorem~\ref{unionorbits}. Then 
$\P_\alpha^*\cap \P_{\dl^{\m,n}_f}^*\neq\emptyset$ and thus
$\alpha=\dl^{\m,n}_f$. Notice that $n$ and $f$ are not unique, but
$\dl^{\m,n}_f$ is and it determines $p^n_f$. It follows that
$\P_\alpha^*\cap Y^{a,b}_{\ell,\m}$ is the orbit of $p^n_f$.

It remains to show that $\P_\alpha^*\cap Y^{a,b}_{\ell,\m}$ is dense
in $\P_\alpha\cap Y^{a,b}_{\ell,\m}$. Consider a point of
$\P_\beta^*\cap Y^{a,b}_{\ell,\m}$ for $\beta\in C^1(G,\frac 12 \Z)$
such that $\P_\beta\subseteq\P_\alpha$. As we have seen above, the
point is on the orbit of $p^q_h$ for certain $q\in\mathbb N$ and $h\in
C^1(G,\Z)$ such that $\beta=\dl^{\m,q}_h$. Let $G_\alpha$
(resp.~$G_\beta$) be the spanning subgraph of $G$ containing the
supports of all the oriented edges $e\in\E$ for which $\alpha_e\in\Z$
(resp.~$\beta_e\in\Z$). Since
$\P_\beta\subseteq\P_\alpha$, we have that $G_\beta$ is a subgraph of
$G_\alpha$. 

Give an orientation $D$ to the set of edges in $G_\alpha$
which are not in $G_\beta$ by saying the $e\in\E(D)$ if
$\beta(e)<\alpha(e)$. It follows from Lemma~\ref{e1e2} that the
orientation is acyclic. Furthermore, let $X_0,\dots,X_m$ be the connected
components of $G_\beta$. Since $\alpha_e=\beta_e$ for each
$e\in\E(X_i)$ for each $i$, it follows from Lemma~\ref{e1e2} that we may
assume the $X_i$ are ordered in such a way that if $e=uv\in\E(D)$
connects $u\in X_i$ and $v\in X_j$ then $i>j$. 

Consider now the character $c\colon C^0(G,\Z)\to\mathbf G_{\mathbf
  m}(\kappa(t))$ that takes $v\in V(X_i)$ to $t^i$ for each
$i=0,\dots,m$. Let it act on $p^n_f$ for each nonzero $t\in\kappa$ and
let $t$ tend to 0. The limit is $p^q_h$. Thus $p^q_h$, and hence its
orbit, is in the
closure of the orbit of $p^n_f$. It follows that $\P_\beta^*\cap
Y^{a,b}_{\ell,\m}$ is in the closure of $\P_\alpha^*\cap
Y^{a,b}_{\ell,\m}$ for each $\beta\in C^1(G,\frac 12 \Z)$
such that $\P_\beta\subseteq\P_\alpha$, and hence $\P_\alpha^*\cap Y^{a,b}_{\ell,\m}$ is dense
in $\P_\alpha\cap Y^{a,b}_{\ell,\m}$.
\end{proof}

\section{Mixed toric tilings III: Equations}\label{sec:eqtiling}

In this section, we work out the equations for the toric tiling $Y_{\ell, \m}^{a,b}$ in $\mathbf R$.  For this, we will first give ``equations'' to the toric tiling
$\mathbf R$. The tiling $Y_{\ell,\m}^{a,b}$ 
will then be described by a 
further set of equations, 
leading to the
determination
of all the  
supporting ``hypercubes'' $\P_\alpha$, and restricting to the equations determining the
$\P^{a,b}_{\ell,\mathfrak m, f}$ on these hypercubes.

\subsection{Equations for $\mathbf R$}\label{eqR}
In this subsection, we explain how to put coordinates on the toric
tiling  $\mathbf R$ associated to  the tiling of $C^1(G, \R)$ by the
hypercubes $\square_{\alpha}$.

\medskip

We fix an orientation $\mathfrak o: E \rightarrow \E$ of the edges of
$G$. We set $E^{\mathfrak o}:=\mathfrak o(E)$. 

\medskip

For each edge $e\in E^\o$ and each integer $i\in \Z$, let $\P^1_{e,i}$ be a copy of the projective line with projective coordinates $(\x_{e,i}:\x_{\overline e,i})$. Consider the product 
$\mathbf{P} = \prod_{e \in E^\o} \prod_{i\in \Z} \P^1_{e,i}$. It can
be equipped with a natural structure of a projective limit of schemes
of finite type, or can be simply seen as the functor of sets it
represents on the category of rings or schemes. For our purposes in
the present article, we find it more convenient 
to just view $\mathbf{P}$ as a set (by looking at $\k$-points).

\medskip

Define the subset $\mathbf R' \subset \prod_{e \in E^\o} \prod_{i\in \Z} \P^1_{e,i} $ by the equations
 \medskip
 
 \begin{itemize}
 \item[(1)] $\x_{e,i}\x_{\ol e, j}=0$ for all pairs $(e,i), (e,j) \in E^\o \times \Z$ with $ j>i$, 
 \end{itemize}
with the requirement that 
\begin{itemize}
  \item[(2)] For each $e\in E^\o$, there exist indices $i,j  \in \Z$ such that $\x_{e,i} \neq 0$ and  $\x_{\overline e, j} \neq 0.$
\end{itemize}
 Alternatively, for each $e\in E^\o$,  denote by $0_e$
(resp.~$\infty_e$) the point of $\prod_{i\in \Z} \P^1_{e,i}$ with
coordinates $(\x_{e,i}:\x_{\overline e,i}) = (0:1)$
(resp.~$(\x_{e,i}:\x_{\overline e,i})= (1:0)$) for every $i$. Set 
$$
\mathbf P_e := \prod_{i\in \Z} \P^1_{e, i} - \{0_e, \infty_e\}.
$$
Then $\mathbf R'$ is the subset of $\prod_{e\in E^{\mathfrak o}}\mathbf
P_e$ given by Equations (1).

Recall from Subsection~\ref{tortilR} the definition of 
$\mathbf R_e$ as a doubly infinite chain of rational 
smooth curves and the characterization of the toric tiling $\mathbf R$
as the product of the $\mathbf R_e$ for $e\in E^{\mathfrak o}$.

\begin{prop} The locally of finite type scheme $\mathbf
  R_e$ is naturally identified
  with the subset of  $\mathbf P_e$ given by the equations
  $\x_{e,i}\x_{\ol e, j}=0$ for all $i,j\in\Z$ with $j>i$.
\end{prop}

\begin{proof} For each $i\in\Z$, set $0_{e,i}:= (0:1) \in \mathbf
  P^1_{e,i}$ and $\infty_{e,i} := (1:0) \in \mathbf P^1_{e,i}$, and put
$$
\P^1_i : = \prod_{j<i}\{0_{e,j}\} \times \mathbf P^1_{e,i} \times
\prod_{j>i} \{\infty_{e,j}\}.
$$
Note that we have a canonical identification $\mathbf P^1_i \simeq
\mathbf P^1$, and under this identification, the point $\infty$ on
$\mathbf P^1_i$ 
is the point 
$\prod_{j < i}\{0_{e,j}\} \times \prod_{j\geq i}\{\infty_{e,j}\}$, 
whereas the point $0$ on $\mathbf P^1_i$ is 
$\prod_{j\leq i}\{0_{e,j}\} \times \prod_{j > i}\{\infty_{e,j}\}$. In particular, the two subsets
$\mathbf P^1_i$ and $\mathbf P^1_{i+1}$ intersect at the point
$0$ on $\mathbf P^1_i$, which is also the point $\infty$ on $\mathbf
P^1_{i+1}$. It follows that $\mathbf R_e$ is identified with the doubly
infinite chain of $\mathbf P^1$ obtained by taking the union of the
$\mathbf P^1_i$ in $\mathbf P_e$.

We need only show that the stated equations define the union of the
$\mathbf P^1_i$. Clearly, each point on $\P^1_i$ for each $i\in\Z$ satisfies the
equations. Conversely, since $\mathbf P_e$ contains neither
$0_e$ nor $\infty_e$, the coordinates of a point on $\mathbf P_e$
satisfy $\x_{e,i}\neq 0$ for some $i\in\Z$ and $\x_{\ol e,j}\neq 0$ for some
$j\in\Z$. Moreover, if the stated equations are satisfied at the point, then there is
a minimum $i$ with $\x_{e,i}\neq 0$ and a maximum $j$ with
$\x_{\ol e,j}\neq 0$. Thus $\x_{e,m}=0$ for $m<i$ and $\x_{\ol e,m}=0$ for
$m>j$ at the point. From the equations we get $\x_{\ol e,m}=0$, whence
$\x_{e,m}\neq 0$ for $m>i$, and
$\x_{e,m}=0$, whence $\x_{\ol e,m}\neq 0$ for $m<j$. Thus 
$\x_{e,m}\neq 0$ if and only if $m\geq i$ and $\x_{\ol e,m}\neq 0$ if
and only if $m\leq j$. So $j\geq i-1$. But if $j\geq
i+1$, then $\x_{e,i}\x_{\ol e,i+1}\neq 0$, contradicting the
equations. Thus $j=i$ or $j=i-1$. It follows that the point lies on $\P^1_i$.
\end{proof}

\begin{prop} The toric tiling $\mathbf R$ associated to the
  decomposition of $C^1(G, \R)$ by the hypercubes $\square_{\alpha}$
  is naturally identified with $\mathbf R'$. 
\end{prop}

\begin{proof} This follows from the previous proposition. 
\end{proof}

From now one we will write $\mathbf R$ for $\mathbf R'$.

\subsection{Equations for $Y_{\ell, \m}^{a,b}$} We keep the notations
of the previous section. Let $\ell: E\rightarrow \mathbb N$ be a
length function, $\m\in C^1(G,\Z)$ be a twisting, and $a$ and $b$ be two characters of $C^1(G, \mathbb Z)$ and $H^1(G, \mathbb Z)$. 

For each $\alpha \in C^1(G, \Z)$ and $\gamma \in H^1(G, \Z)$, 
consider a copy $\P^1_{\alpha, \gamma}$ of the projective line with
projective coordinates $(\p_{\alpha, \gamma} :  \q_{\alpha,
  \gamma})$. As usual, let $0_{\alpha,\gamma}$ be the point on
$\P^1_{\alpha, \gamma}$ given by $\p_{\alpha, \gamma}=0$ and
$\infty_{\alpha,\gamma}$ 
that given by $\q_{\alpha, \gamma}=0$

Consider the product:
$$
\mathfrak Z:=\prod_{\substack{\alpha \in C^1(G, \Z)\\ \gamma \in
    H^1(G, \Z)}}\P^1_{\alpha, \gamma}.
$$
For each $\mathfrak z \in \mathfrak Z$, let  $Y_{\mathfrak z} \subset \bf R \simeq \bf R \times \{\mathfrak z\} \subset \bf R \times \frak Z$ be the subscheme of $\bf R$ given by the equations 
\begin{align*}
\forall&  \alpha \in C^1(G, \Z), \,\, \gamma \in H^1(G, \Z), \, \\
 \p_{\alpha, \gamma} \prod_{\substack{e\in E^\o\\ \gamma_e>0}} \x_{e,\alpha_e}^{\gamma_e}&
\prod_{\substack{e\in E^\o\\ \gamma_e<0}} \x_{\overline e,\alpha_e}^{-\gamma_e} 
= \q_{\alpha, \gamma} \prod_{\substack{e\in E^\o\\ \gamma_e>0}} 
\x_{\overline e,\alpha_e}^{\gamma_e} \prod_{\substack{e\in E^\o\\\gamma_e<0}} 
\x_{e, \alpha_e}^{-\gamma_e}.
\end{align*}

\begin{thm}\label{thm:main2} Let $\ell\: E\to\N$ be a length function, 
$\m\in C^1(G,\Z)$ and let
$$
a\: C^1(G,\Z)\to \k^*\quad\text{and}\quad b\: H^1(G,\Z)\to \k^*
$$
be characters. Let $\mathfrak z\in\mathfrak Z$ with coordinates for
each $\alpha \in C^1(G, \Z)$ and $\gamma \in H^1(G, \Z)$ satisfying:
\begin{itemize}
\item If $\sum_{e\in E^{\mathfrak o}} (\alpha_e \gamma_e\ell_e - \m_e\gamma_e) <0$, then $(\p_{\alpha, \gamma} :  \q_{\alpha, \gamma}) = (1:0)$ 

\medskip

\item If $\sum_{e\in E^{\mathfrak o}} (\alpha_e \gamma_e\ell_e - \m_e\gamma_e) >0$, then $(\p_{\alpha, \gamma} :  \q_{\alpha, \gamma}) = (0:1)$;

\medskip

\item If $\sum_{e\in E^{\mathfrak o}} (\alpha_e \gamma_e\ell_e -
  \m_e\gamma_e) =0$, then $(\p_{\alpha, \gamma} :  \q_{\alpha, \gamma})
  = (1: b(\gamma) \prod_{e\in E^{\mathfrak o}} a_e^{\alpha_e\gamma_e})$.
\end{itemize} 
Then $Y_{\mathfrak z} = Y_{\ell,\m}^{a,b}$ as subsets of $\mathbf R$.
\end{thm}

The rest of this section is devoted to the proof of this theorem. 

\medskip

First we prove that the scheme $Y^{a,b}_{\ell,\m}$ is included in
$Y_{\mathfrak z}$. Let 
$f\in C^0(G, \mathbb Z)$. Recall the definition of $\dl^\m_f\in
C^1(G,\R)$ from Section~\ref{sec:recap}: for each $e=uv$ in $\mathbb E$ we have
that $\dl^m_f(e)$ is the ratio
$\frac{f(v)-f(u)+\m_e}{\ell_e}$ if integer; otherwise 
$\dl^m_f(e)=\lfloor \frac{f(v)-f(u)+\m_e}{\ell_e}\rfloor+\frac{1}{2}$. Also, $G^\m_f$ is 
the spanning subgraph of $G$ supported on the set of
edges at which $\dl^\m_f$ is an integer. Recall from Subsection~\ref{tortilR}
the definition of
$$
\mathbf P_{\dl^\m_f}:=\prod_{e\in E^{\mathfrak o}}\mathbf
  P^1_{e,\dl^\m_f(e)}\subset\mathbf R.
$$
Under the identification in Section~\ref{eqR}, as a subset of $\prod_j\mathbf P^1_{e,j}$,
$$
\mathbf P^1_{e,i}=\prod_{j<i}\{0_{e,j}\} \times \mathbf P^1_{e,i} \times
\prod_{j>i} \{\infty_{e,j}\}\quad\text{if }i\in\Z
$$
and 
$$
\mathbf P^1_{e,i}=\prod_{j<i}\{0_{e,j}\} \times
\prod_{j>i} \{\infty_{e,j}\}\quad\text{if }i\in\frac{1}{2}\Z-\Z.
$$

The toric arrangement $Y^{a,b}_{\ell,\m}$ is supported in the union of
the $\mathbf P_{\dl^\m_f}$ as $f$ ranges in $C^0(G,\Z)$. For a given
$f$, the variety $Y_{a,b}^{\ell,\m} \cap \mathbf P_{\dl^\m_f}$ is
given by the equations 
\[
\forall \textrm{ oriented cycle } \gamma \textrm { in $G^\m_f$},
\]
\[
\prod_{e\in\bar\gamma\cap E^{\mathfrak o}}b_ea_e^{\dl^\m_f(e)}\prod_{e\in \gamma\cap E^{\mathfrak o}} \x_{e, \dl^\m_f(e)} \prod_{e\in\bar\gamma\cap E^{\mathfrak o}}\x_{\ol e, \dl^\m_f(e)}
=\prod_{e\in \gamma\cap E^{\mathfrak o}} b_ea_e^{\dl^\m_f(e)} \prod_{e\in \bar \gamma\cap E^{\mathfrak o}} \x_{e, \dl^\m_f(e)}\prod_{e\in \gamma\cap E^{\mathfrak o}} \x_{\ol e, \dl^\m_f(e)}.
\]
(Here, abusing the notation, $b$ denotes an extension of $b$ to a
character of $C^1(G,\Z)$.)
 
We show first the following claim:

\begin{claim} \label{claim:th1} We have $Y^{a,b}_{\ell,\m}\cap \mathbf P_{\dl^\m_f}=
Y_{\mathfrak z}\cap \mathbf P_{\dl^\m_f}$ for each $f \in C^0(G,
\mathbb Z)$. In particular, $Y^{a,b}_{\ell,\m}\subseteq Y_{\mathfrak z}$.
\end{claim}

\begin{proof} The second statement is a consequence of the first, as
$Y^{a,b}_{\ell,\m}$ is supported in the union of the $\mathbf P_{\dl^\m_f}$.

Fix $f \in C^0(G,\mathbb Z)$. Let $\beta$ be the element of $C^1(G, \mathbb Q)$ given by
 $\beta_e = \frac{f(v) -f(u)+\m_e}{\ell_e}$ for each oriented edge
 $e=uv\in\mathbb E$. Note that we have for each cycle $\gamma \in H^1(G, \mathbb Z)$, 
  \begin{equation}\label{eq:20}
  \sum_{e \in E^{\mathfrak o}} (\beta_e \gamma_e \ell_e - \m_e\gamma_e) = 0.
  \end{equation}  
 
The defining equations of $Y_{\mathfrak z}$ are obtained by taking
$\alpha \in C^1(G,\mathbb Z)$ 
and $\gamma \in H^1(G, \mathbb Z)$, and looking at the sum 
$\sum_{e\in E^{\mathfrak o}} (\alpha_e \gamma_e\ell_e -
\m_e\gamma_e)$: 
the associated equation is 
\begin{equation}\label{PQalpha}
\p_{\alpha, \gamma} \prod_{\substack{e\in E^\o\\ \gamma_e>0}} \x_{e,\alpha_e}^{\gamma_e}
\prod_{\substack{e\in E^\o\\ \gamma_e<0}} \x_{\overline e,\alpha_e}^{-\gamma_e} 
= \q_{\alpha, \gamma} \prod_{\substack{e\in E^\o\\ \gamma_e>0}} 
\x_{\overline e,\alpha_e}^{\gamma_e} \prod_{\substack{e\in E^\o\\\gamma_e<0}} 
\x_{e, \alpha_e}^{-\gamma_e}.
\end{equation}
for the point $(\p_{\alpha,\gamma}:\q_{\alpha, \gamma})$ given in the
statement of the claim, according to the sign of the sum.
 
Suppose first that  $\sum_{e\in E^{\mathfrak o}} (\alpha_e \gamma_e
\ell_e - \m_e\gamma_e) < 0$. 
In this case, the equation corresponding to 
$(\alpha, \gamma)$ is 
$$
\prod_{\substack{e\in E^\o\\ \gamma_e>0}} \x_{e,\alpha_e}^{\gamma_e}
\prod_{\substack{e\in E^\o\\ \gamma_e<0}} \x_{\overline e,\alpha_e}^{-\gamma_e} 
= 0.
$$
From Equation~\eqref{eq:20}, we infer the existence of $e\in E^{\mathfrak o} $ for which 
$\alpha_e\gamma_e < \beta_e \gamma_e$. Then $\gamma_e\neq 0$. If $\gamma_e
>0$ then $\alpha_e < \dl^\m_f(e)$, and if $\gamma_e <0$ then $\alpha_e>
\dl^\m_f(e)$. By definition, each point on 
$\P_{\dl^\m_f}$ has coordinates satisfying $\x_{\ol e,i}=0$ for
$i>\dl^\m_f(e)$ and $\x_{e,i} =0$ for $i<\dl^\m_f(e)$.  So,
in any case, the equation corresponding to 
$(\alpha, \gamma)$ is automatically verified on $\P_{\dl^\m_f}$. The case where 
$\sum_{e\in E^{\mathfrak o}} (\alpha_e \gamma_e \ell_e - \m_e\gamma_e) > 0$ follows similarly.

It remains thus to treat the last case, where  the sum 
$\sum_{e\in E^{\mathfrak o}}(\alpha_e \gamma_e \ell_e - \m_e\gamma_e)$
vanishes. 
It follows from Equation~\eqref{eq:20} that two cases can happen:
 \begin{itemize}
 \item[(i)] For every $e$ in the support of $\gamma$, we have $\alpha_e = \beta_e$.
 \item[(ii)] There exist two edges $e_1$ and $e_2$ in $E^{\mathfrak
     o}$ such that  $\alpha_{e_1}\gamma_{e_1} > \beta_{e_1} \gamma_{e_1}$ and $\alpha_{e_2}\gamma_{e_2} < \beta_{e_2} \gamma_{e_2}$. 
 \end{itemize}
 
In Case (i), all the edges $e$ in the support of $\gamma$ belong to
the subgraph $G_f^\m$, and we have $\alpha_e = \dl^\m_f(e)$. If
$\gamma$ is an oriented cycle, by the definition of the point
$(\p_{\alpha, \gamma}:\q_{\alpha, \gamma})$, the equation corresponding
to $(\alpha, \gamma)$ is precisely the
equation in $\mathbf P_{\dl^\m_f}$ associated to $\gamma$. At any
rate, $\gamma$ is a sum of oriented cycles
supported in $G_f^\m$, whence the equation corresponding
to $(\alpha, \gamma)$ is a product of equations in $\mathbf
P_{\dl^\m_f}$ associated to oriented cycles.

In Case (ii), from $\alpha_{e_1}\gamma_{e_1} > \beta_{e_1} \gamma_{e_1}$ we get:
\begin{itemize}
\item Either $\gamma_{e_1}>0$, in which case $\alpha_{e_1} >
\dl^\m_f(e_1)$ and so $\x_{\ol e_1, \alpha_{e_1}} =0$ on $\mathbf P_{\dl^\m_f}$;
\item Or $\gamma_{e_1}<0$, in which case $\alpha_{e_1} <
\dl^\m_f(e_1)$ and so $\x_{e_1, \alpha_{e_1}} =0$ on $\mathbf P_{\dl^\m_f}$.
\end{itemize}
In either case,
$$
\prod_{\substack{e\in E^\o\\ \gamma_e>0}} \x_{\overline e,\alpha_e}^{\gamma_e}
\prod_{\substack{e\in E^\o\\ \gamma_e<0}} \x_{e,\alpha_e}^{-\gamma_e} 
= 0\quad\text{on $\mathbf P_{\dl^\m_f}$.}
$$
Analogously, from $\alpha_{e_2}\gamma_{e_2} < \beta_{e_2} \gamma_{e_2}$ we get 
$$
\prod_{\substack{e\in E^\o\\ \gamma_e>0}} \x_{e,\alpha_e}^{\gamma_e}
\prod_{\substack{e\in E^\o\\ \gamma_e<0}} \x_{\overline e,\alpha_e}^{-\gamma_e} 
= 0\quad\text{on $\mathbf P_{\dl^\m_f}$.}
$$
In particular, once again Equation~\eqref{PQalpha} 
corresponding to $(\alpha, \gamma)$ 
is automatically verified on $\mathbf P_{\dl^\m_f}$.
\end{proof}

It follows from the above claim that to prove Theorem~\ref{thm:main2} we need only
show that $Y_{\mathfrak z}$ is contained in the union of the $\mathbf
P_{\dl^\m_f}$. This will follow from a series of claims, culminating
with Claim~\ref{finalclaim:main2}.

\begin{claim}\label{claim2}
Let $\alpha \in C^1(G, \mathbb R)$ and $\gamma \in H^1(G, \mathbb
Z)$. Then:
\begin{itemize}
\item If $\sum_{e\in E^{\mathfrak o}} (\alpha_e\gamma_e\ell_e - \m_e\gamma_e) <0$ then 
$$
Y_{\mathfrak z} \subset \bigcup\mathbf P_{\beta},
$$
with the union over the $\beta\in C^1(G,\Z)$ such that
$(\alpha_e - \beta_e) \gamma_e <0$ for some $e\in E^\o$.
\item If $\sum_{e\in E^{\mathfrak o}} (\alpha_e\gamma_e\ell_e - \m_e\gamma_e) >0$ then 
$$
Y_{\mathfrak z} \subset \bigcup\mathbf P_{\beta},
$$
with the union over the $\beta\in C^1(G,\Z)$ such that 
$(\alpha_e - \beta_e) \gamma_e >0$ for some $e\in E^\o$.
\end{itemize}
\end{claim}

\begin{proof} We treat only the first case; the second follows with a similar argument.

Suppose first $\alpha \in C^1(G, \mathbb Z)$. In this case,
$$
\prod_{\substack{e\in E^\o\\ \gamma_e>0}} \x_{e,\alpha_e}^{\gamma_e}
\prod_{\substack{e\in E^\o\\ \gamma_e<0}} \x_{\overline e,\alpha_e}^{-\gamma_e} 
= 0
$$
on $Y_{\mathfrak z}$. Then the claim follows, since for each $i\in\Z$,
$$
\x_{\ol e,i}=0\text{ defines } \bigcup_{\substack{\beta\in 
C^{1}(G,\Z)\\ \beta_e<i}} \mathbf P_{\beta}
\quad\text{and}\quad
\x_{e,i}=0\text{ defines } \bigcup_{\substack{\beta\in
C^{1}(G,\Z)\\ \beta_e> i}} \mathbf P_{\beta}.
$$

In the general case,  let $\tilde \alpha \in C^1(G, \mathbb Z)$
satisfying, for each $e\in E^{\mathfrak o}$,
\[\tilde\alpha_e = 
\begin{cases} \lceil \alpha_e \rceil & \textrm{if } \gamma_e\leq 0\\
\lfloor \alpha_e \rfloor & \textrm{if } \gamma_e > 0.
\end{cases}\]
Clearly,
\[
\sum_{e\in E^{\mathfrak o}} \tilde \alpha_e \gamma_e\ell_e \leq
\sum_{e\in E^{\mathfrak o}} \alpha_e \gamma_e\ell_e <\sum_{e\in
  E^{\mathfrak o}}\mathfrak m_e\gamma_e.
\]
We may thus apply what we have just proved to conclude that 
$$
Y_{\mathfrak z} \subset \bigcup \mathbf P_{\beta},
$$
where the union is over the $\beta\in C^1(G,\Z)$ such that
$(\tilde\alpha_e - \beta_e) \gamma_e <0$ for some $e\in E^\o$. But 
$(\tilde \alpha_e - \beta_e) \gamma_e < 0$ if and only if 
$(\alpha_e - \beta_e) \gamma_e <0$.
\end{proof}

\begin{claim}\label{claim3} Let $\alpha\in C^1(G,\frac 12 \mathbb Z)$
  such that $Y_{\mathfrak z}$ has a point on the interior of $\mathbf
  P_{\alpha}$. Then, for each $\gamma \in H^1(G, \mathbb Z)$, we have 
\[
\Bigl| \sum_{e\in E^{\mathfrak o}} (\alpha_e \gamma_e \ell_e -
\m_e\gamma_e) \Bigr| \leq \frac 12 \bigl(\sum_{\substack{e \in
    E^{\mathfrak o}\\ \alpha_e \not\in \mathbb Z }} \ell_e
|\gamma_e|\bigr).
\]
\end{claim}

\begin{proof}
Suppose by contradiction that
\[
\Bigl| \sum_{e\in E^{\mathfrak o}} (\alpha_e \gamma_e \ell_e -
\m_e\gamma_e) \Bigr| > \frac 12 \bigl(\sum_{\substack{e \in
    E^{\mathfrak o}\\ \alpha_e \not\in \mathbb Z}} 
\ell_e |\gamma_e|\bigr)
\]
for a certain $\gamma \in H^1(G, \mathbb Z)$. 
Suppose first $\sum_{e\in E^{\mathfrak o}}(\alpha_e \gamma_e \ell_e - \m_e\gamma_e) \geq 0$. Then
$$
\sum_{\substack{e\in E^{\mathfrak o}\\ \alpha_e\in \mathbb Z}}
\alpha_e \gamma_e \ell_e
+\sum_{\substack{e\in E^{\mathfrak o}\\ \alpha_e \not\in\mathbb Z\\
\gamma_e \geq 0}} (\alpha_e-1/2) \gamma_e \ell_e 
+ \sum_{\substack{e\in E^{\mathfrak o}\\ \alpha_e \not\in \mathbb Z\\
\gamma_e < 0}} (\alpha_e+1/2) \gamma_e \ell_e  
- \sum_{e\in E^{\mathfrak o}}\m_e\gamma_e >0.
$$
Let $\tilde \alpha \in C^1(G, \mathbb Z)$
satisfying, for each $e\in E^{\mathfrak o}$,
$$
\tilde \alpha_e = \begin{cases}\alpha_e & \textrm{ if } \alpha_e \in \mathbb Z\\
\alpha_e- \frac 12 & \textrm{ if } \alpha_e \not \in \mathbb Z \textrm{ and } \gamma_e\geq0 \\
\alpha_e+\frac 12 & \textrm{ if } \alpha_e \not \in \mathbb Z \textrm{
  and } \gamma_e< 0.
\end{cases}
$$ 
Applying Claim~\ref{claim2}, we get that 
\[ 
Y_{\mathfrak z} \subset \bigcup \mathbf P_{\beta}
\]
where the union is over the $\beta\in C^1(G,\Z)$ such that
$(\tilde\alpha_e - \beta_e) \gamma_e >0$ for some $e\in E^\o$. For
any such $\beta$ we have that either $\beta_e<\alpha_e-1/2$ or
$\beta_e>\alpha_e+1/2$ for some $e\in E^\o$. But in neither case there
is a point on $\mathbf
P_{\beta}$ lying on the interior of $\mathbf P_{\alpha}$,
contradicting the hypothesis on $Y_{\mathfrak z}$.

The case where 
$\sum_{e\in E^{\mathfrak o}}(\alpha_e \gamma_e \ell_e - \m_e\gamma_e) \leq 0$ is treated similarly.
\end{proof}

\begin{claim} Let $\alpha \in C^1(G, \mathbb Z)$ be such that
  $Y_{\mathfrak z}$ has a point on the interior of $\mathbf P_{\alpha}$. 
Then there exists $f\in C^0(G, \mathbb Z)$ such that $\alpha = \dl^\m_f$.
\end{claim}

\begin{proof} By Claim~\ref{claim3}, since $\alpha$ is integer valued, we get 
\[
\forall\, \gamma \in H^1(G, \mathbb Z), \qquad 
\sum_{e\in E^{\mathfrak o}} (\alpha_e \gamma_e \ell_e - \m_e\gamma_e)
=0.
\]
Define $\beta \in C^1(G, \mathbb Z)$ by setting $\beta_e
:=\alpha_e\ell_e - \m_e$ for each 
$e\in \E$. Then 
\[
\sum_{e\in E^{\mathfrak o}} \beta_e \gamma_e =0
\]
for every $\gamma \in H^1(G, \mathbb Z)$. This guarantees the
existence of a function $f\in C^0(G, \mathbb Z)$ such that $\beta_e =
f(v)-f(u)$ for each oriented edge $e= uv$ in $\E$. Equivalently, 
\[
\alpha_e = \frac{f(v)-f(u)+\m_e}{\ell_e}
\]
for each $e= uv$ in $\E$, or $\alpha = \dl^{\m}_f$. 
\end{proof}

A weaker statement holds for the general case of an 
$\alpha \in C^1(G, \frac 12\mathbb Z)$ which is not necessarily
integral. For the claim we use the following lemma.

\begin{lemma}\label{lem:th2} Let $G=(V,E)$ be a connected loopless graph and
  $\mathfrak o$ an orientation. Let $\beta \in C^1(G, \mathbb Z)$ and 
$h\colon E \to \mathbb Z_{\geq 0}$. Suppose that for each $\gamma \in H^1(G, \mathbb Z)$, we have 
\[
\Bigl|\sum_{e\in E^{\mathfrak o}} \beta_e\gamma_e\Bigr| \leq \sum_{e
  \in E^{\mathfrak o}} h(e)|\gamma_e|.
\]
Then there exists an element $\eta \in C^1(G, \mathbb Z)$ that verifies the following two properties:
\begin{itemize}
\item for each $e\in \E$, we have $|\eta_e| \leq h(e)$;
\item for each $\gamma\in H^1(G, \mathbb Z)$, we have 
\[\sum_{e\in E^{\mathfrak o}} \beta_e \gamma_e = \sum_{e\in E^{\mathfrak o}} \eta_e \gamma_e.\] 
\end{itemize}
\end{lemma}

The proof of this lemma is given in Section~\ref{sec:lem:th2}. It will
be used in the proof of Claim~\ref{finalclaim:main2} below.

\begin{claim}\label{finalclaim:main2} Let $\alpha \in C^1(G,
  \frac{1}{2}\mathbb Z)$ 
such that
  $Y_{\mathfrak z}$ has a point on the interior of $\mathbf
  P_{\alpha}$. Then there exists $f\in C^0(G, \mathbb Z)$ such that 
$\mathbf P_\alpha \subseteq \mathbf P_{\dl^\m_f}$.
\end{claim}

\begin{proof} By Claim~\ref{claim3}, we have the inequality for each 
  $\gamma\in H^1(G,\Z)$:
 \begin{equation}\label{eq:th2}
 \Bigl| \sum_{e\in E^{\mathfrak o}} (\alpha_e \gamma_e \ell_e -
 \m_e\gamma_e) \Bigr| \leq \frac 12 \bigl(\sum_{\substack{e \in
     E^{\mathfrak o}\\ \alpha_e \not\in \mathbb Z }} \ell_e |\gamma_e|\bigr).
 \end{equation} 
 Define $\beta\in C^1(G, \mathbb Z)$ by $\beta_e:= 2\alpha_e\ell_e - 2\m_e$ for each oriented edge $e \in \E$. Define the function $h: E \rightarrow \mathbb Z_{\geq 0}$ as follows:
 \[\forall\, e\in E, \qquad h(e) =\begin{cases} \ell_e & \textrm{if } \alpha_e \notin \mathbb Z\\
 0 & \textrm{otherwise}. 
 \end{cases}\]
By these definitions, Inequalities~\eqref{eq:th2} can be rewritten in the form
\[\Bigl|\sum_{e\in E^{\mathfrak o}} \beta_e\gamma_e\Bigr| \leq \sum_{e \in E^{\mathfrak o}} h(e)|\gamma_e|.\]
Applying Lemma~\ref{lem:th2}, we infer the existence of an element
$\eta\in C^1(G, \mathbb Z)$ verifying the following two properties:
\begin{itemize}
\item for each $e\in \E$, we have $|\eta_e| \leq h(e)$;
\item for each $\gamma\in H^1(G, \mathbb Z)$, we have 
\[\sum_{e\in E^{\mathfrak o}} \beta_e \gamma_e = \sum_{e\in E^{\mathfrak o}} \eta_e \gamma_e.\] 
\end{itemize}
This implies the existence of a function $\tilde f : V \rightarrow \mathbb Z$ such that for each oriented edge $e = uv$ in $\E$, we have 
\[\beta_e - \eta_e = \tilde f(v) - \tilde f(u).\]  
Write $\tilde f = 2 f + \epsilon$ for $\epsilon : V \rightarrow
\{0,1\}$ and $f\in C^0(G, \mathbb Z)$. We claim that $\mathbf P_\alpha \subseteq \mathbf P_{\dl^\m_f}$.

Indeed, for each oriented edge $e=uv$ in $\E$, we have 
\begin{equation}\label{eq3}
2\alpha_e\ell_e -\eta_e -2\m_e = 2(f(v)-f(u)) + \epsilon(v) -\epsilon(u),
\end{equation}
which implies that 
\begin{equation}\label{eq4}
\frac{f(v)-f(u)+ \m_e}{\ell_e} = \alpha_e - \frac{\eta_e -\epsilon(u)+\epsilon(v)}{2\ell_e}.
\end{equation}
Then two cases can happen, depending on whether $\alpha_e \in \mathbb Z$ or not:

\medskip

$\bullet$ If $\alpha_e \in\mathbb Z$, then $|\eta_e| \leq h(e)=0$,
and thus $\eta_e =0$. In addition, from Equation~\eqref{eq3} we get 
$\epsilon(v)-\epsilon(u) \equiv 0$ (mod 2), which implies that
$\epsilon(v)-\epsilon(u)=0$, and hence
\[
\alpha_e = \frac{f(v)-f(u)+ \m_e}{\ell_e} = \dl^\m_f(e).
\]

\medskip

$\bullet$ If $\alpha_e \notin \mathbb Z$ then $\alpha_e = a +
\frac 12$ for  an integer $a$, and 
\[
\frac{f(v)-f(u)+ \m_e}{\ell_e} = a + \frac{\ell_e- \eta_e
  +\epsilon(u)-\epsilon(v)}{2\ell_e}.
\]
Also, $|\eta_e|\leq h(e)=\ell_e$. We claim that 
\begin{equation}\label{aalpha}
0 \leq \frac{\ell_e- \eta_e +\epsilon(u)-\epsilon(v)}{2\ell_e} \leq 1.
\end{equation}
Indeed, since $|\epsilon(u) -\epsilon(v)| \leq 1$ this is obviously
true if $|\eta_e|<\ell_e$. Now, if $\eta_e =\pm\ell_e$, from
Equation~\eqref{eq3} we get, since $\alpha_e = a +\frac 12$, 
\[
\epsilon(v) -\epsilon(u) \equiv \ell_e - \eta_e \equiv 0\textrm{ (mod
  2)}.
\]
This implies that $\epsilon(u) =\epsilon(v)$, and then
\[
\frac{\ell_e- \eta_e +\epsilon(u)-\epsilon(v)}{2\ell_e} = 0 \textrm{
  or }1.
\]
It follows that $\dl^\m_f(e)=a+\frac
12$ if both inequalities in \eqref{aalpha} are strict, whereas 
$\dl^\m_f(e)=a$ or $\dl^\m_f(e)=a+1$ otherwise. 
Since $\alpha_e = a +\frac 12$, we get that either $\alpha_e =
\dl^\m_f(e)$ or we have 
$\dl^\m_f(e) \in \mathbb Z$ and $\alpha_e = \dl^\m_f(e) \pm \frac 12$.

\medskip

In any case $\mathbf P_{e,\alpha_e}\subseteq\mathbf P_{e,\dl^\m_f(e)}$ 
for each $e\in E^\o$, and thus $\mathbf P_{\alpha} \subseteq \mathbf P_{\dl^\m_f}$.
\end{proof}

The interiors of the $\mathbf P_\alpha$ for $\alpha\in
C^0(G,\frac{1}{2}\mathbb Z)$ form a stratification of $\mathbf
R$. Thus the above claim finishes the proof of Theorem~\ref{thm:main2}.

\subsection{Proof of Lemma~\ref{lem:th2}}\label{sec:lem:th2}
In this section we prove Lemma~\ref{lem:th2}. We proceed by induction on the quantity $\Sigma(h):=\sum_{e\in E} h(e)$.  If $h = 0$, then the claim holds trivially for $\eta=0$. Suppose now $h$ is a nonzero function and assume the claim holds for all functions $h'$ with $\Sigma(h') < \Sigma(h)$. 

Note that for an oriented cycle $\gamma$ we have
$$
\sum_{e\in E^{\mathfrak o}}\beta_e\gamma_e=\sum_{e\in\gamma}\beta_e\gamma_e=\sum_{e\in\gamma}\beta_e
\quad\text{and}\quad 
\sum_{e\in E^{\mathfrak o}}h(e)|\gamma_e|=\sum_{e\in\gamma} h(e)\gamma_e=\sum_{e\in\gamma} h(e).
$$
Let $A$ be the set of oriented edges $e\in \E$ for which there exists an oriented cycle $\gamma$ in the graph which contains $e$ and such that the equality holds: 
\[
\sum_{f \in\gamma} \beta_f = \sum_{f\in
  \gamma} h(f).
\]
Note that  it might happen for an oriented edge $e\in\mathbb E$ that both  $e$ and $\overline e$ belong to $A$.  Let $A^+$ be the set of all oriented edges in $A$ with $h(e)>0$.

\begin{claim}\label{eitherorA} For each oriented edge $e\in\mathbb E$ with $h(e)>0$,
  the set $A^+$ contains at most one of the two oriented edges $e$ and $\overline e$. 
\end{claim}

\begin{proof} For the sake of a contradiction, let $e$ be an oriented edge with $h(e)>0$, and let $\gamma_1$ and $\gamma_2$ be two oriented cycles such that  
\[
\sum_{f \in \gamma_i} \beta_f = \sum_{f\in\gamma_i} h(f),
\]
and $e\in \gamma_1$ and $\overline e \in \gamma_2$.  Let $\gamma = \gamma_1 + \gamma_2$.  
Adding the two equations above for $\gamma_i$, we get 
\[\sum_{f \in E^{\mathfrak o}} \beta_f\gamma_f = \sum_{f\in\gamma_1}
  h(f) +\sum_{f\in \gamma_2} h(f).\]
Note however that 
$e \not\in \mathrm{supp}(\gamma)$, whence
\[
\sum_{f\in E^{\mathfrak o}} h(f)|\gamma_f|<
\sum_{f\in\gamma_1} h(f) +\sum_{f\in\gamma_2} h(f).
\]
Since by hypothesis
\[
\sum_{f\in E^{\mathfrak o}} \beta_f\gamma_f \leq  \sum_{f\in E^{\mathfrak o}}
h(f)|\gamma_f|,
\]
we get
\[
\sum_{f\in E^{\mathfrak o}} h(f)|\gamma_f|<\sum_{f\in\gamma_1} h(f)
+\sum_{f\in\gamma_2} h(f)=\sum_{f \in E^{\mathfrak o}} \beta_f\gamma_f
\leq \sum_{f\in E^{\mathfrak o}} h(f)|\gamma_f|,\]
which is a contradiction.
\end{proof}

Suppose now there is an oriented edge $e_0$ with $h(e_0)>0$ such that
neither $e_0$ nor $\overline e_0$ belongs to $A$. 
Then for the function $\tilde h: E \rightarrow \mathbb Z_{\geq 0}$ defined by $\tilde h(e_0) := h(e_0)-1$ and $\tilde h(e) := h(e)$ on edges $e\neq e_0$ of $G$, we still have the inequalities: 
 \[
\forall \gamma \in H^1(G, \mathbb Z),\,\,\,\Bigl| \sum_{e \in E^{\mathfrak o}} \beta_e\gamma_e \Bigr|\leq \sum_{e\in E^{\mathfrak o}}\tilde h(e)|\gamma_e|.
\]
Indeed, for each $\gamma\in H^1(G, \mathbb Z)$ write $\gamma=\sum\gamma_i$, where the $\gamma_i$ are oriented cycles such that for each $i$ and $f\in\E$, we have that $\gamma_{i,f}>0$ only if $\gamma_{f}>0$. Then 
$$
\Bigl| \sum_{e \in E^{\mathfrak o}} \beta_e\gamma_e \Bigr|\leq 
\sum_i \Bigl| \sum_{e \in\gamma_{i}} \beta_e\Bigr|\leq
\sum_i\sum_{e\in\gamma_{i}}\tilde h(e)=
\sum_i\sum_{e\in E^{\mathfrak o}}\tilde h(e) |\gamma_{i,e}|=
\sum_{e\in E^{\mathfrak o}}\tilde h(e)|\gamma_e|,
$$
where the existence of $e_0$ is used in the second inequality. Applying then the induction hypothesis to $\tilde h$, we conclude the proof of the theorem.

We may thus assume that $h(e)>0$ for an oriented edge $e$ if and only if $e\in A^+$ or 
$\ol e\in A^+$. Furthermore, by Claim~\ref{eitherorA}, if $e\in A^+$ then $\ol e\not\in A^+$.

\medskip

The following claim completes the proof of our lemma.

\begin{claim} \label{claim4} Let $\eta\in C^1(G,\mathbb Z)$ defined by setting 
$$
\eta_e:=\begin{cases}
h(e)&\text{if }e\in A^+,\\
-h(e)&\text{if }\overline e\in A^+,\\
0&\text{otherwise.}
\end{cases}
$$
Then 
$$
\sum_{e \in E^{\mathfrak o}} \beta_e \gamma_e = 
\sum_{e\in E^{\mathfrak o}} \eta_e\gamma_e
$$
holds for all $\gamma \in H^1(G, \mathbb Z)$.
\end{claim}

\begin{proof} Notice first that 
$$
\sum_{e\in E^{\mathfrak o}} \eta_e\gamma_e=\sum_{e\in A^+}h(e)\gamma_e
$$
for each $\gamma\in H^1(G, \mathbb Z)$. 
Since every $\gamma\in H^1(G, \mathbb Z)$ is a sum of oriented cycles,
we need only prove the stated equality for an oriented cycle $\gamma$.
Then we need to prove that 
$$
\sum_{e \in\gamma} \beta_e  =\sum_{e\in A^+}h(e)\gamma_e.
$$
Suppose by contradiction this is not the case. Up to replacing $\gamma$ by
$\overline\gamma$, if necessary, we may assume that 
$$
\sum_{e \in\gamma} \beta_e  < \sum_{e\in A^+}h(e)\gamma_e.
$$

Let $e_1, \dots, e_k$ be the set 
of all the oriented edges in $\gamma\cap A^+$. 
For each $i=1,\dots, k$,  since $e_i\in A$, there is an oriented cycle $\gamma_i$ with $e_i\in\gamma_i$ for which the equality 
$$
\sum_{e \in\gamma_i} \beta_e = \sum_{e\in\gamma_i} h(e)
$$
holds. Let $\tilde \gamma = \gamma_1+ \dots+\gamma_k - \gamma$. We have

\begin{align*}
\sum_{e\in E^\mathfrak o} \beta_e \tilde\gamma_e &= \Bigl(\sum_{i=1}^k \sum_{e\in \gamma_i} \beta_e\Bigr) - \sum_{e\in\gamma} \beta_e  \\
&= \Bigl(\sum_{i=1}^k \sum_{\substack{e \in\gamma_i}}h(e) \Bigr) -
  \sum_{\substack{e\in \gamma}} \beta_e.
\end{align*}

We may sum now over the $e\in\E$ such that $h(e)>0$, that is, those
$e$ such that $e\in A^+$ or $\overline e\in A^+$. By
Claim~\ref{eitherorA}, if $e\in A^+$ 
then the $\gamma_{i,e}$ are non-negative, whereas if $\overline e\in
A^+$ then the $\gamma_{i,e}$ are nonpositive. It follows that 
$$
\sum_{i=1}^k \sum_{\substack{e \in \gamma_i}}h(e)=\sum_{e\in A^+}h(e)\sum_{i=1}^k \gamma_{i,e}.
$$
Thus
$$
\sum_{e\in E^\mathfrak o} \beta_e \tilde\gamma_e=\sum_{e\in
  A^+}h(e)\sum_{i=1}^k \gamma_{i,e}- \sum_{e\in \gamma} \beta_e
=\sum_{e\in
  A^+}h(e)\tilde\gamma_{e}+\sum_{e\in
  A^+}h(e)\gamma_e - \sum_{e\in \gamma} \beta_e>\sum_{e\in
  A^+}h(e)\tilde\gamma_{e}.
$$

In addition, if $e\in A^+$, since the $\gamma_{i,e}$ are non-negative,
it follows that  $\tilde\gamma_e\geq 0$ if $\gamma_e\leq 0$; but if 
$\gamma_e>0$ then $e=e_i$ for some $i$, whence $\tilde\gamma_e\geq 0$ 
as well. Thus
$$
\sum_{e\in E^{\mathfrak o}}h(e)|\tilde\gamma_{e}|=
\sum_{e\in A^+}h(e)|\tilde\gamma_{e}|=\sum_{e\in A^+}h(e)\tilde\gamma_{e}.
$$
It follows that 
$$
\sum_{e\in E^\mathfrak o} \beta_e \tilde\gamma_e>\sum_{e\in E^{\mathfrak o}}h(e)|\tilde\gamma_{e}|,
$$
a contradiction.
\end{proof}

\section{Mixed toric tilings IV: Degenerations of tori}\label{sec:degtori}

 The aim  of this section is to prove that each arrangement of
 toric varieties $Y_{\ell,\m}^{a,b}$ is an equivariant degeneration of the torus
 ${\mathbf G}^{|V|-1}_{\mathbf m}$. 

Let $G=(V,E)$ be a connected graph without self loops
 
 \subsection{Degenerations of $\Gm$}\label{degG}
 Notations as in Sections \ref{mixtortil} and \ref{sec:eqtiling}. Fix
 $e\in\E$. In this section we construct a degeneration
 of $\Gm$ to the doubly infinite chain $\mathbf R_e$ of smooth
 rational curves. Furthermore, we construct ``finite approximations'' of this
degeneration, which consist of degenerations of $\Gm$ to open
subschemes $\mathbf R^{(n)}_e\subset \mathbf R_e$ for $n\in \mathbb
N$, as defined below. 
 
 Let $n$ be a nonnegative integer number.  
The scheme $\mathbf R^{(n)}_e$ is  the open subscheme of $\mathbf R_e$ defined by 
 \[
\mathbf R^{(n)}_e := \mathbf R_e - \bigcup_{\substack{i\in \mathbb Z\\
    |i|\geq n+1}} \P^1_{e,i}.\] 
 Note in particular that for $n=0$ we have $\mathbf R_e^{(0)} \simeq \Gm$. 
 
 \medskip
 
Recall that we view $\mathbf R_e$ in
 \[
\P_e= \prod_{i\in \Z} \P^1_{e,i} - \{0_e, \infty_e\}
\]
given by the equations $\x_{e,i}\x_{\ol e, j}=0$ for all $i<j$. 
We denoted by $(\x_{e,i}: \x_{\ol e, i})$ the coordinates of $\P^1_{e,i}$
for each $i$, and by $0_{e,i}$ (resp.~$\infty_{e,i}$) the point
on $\P^1_{e,i}$ given by $\x_{e,i}=0$ (resp.~$\x_{\ol e,i}=0$). Also, 
$0_e$ (resp.~$\infty_e$) is the point projecting to $0_{e,i}$ (resp.~$\infty_{e,i}$) for
each $i$.

Now, for each nonnegative integer $n$, let 
\[ 
\P^{(n)}_e := \Big(\P^1_{e,-n}-\{\infty_{e,-n}\}\Big)\times
\prod_{|i|<n} \P^1_{e,i}\times\Big(\P^1_{e,n}-\{0_{e,n}\}\Big).
\]
Then $\mathbf R^{(n)}_e$ is the closed subscheme
of $\P_e^{(n)}$ given by the equations
 \[
\x_{e,i} \, \x_{\ol e, j}
 =0\qquad \forall -n \leq i<j \leq n
\]

We may view $\mathbf R^{(n)}_e$ as the open subscheme of $\mathbf R_e$
given by $\x_{\ol e,-n}\neq 0$ and $\x_{e,n}\neq 0$, as a point on
$\mathbf R_e$ satisfying these inequalities is a point on 
$$
\prod_{i<-n}\{0_{e,i}\}\times \Big(\P^1_{e,-n}-\{\infty_{e,-n}\}\Big)\times
\prod_{|i|<n}
\P^1_{e,i}\times\Big(\P^1_{e,n}-\{0_{e,n}\}\Big)\times
\prod_{i>n}\{\infty_{e,i}\}
$$
satisfying $\x_{e,i} \, \x_{\ol e, j}=0$ for all $-n \leq i<j \leq n$,
and the above product can be identified with $\P^{(n)}_e$.

\begin{prop}\label{prop:gm} Let $B:= \mathrm{Spec}(\k[[t]])$. For each
  $a\in\Gm(k)$, 
define $\mathfrak \x_e(a)\subseteq\P_e\times B$ by
\begin{itemize}
\item[($\star$)] $\x_{e,i} \x_{\ol e, j} = (at)^{j-i} \x_{\ol e,i}
  \x_{e,j}$ 
for all pairs of integers $i,j$ with $i<j$,
\end{itemize} 
and $\mathfrak X^{(n)}_e(a)\subseteq\P^{(n)}_e\times B$ for
each $n\in\Z_{\geq 0}$ by
\begin{itemize}
\item[($\star_n$)] $\x_{e,i} \x_{\ol e, j} = (at)^{j-i} \x_{\ol e,i}
  \x_{e,j}$ 
for all pairs of integers $i,j$ with $-n \leq  i<j\leq n$.
\end{itemize}
Then $\mathfrak X_e(a)$ and $\mathfrak X_e^{(n)}(a)$ for each
$n\in\Z_{\geq 0}$ are flat over $B$
and satisfy:
\begin{itemize}
\item[(1)] Their total spaces are regular.
\item[(2)] Their generic fibers over $B$ are isomorphic to $\Gm$.
\item[(3)] Their special fibers over $B$ are $\mathbf R_e$ and 
$\mathbf R_e^{(n)}$ for each $n\in \mathbb Z_{\geq 0}$, respectively.
\end{itemize}  
\end{prop}

\begin{proof} Property (2) follows from observing that for a point on
  the generic fiber, all the coordinates $\x_{e,i}$ and $\x_{\ol e,
    i}$ are nonzero; once this has been proved, it follows from 
Equations~($\star$) and ($\star_n$) that the ratio
  $\x_{e,0}/\x_{\ol e, 0}\in\Gm$ determines all the other ratios
  $\x_{e,i}/\x_{\ol e, i}$. But suppose by contradiction that for a
  point on the generic fiber of $\mathfrak X_e(a)$ 
we have $\x_{e,i} =0$ for some $i$. Then it follows from Equations
($\star$) that also $\x_{e, j}=0$ for every $j>i$. And, using the
same equations with $i$ and $j$ exchanged, that $\x_{e,j} =0$ for all
$j<i$. This shows that all the coordinates $\x_{e,j}$ are zero. But the
point $0_e$ does not lie on $\P_e$, a contradiction. The same argument works for when $\x_{\ol e,i}=0$, and for 
the scheme $\mathfrak X_e^{(n)}(a)$ for every
$n$. 

Property (3) follows easily, by observing that setting $t:=0$ in
Equations ($\star$) (resp.~($\star_n$)) 
yields the equations defining $\mathbf R_e$ (resp.~$\mathbf
R_e^{(n)}$).

As for flatness and Property (1), it is enough to prove them at a point on
the special fiber. For a point on the special fiber $\mathfrak X_e(a)$, equal to $\mathbf
R_e$ by Property~(2), this means a point
with coordinates $\x_{e,m}=0$ for all $m<\ell$ and $\x_{\ol e,m}=0$ for
all $m>\ell$, for some $\ell$. Furthermore, we may assume that
$\x_{\ol e,\ell}\neq
0$. The parameters at the point are thus the ratios $\x_{e,m}/\x_{\ol e,m}$
for $m\leq\ell$, the ratios $\x_{\ol e,m}/\x_{e,m}$ for $m>\ell$ and
$t$. Equations ($\star$) become the following equations on the parameters:
\begin{align*}
&\frac{\x_{e,m}}{\x_{\ol e,m}}=
(at)^{\ell-m}\frac{\x_{e,\ell}}{\x_{\ol e,\ell}}
\quad\text{for $m<\ell$,}\\
&\frac{\x_{\ol e,m}}{\x_{e,m}}=
(at)^{m-\ell-1}\frac{\x_{\ol e,\ell+1}}{\x_{e,\ell+1}}
\quad\text{for $m>\ell+1$ and}\\
&\frac{\x_{e,\ell}}{\x_{\ol e,\ell}}\frac{\x_{\ol e,\ell+1}}{\x_{e,\ell+1}}=at.
\end{align*}
We obtain that $\x_{e,\ell}/\x_{\ol e,\ell}$
and $\x_{\ol e,\ell+1}/\x_{e,\ell+1}$ are free parameters, the only
ones, and thus $\mathfrak X_e$ is regular. Furthermore, the above last
equation proves that $\mathfrak X_e$ is flat over $B$. The same argument works for
$\mathfrak X_e^{(n)}(a)$, or by simply observing that $\mathfrak
X_e^{(n)}(a)$ is the open subscheme of $\mathfrak X_e(a)$, given by 
$\x_{\ol e,-n}\neq 0$ and $\x_{e,n}\neq 0$. Indeed, these inequalities
impose no conditions on the generic fiber of $\mathfrak X_e(a)$ over
$B$, whereas they extract the special fiber of $\mathfrak
X_e^{(n)}(a)$ over $B$ from that of  $\mathfrak X_e(a)$.
\end{proof}

Observe that the families $\mathfrak X_e(a)$ and the $\mathfrak
X_e^{(n)}(a)$ over $B$ are $\Gm$-equivariant, in the sense that the
natural action of $\Gm$ on $\P_e$, taking a point
with coordinates $(x_{e,i}:x_{\ol e,i})$ to that with coordinates
$(cx_{e,i}:x_{\ol e,i})$ for each $c\in\Gm(k)$, induces naturally
actions on $\mathfrak X_e(a)$ and the $\mathfrak
X_e^{(n)}(a)$ leaving invariant their fibers over $B$. Also, under the
identification of the general fiber of any of these families with
$\Gm$ stated in the proof of Proposition~\ref{prop:gm}, the action on the general
fiber is that given by the group structure of $\Gm$.

Given in addition $\ell\in\mathbb N$, define $\mathfrak X_e(a,\ell)$
(resp.~$\mathfrak X_e^{(n)}(a,\ell)$) as the base change of $\mathfrak
X_e(a)$ (resp.~$\mathfrak X_e^{(n)}(a)$) by the map $B\to B$
given by $t\mapsto t^\ell$. Of course, also $\mathfrak X_e(a,\ell)$
and the $\mathfrak X_e^{(n)}(a,\ell)$ are $\Gm$-equivariant families
over $B$.
 
\subsection{Degeneration of $\mathbf G_{\mathbf m}^{|V|-1}$ to $Y^{a,b}_{\ell,\m}$}
Let $G=(V,E)$ be a connected graph without self loops, and $\ell: E \to \mathbb N$ a length
function. Let 
$a\colon C^1(G, \mathbb Z) \rightarrow \k^*$ and $b\colon H^1(G, \mathbb Z)
\to \k^*$ be characters, and $\m\in C^1(G, \mathbb Z)$. Fix an orientation
$\mathfrak o$ of the edges of the graph. Recall the notation:
$a_e:=a(\chi_{\indm e}-\chi_{\indmbar{\ol e}})$ for each $e\in\E$.

Recall the subscheme 
$Y^{a,b}_{\ell,\m}$ of $\mathbf R$ that we introduced in
Section~\ref{mixtortil}, and the action of  $\mathbf G_{\mathbf
  m}^{|V|-1}$ on $\mathbf R$ that leaves $Y^{a,b}_{\ell,\m}$
invariant in Section~\ref{mixorb}. 
The aim of this section is to describe a $\mathbf G_{\mathbf
  m}^{|V|-1}$-equivariant degeneration 
of $\mathbf G_{\mathbf m}^{|V|-1}$ to $Y^{a,b}_{\ell,\m}$. 

More
precisely, let $B:= \mathrm{Spec}(\k[[t]])$. Let 
$$
\mathbf W := \prod_{e\in E^\mathfrak o}\P_e \subset \prod_{e\in E^\mathfrak o} \prod_{i\in \mathbb Z} \P^1_{e,i}.
$$
And let $\mathfrak Y^{a,b}_{\ell,\m}$ be defined in the product
$\mathbf W\times B$ by the following set of equations:
\begin{itemize}
\item[($\star_{e,i,j}$)] For each edge $e \in E^{\mathfrak o}$ and each pair $(i,j)$ of
  integers with $i<j$:
$$\x_{e,i} \x_{\ol e, j} = (a_et^{\ell_e})^{j-i} \x_{\ol e,i} \x_{e,j}.$$
\item[($\star_{\gamma,\alpha}$)] 
For each $\gamma\in H^1(G,\mathbb Z)$ and each $\alpha \in C^1(G,\mathbb Z)$:
\[
\Bigl(b(\gamma) \prod_{e\in E^{\mathfrak o}} a_e^{\gamma_e\alpha_e}
\Bigr)\cdot 
\prod_{e\in E^{\mathfrak o}} \Bigl(\frac {\x_{\ol e,\alpha_e}}{\x_{e,\alpha_e}}\Bigr)^{\gamma_e} 
= 
t^{\sum_{e\in E^{\mathfrak o}}\gamma_e \alpha_e \ell_e - \m_e\gamma_e}.
\]
\end{itemize}
While Equations ($\star_{e,i,j}$) are clear, by Equations
($\star_{\gamma,\alpha}$) we rather mean:
\begin{itemize}
\item[($\star'_{\gamma,\alpha}$)] For each $\gamma\in H^1(G,\mathbb Z)$ and each $\alpha \in C^1(G,\mathbb Z)$:
\end{itemize}
\medskip
\begin{itemize}
\item If $\sum_{e\in E^\mathfrak o} \gamma_e\alpha_e \ell_e -\m_e\gamma_e \geq 0$ then
\[
\Bigl(b(\gamma) \prod_{e\in E^{\mathfrak o}} a_e^{\gamma_e\alpha_e}
\Bigr)\cdot
\prod_{\substack{e\in E^{\mathfrak o}\\
\gamma_e> 0}} \x_{\ol e,\alpha_e}^{\gamma_e}\cdot\prod_{\substack{e\in E^{\mathfrak o}\\
\gamma_e<0}}\x_{e,\alpha_e}^{-\gamma_e} = t^{\sum_{e\in E^{\mathfrak o}}\gamma_e \alpha_e \ell_e - \m_e\gamma_e} \prod_{\substack{e\in E^{\mathfrak o}\\
\gamma_e> 0}} \x_{e,\alpha_e}^{\gamma_e}\cdot\prod_{\substack{e\in E^{\mathfrak o}\\
\gamma_e<0}}\x_{\ol e,\alpha_e}^{-\gamma_e}. 
\] 
\medskip

\item If $\sum_{e\in E^\mathfrak o} \gamma_e\alpha_e \ell_e -\m_e\gamma_e < 0$ then
\[
\Bigl(b(\gamma) \prod_{e\in E^{\mathfrak o}} a_e^{\gamma_e\alpha_e}
\Bigr)\cdot 
t^{\sum_{e\in E^{\mathfrak o}}\m_e\gamma_e - \gamma_e \alpha_e \ell_e}
\cdot\prod_{\substack{e\in E^{\mathfrak o}\\
\gamma_e> 0}} \x_{\ol e,\alpha_e}^{\gamma_e}\cdot\prod_{\substack{e\in E^{\mathfrak o}\\
\gamma_e<0}}\x_{e,\alpha_e}^{-\gamma_e} = \prod_{\substack{e\in E^{\mathfrak o}\\
\gamma_e>  0}} \x_{e,\alpha_e}^{\gamma_e}\cdot\prod_{\substack{e\in E^{\mathfrak o}\\
\gamma_e<0}}\x_{\ol e,\alpha_e}^{-\gamma_e}. 
\] 
\end{itemize}

For each $v\in V$ (resp.~$e\in\E$), let $\mathbf G_{\mathbf m,v}$ (resp.~$\mathbf G_{\mathbf m,e}$) be a copy of
$\Gm$. In Subsection~\ref{degG} we defined an action of $\mathbf
G_{\mathbf m,e}$ on $\P_e$ for each $e\in\E$. It induces an action of
$\mathbf G_{\mathbf m}^{E}:=\prod_{e\in E^{\mathfrak o}}\mathbf G_{\mathbf m,e}$ on
$\mathbf W$ and on $\mathbf W\times B$. 

There is a natural map of group
schemes
\begin{equation}\label{GVE}
\mathbf G_{\mathbf m}^{V}\to \mathbf G_{\mathbf m}^{E},
\end{equation}
where $\mathbf G_{\mathbf m}^{V}=\prod_{v\in V}\mathbf G_{\mathbf
  m,v}$, defined by sending $(c_v\,|\,v\in V)$ to $(\mu_e\,|\,e\in
E^{\mathfrak o})$, where $\mu_e:=c_v/c_u$ for each $e=uv\in E^{\mathfrak
  o}$. Since $G$ is connected, the kernel of the homomorphism is
$\Gm$, embedded diagonally in $\mathbf G_{\mathbf m}^{V}$. It thus
follows that $\mathbf G_{\mathbf m}^{V}/\Gm$ acts naturally on
$\mathbf W$ and on $\mathbf W\times B$. Furthermore, it leaves
$\mathfrak Y_{\ell, \m}^{a,b}$ invariant.

\begin{thm}\label{thm:degtori} $\mathfrak Y_{\ell, \m}^{a,b}$ is
  naturally a subscheme of the fibered product $\prod_{e\in E^\mathfrak o} \mathfrak X_e(a_e,\ell_e)$ over $B$. It has generic fiber over $B$  isomorphic to $\mathbf G_{\mathbf m} ^V/\Gm$, and special fiber 
equal to $Y^{a,b}_{\ell,\m}$. Furthermore, the action of $\mathbf
G_{\mathbf m}^{V}/\Gm$ on $\mathfrak Y_{\ell, \m}^{a,b}$ is identified
with the action given by the group structure on the generic fiber and
the action defined in Subsection~\ref{action} on the special fiber.
\end{thm}

\begin{proof} Equations ($\star_{e,i,j}$) are those of the fibered product
  $\prod_{e\in E^\mathfrak o} \mathfrak X_e(a_e,\ell_e)$ over $B$, so the
  first statement is clearly true. As for the second statement,
  observe that setting $t=0$ in Equations~($\star_{e,i,j}$) and
  Equations ($\star'_{\gamma,\alpha}$) yield the equations for
  $Y_{\ell,\m}^{a,b}$ in $\mathbf W$, as we proved in
  Theorem~\ref{thm:main2}. This shows that the special fiber of
  $\mathfrak Y_{\ell,\m}^{a,b}/B$ is indeed $Y_{\ell,\m}^{a,b}$. Also,
  that the action of $\mathbf
G_{\mathbf m}^{V}/\Gm$ on the special fiber coincides with that defined in
Subsection~\ref{action} on $Y_{\ell,\m}^{a,b}$ is clear, by comparing
the two definitions.

As for the generic fiber, when $t\neq 0$, Equations ($\star_{e,i,j}$)
imply that all coordinates $\x_{e,i}$ and $\x_{\ol e,i}$ are nonzero,
and the ratios $\x_{e,i}/\x_{\ol e,i}$ are determined by the ratios
$\x_{e,0}/\x_{\ol e,0}$ for $e\in E^\o$, as we have seen in the proof
of Proposition~\ref{prop:gm}. Equations ($\star'_{\gamma,\alpha}$) are
thus equivalent to Equations ($\star_{\gamma,\alpha}$). And those
equations follow from the Equations ($\star_{\gamma,0}$), given
the expressions for the ratios $\x_{e,i}/\x_{\ol e,i}$ in terms of the
ratios $\x_{e,0}/\x_{\ol e,0}$ obtained from Equations
($\star_{e,i,j}$). 

The ratios $\x_{e,0}/\x_{\ol e,0}$ describe the generic fiber of
$\mathfrak Y_{\ell,\m}^{a,b}$ as the subscheme of $\mathbf{G}_{\mathbf m}^E\subset
\prod_{e\in E^{\mathfrak o}} \P^1_{e,0}$ given by the equations:
\[
\forall \,\, \gamma \in H^1(G, \mathbb Z), \qquad 
b(\gamma)\cdot \prod_{e\in E^{\mathfrak o}} 
\Bigl( \frac{\x_{\ol e,0}}{\x_{e,0}}\Bigr)^{\gamma_e}
=t^{-\sum_{e\in E^\o} \m_e\gamma_e}.
\]
Consider the map of group schemes
\[
\mathbf{G}_{\mathbf m}^E \longrightarrow \mathrm{Hom}(H^1(G, \mathbb
Z), \Gm),
\]
which on the level of points is given by sending $\underline
\mu=(\mu_e\,|\,e\in E^{\mathfrak o})$ of $\mathbf{G}_\mathbf m^E$ 
to the morphism $\phi_{\underline \mu}$ defined by 
$$
\phi_{\underline \mu} \bigl(\gamma\bigr) := \prod_{e\in E^{\mathfrak o}} 
\mu_e^{\gamma_e} \qquad \textrm{for each }\gamma \in H^1(G,\mathbb Z).
$$  
Its kernel is naturally identified with the image of the map of group
schemes $\mathbf G_{\mathbf m}^V \to\mathbf{G}_{\mathbf m}^E$ given in \eqref{GVE}.
The kernel of the latter is isomorphic to $\Gm$ via the diagonal embedding.

Thus the natural action (coordinate by coordinate) of $\mathbf
G_{\mathbf m}^E$ on $\mathbf G_{\mathbf m}^E$ induces an action of the
group scheme $\mathbf G_{\mathbf m}^V/\Gm$ on $\mathbf G_\mathbf
m^E$ which is transitive on the subscheme defined by the equations
$$
\forall \,\, \gamma \in H^1(G, \mathbb Z), \qquad \prod_{e\in
  E^{\mathfrak o}} \mu_e^{\gamma_e}=b'(\gamma)
$$
for every character $b'\colon H^1(G, \mathbb Z)\to \k^*$. It
follows that the generic fiber of 
$\mathfrak Y_{\ell,\m}^{a,b}$ is a split torsor over $\mathbf G_\mathbf m^{V} /\mathbf G_m$. 
\end{proof}

By substituting $\prod_{e\in E^{\mathfrak o}}\mathbf P^{(n_e)}_e$ for
$\mathbf W$, for any given $n\: E\to\Z_{\geq 0}$, and restricting the range of the $i, j$ in Equations ($\star_{e,i,j}$)
and the $\alpha_e$ in Equations  ($\star'_{\gamma,\alpha}$) to the
interval $[-n_e,n_e]$ for each $e\in E^{\mathfrak o}$, we obtain a degeneration of $\mathbf G_\mathbf
m^{|V|-1}$ to a ``finite approximation'' of $Y^{a,b}_{\ell,\m}$,
namely to $Y^{a,b}_{\ell,\m}\cap\mathbf R^{(n)}$, where
$\mathbf R^{(n)}:=\prod_{e\in E^\o}\mathbf R^{(n_e)}_e$.

 \vspace{.7cm}

\subsection*{Acknowledgements.}  This project benefited very much from
the hospitality of the Mathematics Department at the \'Ecole Normale
Sup\'erieure (ENS) in Paris and the Instituto de Matem\'atica Pura e
Aplicada (IMPA) in Rio de Janeiro during mutual visits of both authors.
We thank the two institutes and their members for providing for those visits.
We are also specially grateful to the Brazilian-French Network in
Mathematics for providing support for a visit of E.E.~to ENS Paris
and a visit of O.A.~to IMPA. 
 
 \bibliographystyle{alpha}
\bibliography{bibliography}
\end{document}